\documentclass[twoside]{book}
\usepackage{book}
\usepackage{notations}
\usepackage{change}
    
\hypersetup{
pdftitle={Геометрия группы симплектоморфизмов},
pdfauthor={Л. В. Полтерович (Перевод с английского Р. Г. Матвеева и А. М. Петрунина
pод редакцией Л. В. Полтеровича и Я. М. Элиашберга)},
pdfsubject={Дифференциальная геометрия}}

\renewcommand\?[2]{#1}
\begin{document}

\title{Геометрия группы симплектоморфизмов}
\author{Л. В. Полтерович
\\
\\
Перевод с английского
\\
Р. Г. Матвеева и А. М. Петрунина.
\\
\\
Под редакцией\\
Л. В. Полтеровича и Я. М. Элиашберга.
}
\date{}
\maketitle

\thispagestyle{empty}
\noindent\textbf{Л. В. Полтерович}

Геометрия группы симплектоморфизмов/
пер. с англ. Р. Г. Матвеева и А. М. Петрунина;
под ред. Л. В. Полтеровича и Я. М. Элиашберга.

\

Предварительное издание предназначенное исключительно для отлова ляпов. 
Исправления слать по адресу 
\url{petrunin@math.psu.edu}
или
\url{matveev@mis.mpg.de}.
\null
\vfill
Translation from the English language edition:
The Geometry of the Group of Symplectic Diffeomorphisms by Leonid Polterovich
\\
Copyright \textcopyright{} Birkhäuser Verlag AG, 2001.
All Rights Reserved.
\vfill
\noindent{\includegraphics[scale=0.5]{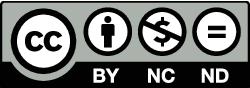}
\vspace*{1mm}
\\
\hbox{\parbox{1\textwidth}
{Лицензия: CC BY-NC-ND 4.0,\\
\texttt{https://creativecommons.org/licenses/by-nc-nd/4.0/}}}

\tableofcontents
\chapter*{Предисловие}

Группа гамильтоновых диффеоморфизмов $\Ham(M,\Omega)$ симплектического
многообразия $(M,\Omega)$ играет
основополагающую роль в геометрии и
механике.
Для геометра, по крайней мере при некоторых предположениях
о многообразии $M$, это связная компонента тождественного отображения
в группе всех симплектических диффеоморфизмов.
С точки зрения
механики $\Ham(M,\Omega)$ является группой всех допустимых движений.
Сколько нужно энергии для реализации данного гамильтонова диффеоморфизма $f$?
Попытка формализовать этот естественный вопрос и ответить на него привела \rindex{Хофер}Х. Хофера \cite{H1} (1990) к удивительному открытию.
Оказалось, что решение этой вариационной задачи можно интерпретировать как {}\emph{геометрическую величину}, а именно как расстояние между $f$ и тождественным
отображением.
Это расстояние связано с канонической биинвариантной метрикой на $\Ham(M,\Omega)$, которую стали называть хоферовской метрикой.
Начиная с работ Хофера, эта новая геометрия интенсивно изучалась в рамках современной
симплектической топологии.
В настоящей книге я опишу некоторые из полученых результатов.

Хоферовская геометрия позволяет изучать различные понятия и задачи из знакомой нам конечномерной геометрии в контексте группы гамильтоновых диффеоморфизмов.
При этом они сильно отличаются от обычного круга задач, рассматриваемых в симплектической топологии, и, таким образом, значительно расширяют горизонты симплектического мира.
Бесконечен ли диаметр $\Ham(M,\Omega)$?
Что там за кратчайшие?
Как найти спектр длин?
В общем случае эти вопросы остаются открытыми.
В некоторых частных случаях ответы найдены и будут нами рассмотрены.

Есть ещё одна, на мой взгляд, даже более важная причина, почему полезно иметь каноническую геометрию на группе гамильтоновых диффеоморфизмов.
Рассмотрим зависящее от времени векторное поле $\xi_t$, $t\in \RR$, на многообразии $M$.
Обыкновенное дифференциальное уравнение
\[\dot x=\xi(x,t)\]
на $M$ определяет поток $f_t\: M \to M$, который отображает точку $x(0)$ в $x(t)$ — значение решения в момент времени $t$.
Траектории потока образуют сложную систему кривых на многообразии.
Обычно, чтобы разобраться в динамике, нужно следить за этими кривыми и изучать их поведение в разных регионах многообразия.
Сменив точку зрения, мы видим, что наш поток становится простым
геометрическим объектом — одной кривой $t \mapsto f_t$ на группе всех диффеоморфизмов многообразия.
Хочется надеяться, что геометрические свойства этой кривой отражают динамику, и в таком случае сложную динамическую систему можно изучать чисто геометрическими средствами. 
Разумеется, за это придётся платить.
Дело в том, что объемлющее пространство — группа
диффеоморфизмов — бесконечномерно.
Более того, появляется другая трудность:
в общем случае у нас нет канонических инструментов для выполнения геометрических измерений на этой группе.
Примечательно, что хоферовская метрика даёт такой инструмент для систем классической механики.
В главах 8 и 11 мы увидим, что такой ход рассуждений оказывается полезным.

Как часто происходит с молодыми быстроразвивающимися областями математики, доказательства некоторых красивых утверждений хоферовской геометрии оказываются технически сложными.
Поэтому я выбрал самые простые нетривиальные случаи основных утверждений (на свой вкус, конечно), не пытаясь представить их в максимальной общности.
По той же причине опущены многие технические детали.
Хотя формально эта книга не требует особых знаний в симплектической топологии (по крайней мере, необходимые определения и формулировки приведены), я советую  два замечательных вводных текста \cite{HZ} и \cite{MS}.
Оба содержат главы по геометрии группы гамильтоновых диффеоморфизмов.
Я постарался избежать повторов.
Книга содержит упражнения, которые предположительно помогут читателю вникнуть в предмет.

Эта книга возникла из двух источников.
Первый — это лекции для аспирантов, а именно миникурсы в университетах Фрайбурга и Уорика, а также курс лекций в Швейцарской высшей технической школе Цюриха.
Вторым источником является обзорная статья \cite{P8}, которая вкратце содержит материал книги.

Теперь я коротко опишу главных персонажей книги.
Диффеоморфизм $f$ симплектического многообразия $(M,\Omega)$
называется \emph{гамильтоновым}, если его можно включить в гамильтонов поток $f_t$ с компактным носителем, удовлетворяющий $f_0=\1$, и $f_1 =f$.
Такой поток определяется гамильтонианом $F\: M \times [0;1] \to \RR$.
На языке классической механики $F$ — энергия механического движения, описывающая $f_t$.
Мы понимаем полную энергию потока как длину соответствующего пути диффеоморфизмов:
\[\length \{f_t\} =
\int_0^1\big(\max_{x\in M}F(x,t)-\min_{x\in M}F(x,t)\big)\, dt 
\]
Определим функцию
\[\rho\: \Ham(M,\Omega) \times \Ham(M,\Omega) \to \RR\]
как
\[\rho (\phi, \psi) = \inf \length \{f_t\},\]
где точная нижняя грань берётся по всем гамильтоновым потокам $\{f_t\}$, которые
порождают гамильтонов диффеоморфизм $f = \phi\psi^{-1}$.
Легко видеть, что $\rho$ неотрицательна, симметрична, обращается в
ноль на диагонали и удовлетворяет неравенству треугольника. 
Более того, $\rho$ биинвариантна по отношению к групповой структуре на $\Ham(M,\Omega)$.
Иными словами, $\rho$ — биинвариантная псевдометрика.
Глубокая теорема говорит, что $\rho$ — настоящая метрика, то есть
$\rho (\phi, \psi)$ строго положительна при $\phi \ne \psi$. 
Метрика $\rho$ и называется хоферовской метрикой.

Группа $\Ham(M,\Omega)$ и хоферовская псевдометрика $\rho$ обсужадются в главах 1 и 2 соответственно.
В главе 3 доказывается, что $\rho$ является настоящей метрикой в случае, когда $M$ --- стандартное симплектическое линейное пространство $\RR^{2n}$.
Наш подход основан на теории Громова голоморфных дисков с лагранжевыми граничными условиями, см.
главу~4.

Затем мы переходим к изучению основных геометрических инвариантов $\Ham(M,\Omega)$.
Существует (пока ещё открытая!) гипотеза о том, что диаметр $\Ham(M,\Omega)$ бесконечен.
В главах 5—7 эта гипотеза доказывается для замкнутых поверхностей.

В главе 8 обсуждается рост однопараметрической подгруппы $\{f_t\}$ группы $\Ham(M,\Omega)$, которая описывает асимптотическое поведение функции $\rho (\1, f_t)$ при $t \to \infty$.
Мы описываем связь между ростом и динамикой $\{f_t\}$ в контексте
теориии инвариантных торов классической механики.

Во многих интересных случаях пространство $\Ham(M,\Omega)$ имеет сложную топологию и, в частности, нетривиальную фундаментальную группу.
Для $\gamma\in\pi_1(\Ham(M,\Omega))$ положим $\nu(\gamma) \z= \inf \length \{f_t\}$, где
точная нижняя грань берётся по всем петлям $\{f_t\}$ гамильтоновых
диффеоморфизмов (то есть по периодическим гамильтоновым потокам),
которые представляют класс петель $\gamma$.
Множество
\[\set{\nu(\gamma)}{\gamma\in\pi_{1}\Ham(M,\Omega)}\]
называется спектром длин группы $\Ham(M,\Omega)$.
В главе 9 мы представляем подход к оценке спектра длин, основанный на
теории симплектических расслоений. 
Важным ингредиентом этого подхода является теория Громова
псевдоголоморфных кривых, обсуждаемая в главе~10. 
В главе 11 приводятся приложения наших результатов к спектру длин в
классической эргодической теории.

В главах 12 и 13 развиваются два разных подхода к теории геодезических на $\Ham(M,\Omega)$.
Один из них элементарный, а другой требует мощного инструмента — гомологий Флоера.
Краткое введение в эту теорию дано в главе 13.

Наконец, в главе 14 мы обсуждаем негамильтоновы симплектические
диффеоморфизмы, они естественно появляются как изометрии хоферовской
геометрии $\Ham(M,\Omega)$. 
Кроме того, мы формулируем и обсуждаем знаменитую \rindex{гипотеза
  потока}\emph{гипотезу потока},
заключающуюся в том, что группа $\Ham(M,\Omega)$ замкнута в группе
всех симплектических диффеоморфизмов с $C^\infty$-топологией.

\paragraph*{Благодарности.}
Я сердечно благодарен \rindex{Аквельд}Мейке Аквельд %
за её незаменимую помощь в напечатании предварительной рукописи,
подготовку рисунков и огромную редакционную работу.
Я очень благодарен Паулю Бирану и Карлу Фридриху Зибургу за их подробные комментарии к рукописи и за улучшения в изложении.
Я признателен 
Рами Айзенбуду,
Диме Гуревичу, 
Мише Энтову, 
Осе Полтеровичу
и Зеэву Руднику за указание на ряд неточностей в предварительной версии книги.
Книга написана во время моего пребывания в Швейцарской высшей технической школе Цюриха в 1997—1998 учебном году и во время моих визитов в Институт высших научных исследований в Бюр-сюр-Иветте в 1998 и 1999 годах.
Я благодарю оба этих института за прекрасную атмосферу для работы. 

\parbf{Добавления к русскому переводу.}
За последние пару десятилетий хоферовская геометрия группы гамильтоновых диффеоморфизмов обогатилась рядом новых интересных достижений.
Для удобства читателя я добавил к переводу комментарии на несколько устаревших замечаний и решенных к настоящему моменту (январь, 2022) проблем и гипотез.
Кроме того, я включил ссылки на новые публикации по теме книги.

\chapter{Знакомство с группой}\label{chap:1}

В этой главе мы обсудим некоторые хорошо известные утверждения о группе гамильтоновых диффеоморфизмов. 

\section[Гамильтоновы диффеоморфизмы]{Гамильтоновы диффеоморфизмы}

Рассмотрим движение частицы единичной массы в $\RR^n(q)$ под
действием потенциальной силы $\Phi(q, t)  \z= - \tfrac{\partial
  U}{\partial q} (q, t)$, здесь $q$ обозначает
координату\footnote{Здесь и далее  $q$ это сокращённая запись для
  $q_1,\dots,q_n$. — \textit{Прим. ред.}} в $\RR^n$. 
Согласно второму закону Ньютона $\ddot q= \Phi (q, t)$.
За исключением некоторых редких случаев, это уравнение явно решить не получается.
Однако можно понять некоторые качественные свойства решений.

Применим небольшую хитрость.
Введём вспомогательную переменную $p = \dot q$ и рассмотрим функцию
$F(p,q,t)= \tfrac {p^2} 2 + U (q, t)$. 
Функция $F$ представляет собой полную энергию частицы (сумму её
кинетической и потенциальной энергий). 
В этих обозначениях приведённое выше уравнение Ньютона можно переписать следующим образом:
\[
\begin{cases}
\quad\dot p &= - \tfrac{\partial F}{\partial q} (p, q, t),\\
\quad\dot q &= \tfrac{\partial F}{\partial p} (p, q, t).
\end{cases}
\]
Эта система дифференциальных уравнений первого порядка называется
\rindex{гамильтонова система}\emph{гамильтоновой системой}. 
Её д\'{о}лжно рассматривать как обыкновенное дифференциальное уравнение в
$2n$-мерном пространстве $\RR^{2n}$ с координатами $p$ и $q$. 
Первый шаг любого качественного исследования состоит в том, чтобы
отказаться думать о явной форме интересующего нас объекта. 
Приняв это к сведению, давайте оставим в стороне выражение для функции $F$, приведённое выше, и сосредоточимся на общем случае гамильтоновых систем, связанных с
более-менее произвольными гладкими функциями энергии $F (p, q, t)$. 
«Более-менее» означает, что мы накладываем определенные ограничения на
поведение $F$ на бесконечности, гарантирующие, что решения
гамильтоновой системы существуют при всех $t\in \RR$. 
Итак, выберем такое $F$ и рассмотрим поток $f_t\: \RR^{2n} \to
\RR^{2n}$, переводящий любое начальное условие $\big(p(0),q(0)\big)$ в
значение $\big(p (t), q (t)\big)$ соответствующего решения в момент
времени $t$. 
Возникающие таким образом диффеоморфизмы $f_t$ будут неформально называться \rindex{механическое движение}\emph{механическими движениями}.
Преимущество такого подхода в том, что эти диффеоморфизмы $2n$-мерного пространства $\RR^{2n}$ обладают следующими замечательными геометрическими свойствами, которые нельзя увидеть в исходном конфигурационном пространстве $\RR^n$.

\begin{thm}[(Теорема Лиувилля)]{Теорема}\label{1.1.A}\rindex{теорема Лиувилля}
Механические движения сохраняют форму объёма 
$\Vol=dp_1\wedge dq_1\wedge\dots\wedge d p_n\wedge dq_n$.
\end{thm}

\begin{thm}[(Более тонкий вариант \ref{1.1.A})]{Теорема}\label{1.1.B}
Механические движения сохраняют 2-форму $\omega = dp_1 \wedge dq_1 +\dots
+ dp_n \wedge dq_n.$
\end{thm}

Обратите внимание, что $\Vol = \tfrac{\omega^n}{n!}$ и, значит, \ref{1.1.B} влечёт \ref{1.1.A}.
Далее, заметим, что при $n = 1$ теоремы \ref{1.1.A} и \ref{1.1.B} равносильны.
Теорема \ref{1.1.B} является простым следствием того, что механические движения происходят из гамильтоновой системы.
Мы приведём её доказательство в разделе \ref{sec:1.3}.

Оба приведённых выше результата получены давно.
Сохранение объёма механическими движениями привлекало большое внимание уже более века назад.
Это свойство послужило основой при создании {}\emph{эргодической теории}, ныне хорошо известной математической дисциплины, изучающей различные свойства повторяемости преобразований, сохраняющих меру.
Однако значение роли инвариантной 2-формы $\omega$ заметили сравнительно недавно.
Насколько я знаю, впервые прямо указал на это \rindex{Арнольд}В. И. Арнольд в 1960-х годах.
Попытки понять разницу между механическими движениями и диффеоморфизмами, сохраняющими объём, породила {}\emph{симплектическую топологию}, которая исследует неожиданные явления жёсткости, возникающие в теории симплектических многообразий и их морфизмов.

Вот пример такого явления, которое обнаружил
\rindex{Сикорав}Ж.-К. Сикорав \cite{S1}%
\footnote{\label{foot:sikorav}В статье Сикорава сформулирована и доказана
  теорема~\ref{3.2.B}, из которой теорема~\ref{1.1.C} сразу
  следует. При этом результат приписывается автором М. Л. Громову, как
  сообщённый в частной переписке. -- \textit{Прим. Ред.}}.
  
Пусть $B^2(r) \subset \RR^2$ — евклидов диск радиуса $r$, ограниченный окружностью $S^1(r)$.
Рассмотрим тор 
\[L_R = 
S^1 (R) \times\dots\times S^1(R) \subset \RR^2 (p_1, q_1) \times\dots\times \RR^2 (p_n, q_n) =\RR^{2n} (p, q)\]
и цилиндр $C_r = B^2 (r) \times \RR^{2n-2}$.

\begin{thm}[(Свойство несжимаемости)]{Теорема}\label{1.1.C}
Не существует механического движения, переводящего $L_R$ в $C_r$, при условии, что $R\z>r$.
\end{thm}

Мы докажем это утверждение в более общем виде в пункте~\ref{3.2.E}.
Отметим, что при $n = 1$ результат очевиден.
Действительно, площадь, ограниченная $S^1 (R)$, больше, чем площадь $B^2 (r)$.
Таким образом, нельзя перевести $S^1(R)$ в $B^2(r)$ преобразованием, сохраняющим площадь.
Однако если $n \ge 2$, то $L_R$ — подмногообразие коразмерности $n \ge 2$, а $C_r$ имеет бесконечный объем.
Таким образом, нет видимой причины, по которой это утверждение должно быть верным.
Более того, подобное утверждение вовсе неверно в категории диффеоморфизмов, сохраняющих объём!

\begin{ex*}{Упражнение}
Найти линейное преобразование $\RR^{2n} \to \RR^{2n}$, сохраняющее объём и переводящее $L_R$ в $C_r$ при произвольных положительных $r$ и $R$.
\end{ex*}

Далее эволюция механической системы будет рассматриваться как кривая в группе всех механических движений и эта кривая будет изучаться геометрическими средствами.
Чтобы всё работало, нам придётся сузить класс механических систем.
Например, неограниченные гамильтонианы, такие как $F(p,q,t) = \tfrac {p^2}2 + U(q,t)$ рассматриваемые выше, для нас будут слишком сложны.
Мы всегда будем предполагать, что гамильтонианы (и, следовательно, соответствующие им механические движения) имеют компактный носитель.
Иными словами, всё движение происходит в ограниченной части нашего пространства.

В этой главе вводится гамильтонова механика на симплектических многообразиях, которая является естественным обобщением модели, описанной выше.
Соответственно, гамильтонов диффеоморфизм — это просто механическое движение, порождённое гамильтонианом с компактным носителем. 

\section{Потоки и пути диффеоморфизмов}

Для начала объясним связь между потоками и дифференциальными уравнениями и дадим геометрическую интерпретацию потоков как путей диффеоморфизмов.
Имея в виду дальнейшие приложениями мы предполагаем, что все рассматриваемые нами объекты имеют компактный носитель.
Однако основные построения, описанные ниже, проходят в более общем случае.

Рассмотрим гладкое многообразие $M$ без края.
Пусть $\phi\: M \to M$ — диффеоморфизм.
Определим его \rindex{носитель}\emph{носитель} \index[symb]{$\supp(\phi)$}$\supp(\phi)$ как замыкание множества всех $x \in M$, что $\phi(x) \ne x$.
Обозначим через $\Diff^c (M)$ группу всех диффеоморфизмов с компактным носителем.
Пусть $I \subset \RR$ — интервал.%
\footnote{Интервал определяется как связное подмножество $\RR$ с непустой внутренней частью.}
\rindex{путь диффеоморфизмов}\emph{Путь в группе диффеоморфизмов} — это отображение 
\[f\: I \to \Diff^c (M),\quad t \mapsto f_t\]
со следующими свойствами:
\begin{itemize}
\item отображение $M \times I \to M$, переводящее $(x, t)$ в $f_t x$, гладкое;
\item существует компактное подмножество $K$ в $M$, содержащее $\supp f_t$ при всех $t \in I$.
\end{itemize}
Такой путь будет обычно обозначаться через $\{f_t\}$.
Отметим, что на замкнутых многообразиях второе условие выполняется автоматически.

Каждый путь диффеоморфизмов порождает семейство векторных полей $\xi_t$, $t \in I$ на $M$ определяемых следующим образом: 
\begin{equation}\tfrac{d}{dt} f_t x = \xi_t (f_t x).
\label{eq:1.2.A}
\end{equation}
Отметим, что это гладкое семейство с компактным носителем: $\xi_t (x) \z= 0$ при всех $x \in M \setminus K$.
Такое семейство будет называться \textit{зависящим от времени векторным полем с компактным носителем} на~$M$.
Приведённое выше соответствие не является инъективным.
В самом деле, если $g$ — произвольный элемент $\Diff^c (M)$, то путь вида $\{f_t g\}$ порождает то же самое зависящее от времени векторное поле $\xi$.
Однако для каждой точки $s \in I$ существует единственный путь $\{f_t\}$, который 
порождает $\xi$ такой, что $f_s$ равно тождественному отображению $\1$.
Этот путь определяется как единственное решение уравнения~(\ref{eq:1.2.A}), которое теперь рассматривается как обыкновенное дифференциальное уравнение с начальным условием $f_s = \1$.
Предположим, что $0 \in I$, и возьмём $s = 0$.
Построенный выше путь $\{f_t\}$ при $f_0 = \1$ называется \rindex{поток}\emph{потоком} зависящего от времени векторного поля $\xi$.
Таким образом, потоки — это просто пути $\{f_t\}$, такие, что $f_0 = \1$.

\section[Классическая механика]{Математическая модель классической механики}\label{sec:1.3}

Роль фазового пространства в классической механике играет \rindex{симплектическое многообразие}\emph{симплектическое многообразие} $(M^{2n},\Omega)$.
Мы считаем, что многообразие $M$ связно, без границы и чётной размерности $2n$, а $\Omega$ — замкнутая дифференциальная 2-форма на $M$.
Форма $\Omega$ предпологается {}\emph{невырожденной}.
То есть, её максимальная степень $\Omega^n$ не обращается в нуль ни в какой точке.
Форма $\Vol =  \tfrac{\Omega^n}{n!}$ называется \rindex{объём}\rindex{форма объёма}\emph{канонической формой объёма} на $(M, \Omega)$.
Полезно иметь в виду два элементарных примера симплектических многообразий:
ориентируемую поверхность с формой площади и линейное пространство $\RR^{2n} (p_1,\z\dots, p_n, q_1,\z\dots, q_n)$ с формой $\omega = \sum^n_{j = 1} dp_j \wedge dq_j$.
Второй пример очень важен с точки зрения классической \rindex{теорема Дарбу}\emph{теоремы Дарбу} \cite{MS}.
Она утверждает, что локально каждое симплектическое многообразие выглядит как $(\RR^{2n}, \omega)$.
Иными словами, для каждой точки $M$ можно выбрать такие локальные координаты $(p, q)$, что в этих координатах $\Omega$ записывается как $\sum^n_{j = 1} dp_j \wedge dq_j$.
Мы называем $(p, q)$ \rindex{канонические координаты}\emph{каноническими локальными координатами}.

Пусть $F$ — гладкая функция на $M$.
Векторное поле $\xi$ на $M$ называется \rindex{гамильтоново векторное поле}\emph{гамильтоновым векторным полем} гамильтониана $F$, если оно поточечно удовлетворяет линейному алгебраическому уравнению $i_\xi \Omega \z= -dF$.
Элементарное рассуждение из линейной алгебры (основанное на невырожденности $\Omega$) даёт, что $\xi$ всегда существует и единственно \cite{MS}.
Иногда $\xi$ обозначают $\sgrad F$ (\rindex{косой
    градиент}\emph{симплектический градиент} $F$).

\begin{ex}{Упражнение}\label{1.3.A}
Докажите, что в канонических локальных координатах $(p, q)$ на
$M$ выполняется равенство \index[symb]{$\sgrad$} $\sgrad F = (-\tfrac{\partial F}{\partial q},\tfrac{\partial F}{\partial p})$
\end{ex}

\begin{ex}{Упражнение}\label{1.3.B}
Пусть  $\phi\: M \to M$ — симплектический диффеоморфизм
(то есть $\phi^\ast \Omega = \Omega$).
Докажите, что $\sgrad (F \circ \phi^{-1}) \z= \phi_\ast \sgrad F$ для любой функции $F$ на $M$.
Это свойство, конечно же, отражает то, что операция $\sgrad$ определяется бескординатным образом.
\end{ex}

В классической механике энергия определяет эволюцию системы.
Энергия — это семейство функций $F_t$ на $M$, которое зависит от
дополнительного временного параметра $t$.
Время $t$ определено на некотором интервале $I$.
Эквивалентно, можно рассматривать энергию как одну единственную функцию $F$ на $M \times I$.
Мы используем оба варианта на протяжении всей книги, сохраняя обозначение $F_t (x) = F (x, t)$.
Традиционно $F$ называется \rindex{гамильтониан}\emph{гамильтонианом (зависящим от
  времени)}.

{\sloppy

Эволюция системы описывается \rindex{уравнение Гамильтона}\emph{уравнением Гамильтона} $\dot x \z= \sgrad F_t (x)$.
В локальных канонических координатах $(p, q)$ на $M$ уравнение
Гамильтона имеет знакомый вид (ср. с \ref{1.3.A})
\[
\begin{cases}
\quad\dot p &= - \tfrac{\partial F}{\partial q} (p, q, t),\\
\quad\dot q &= \tfrac{\partial F}{\partial p} (p, q, t).
\end{cases}
\]

}

Введём линейное функциональное пространство \index[symb]{$\A$}\index[symb]{$\A(M)$}$\A=\A(M)$, которое далее играет важную роль.
Если $M$ замкнуто, определим $\A(M)$ как пространство всех гладких функций на $M$ с нулевым средним относительно канонической формы объёма.
Если $M$ открыто, то определим $\A(M)$ как пространство всех гладких функций с компактным носителем.

\begin{ex}{Определение}\label{1.3.C}
Пусть $I \subset \RR$ — интервал.
Гамильтониан $F$ (зависящий от времени) на $M \times I$ называется \rindex{нормализованный гамильтониан}\textbf{нормализованным}, если $F_t$ принадлежит $\A$ при всех $t$.
Если $M$ открыто, то мы дополнительно требуем, чтобы существовало компактное подмножество $M$, содержащее носители функций $F_t$ при всех $t \in I$.
\end{ex}

Далее мы рассматриваем только нормализованные гамильтонианы.
Приведём пару соображений в пользу этого соглашения.

Прежде всего, на открытых многообразиях необходимо наложить некоторые ограничения на поведение гамильтонианов на бесконечности.
В противном случае решения гамильтонова уравнения могут
убежать на бесконечность за конечное время и, таким образом, гамильтонов поток может быть не определен.
Важной особенностью приведенного выше определения является то, что зависящее от времени гамильтоново векторное поле $\sgrad F_t$ нормализованного гамильтониана $F$ имеет компактный носитель.
Таким образом, когда $I$ содержит ноль, такое поле определяет поток с
компактным носителем, то есть Определение~\ref{1.3.C} вписывается в идеологию предыдущего раздела.

Во-вторых, как на открытых, так и на замкнутых многообразиях отображение, переводящее функцию из $\A$ в его гамильтоново векторное поле, инъективно.
Действительно, гамильтоново векторное поле определяет соответствующий гамильтониан однозначно с точностью до аддитивной константы — ясно, что наша нормализация запрещает добавлять константы!
Это свойство нормализованных гамильтонианов будет полезно.

\section[Группа гамильтоновых диффеоморфизмов]{Группа гамильтоновых\\диффеоморфизмов}\label{1.4}

Пусть $F\: M \times I \to \Rbb$ — нормализованный гамильтониан, зависящий от времени.
Предположим, что $I$ содержит ноль.
Рассмотрим поток $\{f_t\}$ зависящего от времени векторного поля $\sgrad F_t$.
Мы будем говорить, что $\{f_t\}$ — \rindex{гамильтонов поток}\emph{гамильтонов поток}, порождённый $F$.
Каждый диффеоморфизм $f_a$, $a \in I$, называется \rindex{гамильтонов диффеоморфизм}\emph{гамильтоновым диффеоморфизмом}.
Из определения ясно, что гамильтоновы диффеоморфизмы имеют компактный носитель.

\begin{ex}[(репараметризация потоков)]{Упражнение}\label{1.4.A}
\rindex{параметризация потока}
Пусть $\{f_t\}, t \z\in [0; a]$ — гамильтонов поток, порождённый нормализованным гамильтонианом $F (x, t)$.
Докажите, что $\{f_{at}\}$, $t \in [0; 1]$, тоже является гамильтоновым потоком, порождённым гамильтонианом $aF (x, at)$.
Следовательно, любой гамильтонов диффеоморфизм на самом деле является отображением некоторого гамильтонова потока с единичным временем.
Более общо, покажите, что для любой гладкой функции $b (t)$ с $b (0) =
0$ поток $\{f_{b (t)}\}$ является гамильтоновым потоком, нормализованный
гамильтониан которого равен $\tfrac{\d b}{\d t} (t) F (x, b (t))$.
\end{ex}

Ключевым свойством гамильтоновых диффеоморфизмов является то, что они сохраняют симплектическую форму $\Omega$.
В самом деле, пусть $\xi$ — гамильтоново векторное поле функции $F$ на $M$.
Все, что нам нужно проверить, — это равенство нулю производной Ли $L_\xi \Omega$.
Его видно из следующих вычислений: 
\[L_\xi \Omega = i_\xi (\d\Omega) + \d (i_\xi \Omega) = -\d\d F = 0.\]
Обозначим через \index[symb]{$\Ham (M, \Omega)$}$\Ham (M, \Omega)$ множество всех гамильтоновых диффеоморфизмов.

Путь диффеоморфизмов в $\Ham (M, \Omega)$ называется \rindex{гамильтонов путь}\emph{гамильтоновым путём}.
Естественно назвать поток со значениями в $\Ham (M, \Omega)$ гамильтоновым потоком.
Однако, немного подумав, мы понимаем, что попали в беду.
Действительно, гамильтоновы потоки уже определены иначе.
Ведь совсем неясно, что векторное поле, соответствующее гамильтонову пути, является гамильтоновым векторным полем!
К счастью, это правда.
Это чрезвычайно важное утверждение доказал \rindex{Баньяга}Баньяга \cite{B1}.
Вот его точная формулировка.

\begin{thm}{Предложение}\label{1.4.B}
Для каждого гамильтонова пути $\{f_t\}$, $t \z\in I$, существует (зависящий от времени) нормализованный гамильтониан $F\: M \times I \to \RR$ такой, что 
\[\tfrac \d{\d t} f_t x = \sgrad F_t (f_t x)\]
при всех $x \in M$ и $t \in I$.
\end{thm}

Функция $F$ называется нормализованным гамильтонианом
пути $\{f_t\}$.

Обсудим этот результат.
Сначала предположим, что многообразие $M$ удовлетворяет следующему
топологическому условию: его первая группа когомологий де Рама с компактными носителями равна нулю, то есть $H_{\mathrm{comp}}^1(M;\RR) = 0$.
Чтобы иметь перед собой пример, представьте себе двумерную сферу или линейное пространство.
В этом случае приведённое выше предложение доказывается очень легко.
Обозначим через $\xi_t$ векторное поле, порождённое $f_t$.
Поскольку гамильтоновы диффеоморфизмы сохраняют симплектическую форму $\Omega$, имеем $L_{\xi_t} \Omega = 0$.
Следовательно, $\d i_{\xi_t} \Omega = 0$, то есть $i_{\xi_t} \Omega$ — замкнутая форма.
Ввиду нашего топологического условия форма $i_{\xi_t} \Omega$ точна.
Следовательно, существует единственное гладкое семейство функций $F_t (x) \in \A$ такое, что $-dF_t = i_{\xi_t} \Omega$.
Отсюда следует, что $F (x, t)$ — нормализованный гамильтониан $\{f_t\}$, что завершает доказательство.
Однако если $H^1_{\mathrm{comp}} (M;\RR) \ne 0$, то неясно, точны ли замкнутые формы $i_{\xi_t} \Omega$.
В этом случае нужно дополнительно использовать, что каждое отдельное $f_t$ гамильтоново.
Это требует некоторых новых идей, см. \cite{B1}, \cite{MS}.

{\sloppy 

\begin{ex}{Замечание}\label{1.4.C}
Пусть \index[symb]{$\Symp(M,\Omega)$}$\Symp(M,\Omega)$ — группа всех диффеоморфизмов $f$ с компактным носителем на $M$, сохраняющих $\Omega$, то есть $f^\ast\Omega = \Omega$.
Такие диффеоморфизмы называются \rindex{симплектоморфизм}\emph{симплектоморфизмами}.
Обозначим через \index[symb]{$\Symp_0(M,\Omega)$}$\Symp_0(M,\Omega)$ компоненту линейной связности единицы в $\Symp(M,\Omega)$.
По определению она содержит те симплектоморфизмы $f$, которые можно соединить путем симплектоморфизмов с тождественным отображением.
Рассуждая точно так же, как в доказательстве \ref{1.4.B}, мы видим, что если $H^1_{\mathrm{comp}} (M;\RR) = 0$, то такой путь является гамильтоновым.
В этом случае
\[\Symp_0 (M, \Omega) = \Ham (M, \Omega).\]
Если $H^1_{\mathrm{comp}} (M;\RR) \ne 0$, то последнее равенство может не выполняться.
Например, рассмотрим двумерный тор $\TT^2=\RR^2(p,q)/\ZZ^2$ с формой площади $\Omega = dp \wedge dq$.
Сдвиг $(p, q) \to (p + a, q)$, очевидно, лежит в $\Symp_0 (\TT^2, \Omega)$.
Однако можно показать, что для нецелого $a$ он не гамильтонов.
С другой стороны, разница между $\Ham$ и $\Symp_0$ не слишком велика, и её можно описать довольно просто, см. \ref{sec:14.1}.
\end{ex}

}

Следующее предложение содержит замечательную элементарную формулу, которая многократно используется ниже.

\begin{thm}[(Гамильтониан произведения)]{Предложение}\label{1.4.D}
Рассмотрим два гамильтонова пути $\{f_t\}$ и $\{g_t\}$.
Пусть $F$ и $G$ — их нормализованные гамильтонианы.
Тогда путь произведения $h_t = f_t g_t$ является гамильтоновым путём, порождённым нормализованным гамильтонианом 
\[H(x,t) = F(x,t) + G(f_t^{-1} x, t).\]
\end{thm}

{\sloppy 

Докажем эту формулу.
Нам дано, что 
\[\tfrac{\d}{\d t} f_t = \sgrad F_t
\quad\text{и}\quad
\tfrac{\d}{\d t}g_t = \sgrad G_t
\]
Таким образом, 
\[\tfrac{\d}{\d t} (f_t g_t) = \sgrad F_t + f_{t\ast} \sgrad G_t.\]
Согласно упражнению \ref{1.3.B} второе слагаемое правой части равно $\sgrad  (G \circ f_t^{-1})$.
Отсюда 
\[\tfrac{\d}{\d t} h_t = \sgrad  (F_t + G \circ f_t^{-1}) = \sgrad H_t.\]
Формула доказана.

}

Теперь мы в состоянии обосновать название раздела.

\begin{thm}{Предложение}
Множество гамильтоновых диффеоморфизмов является группой относительно композиции.
\end{thm}

Действительно, возьмём два гамильтоновых диффеоморфизма $f$ и $g$.
Ввиду \ref{1.4.A} можно записать $f = f_1$ и $g = g_1$ для некоторых гамильтоновых потоков $\{f_t\}$, $\{g_t\}$, определённых для $t \in [0; 1]$.
Предложение \ref{1.4.D} означает, что путь $\{f_t g_t\}$ является гамильтоновым потоком.
Таким образом, его отображение $f g$ в единичное время является гамильтоновым диффеоморфизмом.
В частности, множество $\Ham (M, \Omega)$ замкнуто относительно композиции диффеоморфизмов.
Осталось проверить, что $f^{-1}$ — гамильтонов диффеоморфизм.
Это следует из следующего упражнения.

\begin{ex*}{Упражнение} Покажите, что путь $\{f_t^{-1}\}$ является гамильтоновым потоком, порождённым гамильтонианом $-F (f_t x, t)$.
\emph{Подсказка:} продифференцируйте тождество $f_t \circ f_t^{-1} = \1$ по $t$ и рассуждайте, как в доказательстве \ref{1.4.D}.
\end{ex*}

\begin{ex}{Замечание}\label{1.4.F}
В дифференциальной геометрии обычно имеют дело с группами
преобразований, сохраняющих определённые структуры на многообразии
(например, группу $\Symp(M,\Omega)$ всех симплектоморфизмов).
Группа гамильтоновых диффеоморфизмов не обладает таким понятным
описанием (гамильтоновы диффеоморфизмы не определяются как морфизмы в
определённой естественной категории).
Это приводит к очень неожиданным сложностям. 
Например, следующий вопрос (известный как \rindex{гипотеза
  потока}\emph{гипотеза потока},%
\footnote{
Гипотезу потока доказал Оно \cite{O06}.
Однако вопрос о замкнутости группы гамильтоновых диффеоморфизмов в $C^{0}$-топологии
ещё далёк от разрешения.\dpp}
см. главу~\ref{chap:14}) все ещё
открыт для большинства симплектических многообразий $M$.
Пусть $M$ замкнуто.
Предположим, что некоторая последовательность гамильтоновых
диффеоморфизмов $C^\infty$-сходится к симплектическому диффеоморфизму
$f$.
Является ли $f$ гамильтоновым?
\end{ex}

{\sloppy 

Чрезвычайно полезно думать о $\Ham (M, \Omega)$ в терминах теории групп Ли.
Эта точка зрения является основной для развития необходимой нам геометрической интуиции.
Давайте разработаем этот язык.
Мы будем рассматривать $\Ham (M, \Omega)$ как подгруппу Ли группы всех диффеоморфизмов $M$.
Таким образом, алгебра Ли%
\footnote{\label{footnote} Как векторное пространство алгебра Ли по определению является касательным пространством к группе в единице.
Касательные пространства к группе во всех остальных точках отождествляется с алгеброй Ли с помощью {}\emph{правых сдвигов} группы.}
группы $\Ham (M, \Omega)$ — это просто алгебра всех векторных полей $\xi$ на $M$ вида 
\[\xi(x)=\left.\frac{\d}{\d t}\right|_{t=0}f_tx\]
где $\{f_t\}$ — гладкий путь в $\Ham (M, \Omega)$ с $f_0 = \1$.
Каждое такое поле гамильтоново.
В самом деле, $\xi = \sgrad F_0 (x),$ где $F (x, t)$ — (единственный!) нормализованный гамильтониан, порождающий путь, и $F_0 (x) = F (x, 0)$.
Отметим, что $F_0 \in \A$.
Наоборот, для любой функции $F \in \A$ векторное поле $\sgrad F$ по определению является производной в точке $t = 0$ соответствующего гамильтонова потока.
Мы заключаем, что \rindex{алгебра Ли группы $\Ham (M, \Omega)$}\emph{алгебру Ли группы $\Ham (M, \Omega)$ можно отождествить с $\A$}.

\begin{ex}{Упражнение}\label{1.4.G}
Покажите, что на этом языке вектор, касательный к гамильтонову пути $\{f_t\}$ в точке $t = s$, является функцией $F_s \in \A$.
\emph{Подсказка:} идентификация касательных пространств группы осуществляется с помощью правого сдвига (см. сноску \ref{footnote}%
).
Таким образом, рассматриваемый касательный вектор отождествляется с
вектором, касательным к пути $\{f_tf_s^{-1}\}$ в точке $t \z= s$.
\end{ex}

Следующее важное понятие — присоединённое действие группы Ли на её алгебре Ли.
Напомним, что эта операция определяется следующим образом.
Выберем элемент $f$ группы $\Ham (M, \Omega)$ и элемент~$G$ алгебры Ли $\A$.
Пусть $\{g_t\}$, $g_0 = \1$, — путь на группе, касающийся $G$.
В нашем случае условие касания означает, конечно, что нормализованный гамильтониан потока $\{g_t\}$ в момент времени $t = 0$ равен $G$.
По определению \index[symb]{$\Ad_f G$}
\[\Ad_f G = \left.\frac{\d}{\d t}\right|_{t=0} fg_tf^{-1}.\]
Продифференцировав, получаем, что векторное поле в правой части равно $f_\ast \sgrad G$, а это в точности $\sgrad  (G\circ f^{-1})$ ввиду упражнения \ref{1.3.B}.
Возвращаясь к нашему отождествлению, получаем, что 
\[\Ad_f G = G \circ f^{-1}.\]
Таким образом, \rindex{присоединённое действие}\emph{присоединённое действие $\Ham (M, \Omega)$ на $\A$ — это просто обычное действие диффеоморфизмов на функциях}.

Наконец, давайте обсудим скобку Ли на $\A$.
Выберем два элемента $F,G \in \A$, и пусть $\{f_t\}$, $f_0 = \1$, --- гамильтонов путь, касающийся $F$ в~$0$.
Скобка Ли \index[symb]{$\{F, G\}$}$\{F, G\}$ элементов $F$ и $G$ называется \rindex{скобка Пуассона}\emph{скобкой Пуассона} и определяется следующим образом: 
\[\{F, G\}=\left.\frac{\d}{\d t}\right|_{t=0}Ad_{f_t}G\]
Вычисляя выражение в правой части, получаем, что 
\[\{F, G\} = -\d G (\sgrad F) = \Omega (\sgrad G, \sgrad F).\]
Отметим, что в терминах векторных полей скобка Ли с точностью до знака совпадает с обычным коммутатором.
Читателю предлагается убедиться в том, что 
\[[\sgrad F, \sgrad G] = -\sgrad  \{F, G\},\]
где коммутатор $[X, Y]$ двух векторных полей определяется как $L_{[X, Y]} \z= L_X L_Y -L_Y L_X$.

\begin{framed}
\parbf{Внимание:} 
Часто используются противоположные знаки в следующих определениях:
гамильтоново векторное поле,
скобка Пуассона,
коммутатор векторных полей
и кривизна связности.
\end{framed}

\begin{ex}{Пример}\label{1.4.H}
Рассмотрим единичную сферу $S^2$ в евклидовом пространстве $\RR^3$.
Пусть $\Omega$ — индуцированная форма площади на сфере.
Группа $\SO(3)$ действует на $S^2$ диффеоморфизмами, сохраняющими площадь.
Поскольку группа $\SO(3)$ линейно связна, она содержится в $\Symp_0 (S^2)$.
Применяя \ref{1.4.C}, получаем, что $\Symp_0 (S^2) \z= \Ham (S^2)$, а значит, $\SO (3)$ является подгруппой в $\Ham (S^2)$.
В частности, каждый элемент алгебры Ли $\so(3)$ однозначно представляется нормализованным гамильтонианом на $S^2$.
Давайте подробно опишем это соответствие.
Отождествим $\so (3)$ с $\RR^3$ следующим образом.
Каждый вектор $a \in \RR^3$ рассматривается как кососимметричное преобразование $x \mapsto [x, a]$ пространства, где квадратные скобки обозначают стандартное векторное произведение.
Отождествим касательную плоскость к $S^2$ в точке $x$ с ортогональным
дополнением к $x$ в объемлющем пространстве.
По тавтологическим причинам гамильтоново векторное поле $v$ потока $x \z\mapsto \exp (ta) x$ на сфере задаётся формулой $v (x) = [x, a]$.
Мы утверждаем, что \textit{соответствующий нормализованный гамильтониан является
функцией высоты $F(x)=\langle a,x\rangle$} в направлении вектора $a$, 
где угловые скобки обозначают скалярное произведение.

Прежде всего, отражение относительно ортогонального дополнения к $a$ переводит $F$ в $-F$, таким образом, $F$ имеет нулевое среднее.
Далее заметим, что $\Omega (\xi, \eta) = \langle\eta, [x, \xi]\rangle$ при $\xi, \eta \in \T_x S^2$.
Обозначим через $a'$ ортогональную проекцию $a$ на $\T_x S^2$.
Таким образом, $\Omega (\xi, v(x)) \z= \langle[x, a'], [x, \xi]\rangle$ при любом $\xi \in \T_x S^2$.
Поскольку векторное произведение с $x$ является ортогональным преобразованием касательного пространства $\T_x S^2$, последнее выражение равно $\langle a', \xi\rangle$,
а это в точности $\d F(\xi)$, и утверждение следует.
\end{ex}

\section[\texorpdfstring{Алгебраические свойства $\Ham(M,\Omega)$}{Алгебраические свойства Ham(M,Ω)}]%
{Алгебраические свойства $\bm{\Ham(M, \Omega)}$}

Алгебраические свойства группы гамильтоновых диффеоморфизмов изучались \rindex{Баньяга}А. Баньягой \cite{B1,B2}.
В частности, он доказал следующий поразительный результат.
Напомним, что группа $D$ называется \rindex{простая группа}\emph{простой}, если каждая её нормальная подгруппа тривиальна, то есть она либо $\{\1\}$, либо вся $D$.

\begin{thm}{Теорема}\label{1.5.A}
Пусть $(M, \Omega)$ — замкнутое симплектическое многообразие.
Тогда группа $\Ham (M, \Omega)$ проста.
\end{thm}

Существует версия этого утверждения и для открытого $M$.

Обратите внимание на то, что абелева группа проста тогда и только тогда, когда каждый элемент порождает всю группу (значит, это конечная циклическая группа простого порядка).
Таким образом, вообще говоря, простые группы далеки от абелевых.
Ниже приводится элементарное утверждение, поясняющее этот принцип для группы гамильтоновых диффеоморфизмов.
Нам он потребуется в следующей главе.

\begin{thm}{Предложение}\label{1.5.B}
Пусть $(M, \Omega)$ — симплектическое многообразие и $U \subset M$ — непустое открытое подмножество.
Тогда существуют такие $f, g \in \Ham (M, \Omega)$, что $\supp (f)$, $\supp (g) \subset U$ и $f g \ne gf$.
\end{thm}

Доказательство основано на следующем предложении.

\begin{thm}{Предложение}\label{1.5.C}
Пусть $\{f_t\}$ и $\{g_t\}$ — гамильтоновы потоки, порождённые независимыми от времени нормализованными гамильтонианами $F$ и $G$ соответственно.
Если $f_t g_t = g_t f_t$ при всех~$t$, то $\{F, G\} = 0$.
\end{thm}

\parit{Доказательство.}
Согласно \ref{1.4.D} гамильтонианы, соответствующие потокам $f_t g_t$ и $g_t f_t$ равны
\[F(x)+G(f_t^{-1} (x))
\quad\text{и}\quad
G (x) + F (g_t^{-1}(x)).
\]
Поскольку оба определяют один и тот же поток, получаем, что 
\[F (x) + G (f_t^{-1} (x)) = G (x) + F (g_t^{-1} (x))\]
при всех $t$.
Продифференцировав по $t$, получаем
\[\d G(-\sgrad F)=\d F(-\sgrad G),\]
и, значит $\{F, G\} = \{G, F\}$.
В силу антикоммутативности скобки Ли получаем $\{F, G\} = 0$.
\qeds

\parit{Доказательство \ref{1.5.B}.} 
Выберем точку $x\in U$ и касательные векторы $\xi, \eta \in \T_x U$, что $\Omega (\xi, \eta) \ne 0$.
Далее выберем ростки функций $F$ и $G$ (см. упражнение ниже), для которых $\sgrad F (x) = \xi$, $\sgrad G (x) = \eta$.
Продолжим эти функции нулём за пределы $U$.
Если $M$ открыто, то задача решена.
Если $M$ замкнуто, добавим константу, чтобы гарантировать, что $F$ и $G$ имеют нулевое среднее.
Таким образом, функции $F$ и $G$ принадлежат $\A$.
Кроме того, они постоянны вне $U$, поэтому носители соответствующих гамильтоновых диффеоморфизмов $f_t$ и $g_t$ лежат в $U$.
Поскольку $\{F, G\} \ne 0$, мы видим, что при некотором $t$ диффеоморфизмы $f_t$ и $g_t$ не коммутируют.
\qeds

\begin{ex*}{Упражнение}
Используя локальные канонические координаты в $x$, докажите, что $F$ и $G$ в доказательстве \ref{1.5.B} действительно существуют.
\end{ex*}

Завершим этот раздел формулировкой следующего результата \rindex{Баньяга}А. Баньяги \cite{B2}.

\begin{thm}{Теорема}\label{1.5.D}
Пусть $(M_1, \Omega_1)$ и $(M_2, \Omega_2)$ — два замкнутых симплектических многообразия, группы гамильтоновых диффеоморфизмов которых изоморфны.
Тогда многообразия \rindex{конформно симплектоморфные многообразия}\emph{конформно симплектоморфны},
то есть существуют диффеоморфизм $f\: M_1 \to  M_2$ и число $c \ne 0$, что $f^\ast \Omega_2 = c\Omega_1$. 
\end{thm}

Иными словами, алгебраическая структура группы гамильтоновых диффеоморфизмов определяет симплектическое многообразие с точностью до множителя.

\chapter{Знакомство с геометрией}\label{chap:2}

В этой главе обсуждаются биинвариантные финслеровы метрики на группе
гамильтоновых диффеоморфизмов и определяется хоферовская геометрия. 

\section{Вариационная задача}\label{2.1}

Сколько нужно энергии для получения данного гамильтонова диффеоморфизма $\phi$? 
Этот естественный вопрос можно формализовать следующим образом.
Рассмотрим всевозможные гамильтоновы потоки $\{f_t\}$, $t \in [0; 1]$ такие, что $f_0 = \1$ и $f_1 = \phi$.
Для каждого потока возьмём, соответсвующий ему нормализованный гамильтониан $F_t(x)$ и «измерим его величину».
Затем минимизируем результат измерения по всем таким потокам.
Остаётся понять, что значит «измерить его величину».
Напомним, что при любом $t$ функция $F_t$ является элементом алгебры Ли $\A$.
Возьмём любую безкоординатную норму $\|\ \|$ на $\A$,
то есть потребуем, чтоб
\begin{equation}
 \|H \circ \psi^{-1}\|
= \| H \|
\quad\text{при всех}\quad
H \in \A \quad\text{и}\quad  \psi \in \Ham (M, \Omega).
\label{2.1.A}
\end{equation}
Теперь определим «величину гамильтониана» как  $\int_0^1\| F_t \|\, dt.$
Собрав все части этой процедуры, получаем следующую вариационную задачу 
\begin{equation}
\inf\int_0^1 \| F_t \|\,dt, 
\label{2.1.B}
\end{equation}
где $\phi$ фиксировано, а точная нижняя грань берётся по всем потокам
$\{f_t\}$, описанным выше.

\section[\texorpdfstring{Биинвариантные геометрии на $\Ham(M, \Omega)$}{Биинвариантные геометрии на Ham(M,Ω)}]{Биинвариантные геометрии\\ на $\bm{\Ham(M, \Omega)}$}\label{2.2}

Приведённую вариационную задачу можно переформулировать в чисто геометрических терминах.
Для этого надо вспомнить понятие \rindex{финслерова
  структура}\emph{финслеровой структуры} на многообразии. 
Мы говорим, что многообразие $Z$ наделено финслеровой структурой, если его касательные пространства $\T_z Z$ оснащены нормой, которая гладко зависит от точки $z \in Z$.
Конечно же, римановы структуры являются частным случаем этого понятия.
Однако в общем случае, нормы не обязаны задаваться квадратичной формой.
Финслерову структуру можно использовать при определения длины кривой
точно так же, как риманову 
\[\length \{z(t)\}_{t\in [a; b]} =\int_a^b \norm(\dot z (t))\,dt.\]
Кроме того, можно ввести {}\emph{расстояние} между двумя точками $z$ и $z'$ в $Z$ как точную нижнюю грань длин всех кривых, соединяющих $z$ с~$z'$.

Вернёмся к задаче предыдущего раздела.
Поскольку все касательные пространства группы $\Ham(M, \Omega)$ можно
отождествить с $\A$ (см. раздел \ref{1.4}), каждый выбор нормы
$\|\ \|$ на $\A$ определяет финслерову структуру на группе.
Таким образом, можно определить длину гамильтонова пути и расстояние между двумя гамильтоновыми диффеоморфизмами.
В частности, длина гамильтонова пути $\{f_t\}$, $t \in [a; b]$ с нормализованным гамильтонианом $F$ определяется следующим образом \index[symb]{$\length\{f_t\}$}
\[\length\{f_t\}=\int_a^b \| F_t \|\,dt.\]
А расстояние между двумя гамильтоновыми диффеоморфизмами $\phi$ и $\psi$ определяется как \index[symb]{$\rho$}
\[\rho (\phi, \psi) = \inf \length \{f_t\},\] 
где точная нижняя грань берётся по всем гамильтоновым путям $\{f_t\}$, $t \in [a; b]$ с $f_a = \phi$ и $f_b = \psi$.
Конечно, длина пути не зависит от параметризации, таким образом, в приведённом определении расстояния можно считать, что $a = 0$ и $b = 1$.
На этом языке, решение вариационной задачи
\ref{2.1.B} — это в точности расстояние $\rho (\1, \phi)$!

Следующие свойства $\rho$ легко проверить,
они даются как упражнения:
\begin{itemize}
\item $\rho (\phi, \psi) = \rho (\psi, \phi)$;
\item неравенство треугольника: $\rho (\phi, \psi) + \rho (\psi, \theta) \ge \rho (\phi, \theta)$;
\item $\rho (\phi, \psi) \ge 0$.
\end{itemize}

Напомним теперь об условии \ref{2.1.A}, которое было наложено на норму $\|\ \|$ в разделе \ref{2.1}.
На геометрическом языке это означает, что норма инвариантна относительно присоединённого действия группы на её алгебре Ли (см. \ref{1.4}).
Далее мы имеем дело только с такими нормами.

\begin{ex*}{Упражнение}
Покажите что из \ref{2.1.A} следует, что функция $\rho$ — \rindex{биинвариантная (псевдо)метрика}\emph{биинвариантна},%
\footnote{Без \ref{2.1.A} мы получаем лишь правоинваринатность $\rho$, то есть $\rho (\phi, \psi) \z= \rho (\phi\theta, \psi\theta)$.
Такие метрики играют важную роль в гидродинамике, см. \cite{AK}.}
то есть
\[\rho (\phi, \psi) = \rho (\phi \theta, \psi \theta) = \rho (\theta\phi, \theta\psi)\]
при всех $\phi, \psi, \theta \in \Ham (M, \Omega)$.
\end{ex*}

Было бы честней пока назвать функцию $\rho$ \rindex{псевдометрика}\emph{псевдометрикой}.
Действительно, как мы уже знаем, она удовлетворяет всем аксиомам метрики, за исключением, возможно, {}\emph{невырожденности}: 
\begin{equation}
\rho (\phi, \psi)> 0
\quad\text{при}\quad
\phi \ne \psi.
\label{eq:2.2.A}
\end{equation}

Даже в конечномерном случае, невырожденность проверить не так просто.
В этом случае, доказательство использует локальную компактность многообразия.
При этом наша группа $\Ham (M, \Omega)$ бесконечномерна и у неё нет никаких свойств компактности.
Таким образом, у нас пока нет причин верить, что условие \ref{eq:2.2.A} выполнено.
Скоро мы увидим, что свойство невырожденности $\rho$ весьма чувствительно к выбору нормы~$\|\ \|$.

\section[\texorpdfstring{Выбор нормы: $L_p$ или $L_\infty$}{Выбор нормы: Lₚ или L∞}]{Выбор нормы: \bm{$L_p$} или \bm{$L_\infty$}}

Среди норм на $\A$, удовлетворяющих предположению инвариантности
\ref{2.1.A} существует очень естественный класс, который включает
нормы $L_p$, $p \in [1;\infty)$
\[\|H\|_p=\left(\int_M|H|^p\Vol\right)^{\frac1p},\]
и $L_\infty$-норму 
\[\|H\|_\infty = \max H - \min H.\]
Обозначим через $\rho_p$ и $\rho_\infty$ соответствующие псевдометрики.

\begin{thm}{Теорема}\label{2.3.A}
Псевдометрика $\rho_p$ вырождена при всех конечных $p \in [1;\infty)$.
Более того, если многообразие замкнуто, то все такие $\rho_p$ тождественно обращаются в нуль.
\end{thm}

Этот результат получен в \cite{EP}.
Доказательство представлено в этой главе (см. также книги \cite{HZ}, \cite{MS}, \cite{AK}).
Следующая теорема показывает разительный контраст между $L_p$ и $L_\infty$.

\begin{thm}{Теорема}\label{2.3.B}
Псевдометрика $\rho_\infty$ невырождена.
\end{thm}

Эту теорему%
\footnote{Историческое отступление, приведённое ниже, отражает моё личное мнение.
Допускаю, что другие участники этих событий видят дело иначе.}
сформулировал Хофер, и он доказал её в случае $M = \RR^{2n}$
использованием бесконечномерных вариационных методов в \cite{H1}.
\rindex{Витербо}Витербо в \cite{V1} вывел её в случае $M = \RR^{2n}$,
используя разработанную им теорию производящих функций.
Для обоих, Хофера и Витербо, толчком послужил вопрос, заданный
\rindex{Элиашберг}Элиашбергом в частной беседе.
В \cite{P1} это утверждение было распространено на широкий класс
симплектических многообразий с «хорошим» поведением на бесконечности
и, в частности, на все замкнутые симплектические многообразия, у
которых класс когомологий симплектической формы является
рациональным.
Подход \cite{P1} основан на теории псевдоголоморфных кривых \rindex{Громов}Громова.
Наконец \cite{LM1} Лалонд и Макдафф доказали \ref{2.3.B} в полной
общности, используя теорию Громова. 
К настоящему времени найдены другие доказательства различных частных случаев этой теоремы, см., например, \rindex{Чеканов}\rindex{Шварц}\cite{Ch,O3,Sch3}.
Первоначальное доказательство Хофера подробно представленo в книге
\cite{HZ}.
Доказательство \rindex{Лалонд}Лалонда и Макдафф изложены в книге \cite{MS} и в обзоре \cite{L}.
Ниже приводится другое доказательство в случае $M = \RR^{2n}$, которое следует из \cite{P1}.
Все известные доказательства основаны на «жёстких» методах.%
\footnote{Более того, на мой вкус, все доказательства далеко не
  прозрачны.
  Рассуждение, представленное в главе 3, не является исключением.
  Я твердо верю, что будет найдено {}\emph{правильное
    объяснение} этому фундаментальному утверждению.}%
\footnote{  
  Концептуальное доказательство невырожденности метрики Хофера
  было получено сочетанием теории Флоера (кратко обсуждаемой далее в
  этой книге) и теории модулей персистентности и штрих-кодов, возникшей в
  задачах топологического анализа данных, и введённой в симплектическую
  топологию Полтеровичем и Шелухиным в~\cite{PS16}. Мы отсылаем
  читателя к~\cite{PS16,UZ,PRSZ} для ознакомления с этими
  идеями.\dpp}

\section{Энергия смещения}\label{sec:2.4}

Какие инвариантные нормы $\|\ \|$ на $\A$ (то есть, удовлетворяющие \ref{2.1.A}) приводят к невырожденным метрикам $\rho$?
Ниже описана очень полезная переформулировка этого вопроса, которая в конечном счёте позволит доказать теорему \ref{2.3.A}.
Она основана на понятии энергии смещения — прекрасной идее, введённой \rindex{Хофер}Хофером \cite{H1}.
Пусть $\rho$ — биинвариантное псевдометрика на $\Ham (M, \Omega)$, и пусть $A$ — ограниченное подмножество $M$.

\begin{ex*}{Определение}
\rindex{энергия смещения}\emph{Энергия смещения} $A$ определяется как \index[symb]{$\e(A)$}
\[\e(A) = \inf \set{\rho (\1, f)}{f \in \Ham (M, \Omega),\  f (A) \cap A = \emptyset}.\]
\end{ex*}

Множество таких $f$ может быть пусто.
Напомним, что нижняя грань пустого множества равна $+\infty$.
Если $\e (A) \ne 0$, то мы говорим, что $A$ имеет положительную энергию смещения.

Отметим пару очевидных, но важных свойств $\e$.
Прежде всего, $\e$ является монотонной функцией подмножеств: если $A \subset B$, то $\e(A) \le \e(B)$.
Во-вторых, $\e$ — инвариантна, то есть, $\e(A) = \e(f(A))$ для любого гамильтонова диффеоморфизма $f$ на $M$.
Доказательства предоставляются читателю.

\begin{ex*}{Пример}
Рассмотрим $(\RR^2, \omega)$ и возьмём открытый квадрат $A$, рёбра которого имеют длину $u$ и параллельны осям координат.
Оценим энергию смещения $A$ относительно расстояния $\rho_\infty$.
Рассмотрим гамильтониан $H (p, q) = up$.
Соответствующая гамильтонова система имеет вид 
\[
\begin{cases}
\quad\dot q &= u,
\\
\quad\dot p &= 0.
\end{cases}
\]
Следовательно, за единичное время, $(p, q)$ переходит в $(p, q + u)$.
Обратите внимание, что движение квадрата происходит в прямоугольнике
$K = \overline{A \cup h (A)}$. 
Рассмотрим \rindex{срезка}\emph{срезку} $F$ гамильтониана $H$ вне малой окрестности $K$.%
\footnote{Пусть $Y$ — замкнутое подмножество многообразия $Z$, и пусть
  $H$ — гладкая функция, определённая в окрестности $V$ множества
  $Y$. 
  Под срезкой $F$ функции $H$ вне малой окрестности $Y$ понимается следующее.
  Выберем окрестность $W\supset Y$, замыкание которой содержится в $V$.
  Возьмём гладкую функцию $a\colon Z \to [0; 1]$, которая равна $1$ на
  $W$ и обращается в нуль вне $V$. 
  Определим $F$ как $aH$ на $V$ и продолжим её нулём вне $V$.}
Обратите внимание, что $F$ — нормализованный гамильтониан (в отличие от $H$).
Поскольку $F = H$ на $K$, гамильтонов поток $f$, порождённый $F$,
по-прежнему смещает квадрат $A$. 
Срезку всегда можно выполнить так, что $L_\infty$-норма $F$ была произвольно близка к осцилляции $H$ на $K$.
Осцилляция определяется как 
\[\max_K H - \min_K H = u^2 - 0 = u^2.\]
Таким образом, получаем, что 
\[\e (A) \le u^2 = \area (A).\]
\end{ex*}

Обратите внимание, что квадрат симплектоморфен диску той же площади в $\RR^2$ (это неверно в более высоких размерностях!).
Таким образом, мы доказали, что $\e (B^2 (r)) \le \pi r^2$.
Глубокий результат \rindex{Хофер}Хофера \cite{H1} утверждает, что на самом деле равенство
достигается во всех размерностях, то есть, $\e (B^{2n} (r)) = \pi r^2$.
Обобщение на произвольные симплектические многообразия можно найти в \rindex{Лалонд}\cite{LM1}.
Нижняя оценка на $\e(B^{2n} (r))$ дана в следующей главе.

В общем случае, если энергия смещения всех непустых открытых подмножеств относительно некоторой биинвариантной псевдометрики $\rho$ положительна, то $\rho$ невырождена.
В самом деле, любой $f \in \Ham (M, \Omega)$ такой, что $f \ne \1$,
должен сместить некоторый маленький шарик $A \subset M$.
Таким образом, получаем, что $\rho (\1, f) \ge \e (A)> 0$.
На самом деле верно и обратное.

\begin{thm}[(\cite{EP})]{Теорема}\label{2.4.A}
Если $\rho$ невырождено, то $\e (A)> 0$ для любого непустого открытого подмножества $A$.
\end{thm}

В доказательстве потребуется следующая лемма.

\begin{thm}{Лемма}\label{2.4.B}
Пусть $A \subset M$ — непустое открытое подмножество.
Тогда $\e (A) \ge \tfrac14 \rho (\1, [\phi, \psi])$ при любых $\phi, \psi \in \Ham (M, \Omega)$ таких, что $\supp (\phi) \subset A$ и $\supp (\psi) \subset A$.
\end{thm}

Здесь $[\phi, \psi]$ обозначает коммутатор $\psi^{-1} \phi^{-1} \psi\phi$.
Теорема сразу следует из леммы ввиду предложения \ref{1.5.B} предыдущей главы.

\parit{Доказательство \ref{2.4.A}.} 
Согласно \ref{1.5.B} существуют $\phi$, $\psi$ с носителем в $A$ такие, что $[\phi, \psi] \ne \1$.
Поскольку $\rho$ невырождено, $\e (A) \z\ge \tfrac14 \rho (\1, [\phi, \psi])> 0$.
\qeds

\parit{Доказательство \ref{2.4.B}.}
Введём временное обозначение $|x|\z\df\rho(\1,x)$
для $x$ из группы $\Ham (M, \Omega)$.
Поскольку $\rho$ биинвариантна,
\[
|xyz|
\ge
|xz|-|y|
\quad\text{и}\quad
|x^{-1}|=|x|
\]
при любых $x,y,z\in \Ham (M, \Omega)$.

Предположим, что существует $h \in \Ham (M, \Omega)$ такой, что $h (A) \z\cap A \z= \emptyset$ (если такого $h$ нету, то дело сделано — в этом случае $\e (A) \z= + \infty$).
Поскольку $h$ смещает с себя множество $A$ содержащее носители $\phi$ и $\psi$,
носители $\phi$ и $\psi^{h}=h^{-1} \psi h$ не пересекаются.
В частности, $[\phi,\psi^{h}]=\1$.
Так как $h$ встречается в $[\phi,\psi^{h}]$ четыре раза, получаем
\[
0=\big|[\phi,\psi^{h}]\big|
\ge
\big|[\phi,\psi]\big| - 4\cdot|h|.
\]
\qedsf

Напомним, что теорема \ref{2.3.A} утверждает, что $L_p$-норма задаёт вырожденную псевдометрику при $p <\infty$;
более того, для замкнутых многообразий она тождественно равна нулю.
Давайте докажем это утверждение и за одно увидим, почему рассуждение не проходит в $L_\infty$-случае.

\parit{Доказательство \ref{2.3.A}.}
Покажем, что энергия смещения маленького шарика обнуляется.
Тогда вырождение псевдометрики $\rho_p$ следует из \ref{2.4.A}.
Пусть $U$ — открытое подмножество в $M$ с каноническими координатами $(x, y)$.
В этих координатах симплектическая форма $\Omega$ равна $\sum dx_i\wedge dy_i$.
Не умаляя общности можно предположить, что $U$ содержит шар $\sum(x^2_j + y_j^2) <10$.
Пусть $A \subset U$ — шар с тем же центром радиуса $\tfrac1{10}$.
Рассмотрим (частично определённый) поток $h_t$, $t \in [0; 1]$ на $U$, который представляет собой простой сдвиг на $t$ по координате $y_1$.
Такой сдвиг порождается (ненормализованым!) гамильтонианом $H (x, y) = x_1$ на $U$.
Ясно, что $h_1 (A) \cap A = \emptyset$.
Пусть $S_t$ это сфера $h_t (\partial A)$.
Рассмотрим новый (зависящий от времени) нормализованный гамильтониан
$G_t = F_t + c_t$, где $F_t$ срезка $H$ вне небольшой окрестности $S_t$, а $c_t$ является (зависящей от времени) постоянной.
Конечно же, $c_t = 0$, если многообразие $M$ открыто, и $c_t = -\Vol(M)^{-1}\int_M F_t \Vol$, если $M$ замкнуто.
Поскольку при любом $t$ функция $G_t$ совпадает с $H$ вблизи $S_t$ с точностью до аддитивной константы, заключаем, что $\sgrad G_t = \sgrad H$ вблизи $S_t$.
Следовательно, поток $\{g_t\}$ гамильтониана $G$ удовлетворяет условию $g_t (\partial A) = h_t (\partial A)$ и, значит, $g_1 (\partial A) \cap \partial A = \emptyset$.
Но отсюда, очевидно, следует, что $g_1 (A) \cap A = \emptyset$.
Заметим теперь, что с помощью отсечения вне очень малых окрестностей $S_t$ можно добиться того, чтобы $L_p$-норма каждой функции $G_t$ стала произвольно малой (в этом отличие от $L_\infty$-нормы!).
Следовательно, $L_p$-энергия смещения $A$ обращается в нуль, что завершает доказательство вырождения псевдометрики $\rho_p$.

Обратимся теперь ко второму утверждению теоремы, где мы предполагаем, что многообразие $M$ замкнуто.
Рассмотрим множество
\[\mathcal{G} = \set{g \in \Ham (M, \Omega)}{\rho (\1, g) = 0}.\]
Выберем $f, g \in \mathcal{G}$.
Конечно же $g^{-1} \in \mathcal{G}$.
Далее из неравенства треугольника следует, что 
\[\rho (\1, f g) = \rho (f^{-1}, g) \le \rho (\1, f) + \rho (\1, g) = 0,\]
так что $\mathcal{G}$ является подгруппой $\Ham (M, \Omega)$.
Поскольку $\rho$ биинвариантна, мы знаем, что если $f \in \mathcal{G}$ и $h \in \Ham (M, \Omega)$, то $hf h^{-1} \in \mathcal{G}$ и, значит, $\mathcal{G}$ — нормальная подгруппа.
По теореме Баньяги \ref{1.5.A}, группа $\Ham (M, \Omega)$ проста.
Следовательно, либо $\mathcal{G} = {\1}$, либо $\mathcal{G} = \Ham (M, \Omega)$.
Мы уже доказали, что $\rho_p$ вырождено, и значит $\mathcal{G} \ne {\1}$.
Таким образом, $\mathcal{G}$ совпадает со всей группой $\Ham (M, \Omega)$.
И значит $\rho_p$ тождественно обращается в нуль.
\qeds

\begin{ex*}{Упражнение}
Докажите, что энергия смещения $S^{2n - 2} \subset \RR^{2n-1} \z\subset \RR^{2n}$ относительно $\rho_\infty$ равна нулю.
\end{ex*}

С другой стороны, как будет видно в \ref{2.4.B} следующей главы,
существуют подмногообразия половинной размерности в $\RR^{2n}$,
которые имеют положительную энергию смещения (ср. с \ref{1.1.C}).

\begin{ex*}{Открытая задача}
Какие инвариантные нормы на $A$ порождают невырожденные функции расстояния $\rho$?
Верно ли, что такие нормы всегда ограничены снизу величиной 
$\const\|\ \|_\infty$?%
\footnote{
Эта задача решена в статье~\cite{BO11} Буховским и Островером и
в~\cite{OW05} Островером и Вагнером, а именно, такая норма обязательно
эквивалентна $\|\cdot\|_{\infty}$. Смотри также~\cite{L20}.\dpp}
Здесь сложность в том, что пока нет классификации $\Ham (M, \Omega)$-инвариантных норм.
Возможный подход состоит в исследовании срезок.
Если срезка способна произвольно уменьшить норму, то по приведённому
выше рассуждению псевдометрика $\rho$ вырождена.
\end{ex*}

\begin{ex*}[\cite{EP}]{Открытая задача} 
Вполне естественно рассматривать отдельно положительную и отрицательную части метрики~$\rho_\infty$.
А именно, положим
\begin{align*}
\rho_+ (\1, f)
&= \inf\int_0^1 \max_{x} F_t(x)\,dt
\intertext{и}
\rho_- (\1, f) 
&= \inf \int_0^1\big(-\min_{x} F_t(x)\big)\,dt.
\end{align*}
Тогда очевидно, что
\[\rho (\1, f) \ge \rho_+ (\1, f) + \rho_- (\1, f).\]
Однако во всех известных мне примерах выполняется равенство!
Было бы интересно доказать общий случай или найти контрпример.%
\footnote{%
  Наряду с многими другими интересными результатами о
  $\rho_{\pm}$, Макдафф \cite{McD00} показала, что $\rho \neq
  \rho_{+}+\rho_{-}$, если $(M,\Omega)$ — симплектическое раздутие
  в одной точке комплексной проективной плоскости, оснащённое
  определённой симплектической структурой.\dpp}
Заметим, что из \rindex{Витербо}\cite{V1} следует, что на $\Ham
(\RR^{2n})$ сумма $\rho_+ + \rho_-$  определяет биинвариантную
метрику. 
Насколько мне известно, в общем случае, для симплектических
многообразий похожих утверждений пока нет. 
\end{ex*}

\begin{ex*}{Соглашение}
В дальнейшем если не сказано противное, то мы используем обозначение \index[symb]{$\lVert\ \rVert$}\index[symb]{$\|\ \|$}$\|\ \|$ для $L_\infty$-нормы на $\A(M)$.
Через $\rho$ обозначим метрику $\rho_\infty$ и называем её \rindex{хоферовская метрика}\emph{хоферовской метрикой}.
Величину $\rho (\1, f)$ называем \rindex{хоферовская норма}\emph{хоферовской нормой} $f$.
Через $\length \{f_t\}$ обозначим длину гамильтонова пути $\{f_t\}$ относительно $L_\infty$-нормы (см. \ref{2.2}).
\end{ex*}

\chapter{Лагранжевы подмногообразия}\label{chap:3}

Цель этой и следующей глав — доказать невырожденность хоферовской метрики на $\RR^{2n}$.
Мы используем подход из \cite{P1}.
Для этого нам потребуется понятие лагранжевых подмногообразий в симплектических многообразиях.
Лагранжевы подмногообразия играют важную роль в симплектической топологии, а также в её приложениях в механике и вариационном исчислении.
Они будут многократно появляться в этой книге далее.

\section{Определения и примеры}\label{3.1}

\begin{ex*}{Определение}
Пусть $(M^{2n}, \Omega)$ — симплектическое многообразие.
Подмногообразие $L \subset M$ называется \rindex{лагранжево многообразие}\emph{лагранжевым}, если $\dim L \z= \tfrac12 \dim M = n$ и $\Omega|_{\T L} \equiv 0$.
Вложение (или погружение) $f\colon L^n \z\to M^{2n}$ называется лагранжевым, если $f^\ast \Omega \equiv 0$.
\end{ex*}

Перечислим некоторые важные примеры лагранжевых подмногообразий.

\begin{ex}[Кривые на поверхностях.]{}\label{3.1.A} 
Пусть $(M^2, \Omega)$ ориентированная поверхность с формой площади.
Тогда каждая кривая лагранжева (поскольку касательное пространство к кривой одномерно и $\Omega$ обнуляется на паре пропорциональных векторов).
\end{ex}

\begin{ex}[Расщеплённый тор.]{}\label{3.1.B} (ср. с \ref{1.1.C}).
Тор $S^1 \times\dots\times S^1
\z\subset 
\RR^2\times\z\dots\times \RR^2 
= 
\RR^{2n}$ лагранжев (2-форма $\Omega$ на $\RR^{2n}$ расщепляется).
\end{ex}

\begin{ex}[Графики 1-форм в кокасательных расслоениях.]{}\label{3.1.C}\ 
Этот пример играет важную роль в классической механике.
Пусть $N^n$ — произвольное многообразие.
Рассмотрим кокасательное расслоение $M \z= \T^\ast N$ с естественной
проекцией $\pi\colon \T^\ast N \to N$, $(p, q) \z\mapsto q$.
Определим 1-форму $\lambda$ на $M$ — так называемую \rindex{форма Лиувилля}\emph{форму Лиувилля}.
Для $q \in N$, $(p, q) \in \T^\ast N$ и $\xi \in \T_{(p, q)} \T^\ast
N$ положим $\lambda (\xi) \z= \langle p, \pi_\ast \xi \rangle$, где
$\langle\ ,\ \rangle$ — естественное спаривание
спаривание пространств $\T_qN$ и $T^\ast_qN$. 
Мы утверждаем, что $\Omega = d\lambda$ — симплектическая форма на $\T^\ast M$.
Используя локальные координаты $(p_1,\dots, p_n, q_1,\dots, q_n)$ на
$\T^\ast N$, положим $\xi = (\dot p_1 ,\dots, \dot p_n, \dot
q_1,\dots,\dot q_n)$; 
тогда $\pi_\ast \xi = (\dot q_1,\dots, \dot q_n)$.
В этих обозначениях, $\langle p, \pi_\ast \xi\rangle=\sum p_i \dot
q_i$, откуда следует, что $\lambda=\sum p_i dq_i$.
Следовательно, $\Omega = d\lambda =\sum dp_i \z\wedge dq_i$ — узнaём
стандартную симплектическую форму на $\RR^{2n}$ и утверждение
следует. 
\end{ex}

\begin{ex*}{Упражнение}\label{1-form-lagrange}
Покажите, что график 1-формы $\alpha$ на $N$ является лагранжевым
подмногообразием в $\T^\ast N$ тогда и только тогда, когда форма
$\alpha$ замкнута. 
\noindent\textit{Подсказка:} Докажите, что прообраз формы Лиувилля на
$\T^\ast M$ при естественной параметризации графика любой 1-формы
$\alpha$ на $M$ равен самой форме $\alpha$.
\end{ex*}

\begin{ex}[Симплектоморфизмы как лагранжевые подмногообразия.]{}\label{3.1.D}\\ 
Пусть $f\colon (M, \Omega) \to (M, \Omega)$ — диффеоморфизм.
Рассмотрим новое симплектическое многообразие $(M \times M, \Omega \oplus -\Omega)$.
Оставим как упражнение читателю показать, что график$(f) \subset (M \times M, \Omega \oplus -\Omega)$ является лагранжевым тогда и только тогда, когда $f$ — симплектоморфизм.
\end{ex}

\begin{ex}[Лагранжева надстройка.]{}\label{3.1.E}
\rindex{надстройка}
Пусть $L\subset(M,\Omega)$ — лагранжево подмногообразие.
Рассмотрим петлю гамильтоновых диффеоморфизмов $\{h_t\}$, $t\in S^1$,
$h_0=h_1=\1$, порождённую 1-периодическим гамильтонианом
$H_{t}(x)$. 
\end{ex}

\begin{thm*}{Предложение}
Пусть $M$ и $L$ те же, что выше, а $(r, t)$ — координаты на $\T^\ast S^1 = \RR \times S^1$.
Рассмотрим симплектическое многообразие $M \times \T^\ast S^1$ с симплектической формой $\sigma = \Omega + dr \wedge dt$.
Тогда
\[\phi\colon L \times S^1 \to M \times \T^\ast S^1,
\quad
(x, t) \mapsto \big(h_t (x), -H_{t} (h_t (x)), t\big)\]
— лагранжево вложение.
\end{thm*}

\parit{Доказательство.}
Достаточно доказать, что $\phi^\ast \sigma$ обращается в нуль на парах $(\xi, \xi')$ и на парах  $(\xi, \tfrac{\partial}{\partial t})$ при $\xi, \xi' \in \T L$ и $\tfrac{\partial}{\partial t} \in \T S^1$.
Вычислим 
\begin{align*}
\phi_\ast \xi
&= h_{t\ast} \xi - \langle dH_t, h_{t\ast} \xi\rangle
\tfrac{\partial}{\partial r},
\\
\phi_\ast \xi'
&= h_{t\ast} \xi' - \langle dH_t, h_{t\ast} \xi'\rangle
\tfrac{\partial}{\partial r}.
\end{align*}
Поскольку $L$ лагранжево, $\phi^\ast \sigma (\xi, \xi') = \Omega (h_{t\ast} \xi, h_{t\ast} \xi') = \Omega (\xi, \xi') = 0$.
Более того,
\begin{align*}
\phi_\ast\tfrac{\partial}{\partial t}
&= \sgrad H_t - 
\left( \langle dH_t, \sgrad H_t \rangle +\tfrac{\partial H}{\partial t}\right) \tfrac{\partial}{\partial r} +\tfrac{\partial}{\partial t}=
\\
&=\sgrad H_t - \tfrac{\partial H}{\partial t}\tfrac{\partial}{\partial r} +\tfrac{\partial}{\partial t},
\end{align*}
так что
\begin{align*}
\phi^\ast \Omega (\xi, \tfrac{\partial}{\partial t})
&= \Omega (h_{t\ast} \xi, \sgrad H_t) + dr \wedge dt (-\langle dH_t, h_{t\ast} \xi\rangle \tfrac{\partial}{\partial r},\tfrac{\partial}{\partial t}) =
\\
&=dH_t (h_{t\ast} \xi) - \langle dH_t, h_{t\ast} \xi\rangle = 0.
\end{align*}
\qedsf

\section[Класс Лиувилля]{Класс Лиувилля лагранжевых\\подмногообразий в $\bm{\RR^{2n}}$}
\label{3.2}
Пусть $L \subset (\RR^{2n}, dp \wedge dq)$ — лагранжево подмногообразие.
Рассмотрим сужение $\lambda|_{\T L}$ формы Лиувилля 
\[\lambda = p_1 dq_1+\dots+ p_n dq_n\]
на $L$.
Ясно, что $d (\lambda|_{\T L}) = \Omega|_{\T L} = 0$.
Класс когомологий \index[symb]{$\lambda_L$}$\lambda_L\in H^1(L;\RR)$ этой замкнутой 1-формы называется \rindex{класс Лиувилля}\emph{классом Лиувилля} лагранжева подмногообразия $L$.
Аналогично для лагранжева вложения или погружения $\phi\colon L \to \RR^{2n}$ класс Лиувилля определяется как $[\phi^\ast \lambda]$.
Класс Лиувилля лагранжева подмногообразия можно интерпретировать геометрически следующим образом.
Для 1-цикла $a\in H^1 (L)$ выберем 2-цепь $\Sigma$ в $\RR^{2n}$ такую, что $\partial\Sigma = a$.
Тогда
\[(\lambda_L, a) = \int_a\lambda_L = \int_\Sigma\Omega.\]
Это число не зависит от выбора $\Sigma$.
Из-за этой формулы значение $(\lambda_L, a)$ иногда называют \rindex{симплектическая площадь}\emph{симплектической площадью} класса $a$.
Обобщение этой конструкции на произвольное лагранжево многообразие $L$ симплектического многообразия $M$ даёт естественный гомоморфизм $H_2 (M, L; \ZZ) \to \RR$.
Важным свойством класса $\lambda_L$ является то, что он инвариантен относительно симплектоморфизмов $\RR^{2n}$,
то есть
$f^\ast \lambda_{f (L)} = \lambda_L$.

\begin{thm}[(\cite{G1})]{Теорема}\label{3.2.A}\rindex{Громов}
Пусть $L \subset \RR^{2n}$ — замкнутое лагранжево подмногообразие.
Тогда $\lambda_L \ne 0$.
\end{thm}

Подчеркнём, что $L$ вложено.
Для лагранжевых погружений это утверждение, вообще говоря, неверно.
В случае $n = 1$ это можно увидеть следующим образом.
Ясно, что любая замкнутая вложенная кривая ограничивает область положительной площади, но, например, погруженная восьмёрка может ограничивать нулевую площадь.
Это явление отражает «жёсткость лагранжевых вложений».

\begin{ex*}{Определение}
Замкнутое лагранжево подмногообразие $L \z\subset (\RR^{2n}, \omega)$
называется \rindex{рациональное подмногообразие}\emph{рациональным}, если $\lambda_L (H^1 (L; \ZZ))$ —
дискретная подгруппа в $\RR$.
В таком случае обозначим её положительную образующую через \index[symb]{$\gamma(L)$}$\gamma(L)$.
\end{ex*}

\begin{ex*}{Пример}
Расщеплённый тор $L = S^1 (r) \times\dots\times S^1 (r) \subset
\RR^{2n}$ рационален. 
Действительно, поскольку каждая окружность $S^1 (r)$ имеет
симплектическую площадь $\pi r^2$, получаем $\gamma (L) = \pi r^2$. 
Однако тор $S^1(1)\times S^1(\sqrt[3]{2}) \subset \RR^4$ не является
рациональным.
Симплектические площади двух окружностей равны $\pi$  и
$\sqrt[3]{4}\pi$ соответственно, и они порождают плотную подгруппу в
$\RR$.
\end{ex*}
\?{}{Лёне: Яша заметил, что теорема в статье Сикорава приписывается
  Громову (private communication), но статья действительно содержит и теорему и док-во.}
\begin{thm}[(\cite{S1})]{Теорема}%
  \hspace{-0.6em}\footnote{Смотри сноску на
    странице~\pageref{foot:sikorav}.}
  \label{3.2.B}\rindex{Сикорав}
  Пусть $L \subset B^2 (r) \times \RR^{2n - 2}$ — замкнутое
  рациональное лагранжево подмногообразие. 
  Тогда $\gamma (L) \le \pi r^2$.
\end{thm}

\begin{wrapfigure}[7]{o}{25 mm}
\vskip-3mm
\centering
\includegraphics{mppics/pic-1}
\caption{}\label{pic-1}
\vskip0mm
\end{wrapfigure}

Предположение, что $L$ вложено, необходимо.
На рис. \ref{pic-1} показано лагранжево погружение произвольной
симплектической площади.

Наш следующий результат даёт нижнюю оценку на \rindex{энергия смещения!лагранжева подмногообразия}\emph{энергию смещения} $\e (L)$ рационального лагранжева подмногообразия $L$ относительно хоферовской метрики.

\begin{thm}{Теорема}\label{3.2.C}
  Пусть $L \subset \RR^{2n}$ — замкнутое рациональное лагранжево
  подмногообразие.
  Тогда $\e (L) \ge \tfrac12 \gamma (L)$.
\end{thm}

Теоремы \ref{3.2.A} и \ref{3.2.B} доказаны в следующей главе.
Теорема \ref{3.2.C} следует из \ref{3.2.B} (см. раздел \ref{3.3}).
Выведем некоторые следствия из этих результатов.

\begin{ex}[Невырожденность хоферовской метрики.]{}\label{3.2.D}
Из теоремы \ref{3.2.C} следует, что хоферовская метрика на $\Ham (\RR^{2n}, \omega)$ невырождена.
Действительно, каждый шар 
$B^{2n}(r) = \{p_1^2 +\z\dots\z+ p_n^2 \z+ q_1^2+\z\dots\z+ q_n^2 \le r^2\}$
содержит рациональный расщеплённый тор 
\[
S^1(\tfrac r{\sqrt{n}}) \times\dots\times S^1(\tfrac r{\sqrt{n}})
=
\{p_1^{2}+q_1^2=\dots=p_n^{2}+q_n^2=\tfrac{r^2}{n}\}.
\]
Таким образом, $\e (B^{2n} (r)) \ge \tfrac{\pi r^2}{2n}> 0$, и, следуя рассуждению в \ref{sec:2.4}, получаем желаемую невырожденность метрики $\rho$.
Эта оценка не точна.
\rindex{Хофер}Хофер доказал \cite{H1}, что $\e (B^{2n} (r)) = \pi r^2$.
\end{ex}

\begin{ex}[Свойство несжимаемости.]{}\label{3.2.E}\rindex{несжимаемость}
Отметим, что $\gamma (L)$ — симплектический инвариант, то есть $\gamma (f (L)) = \gamma (L)$ для любого симплектоморфизма $f\colon \RR^{2n} \to \RR^{2n}$.
Таким образом, из теоремы \ref{3.2.B} следует теорема о несжимаемости \ref{1.1.C};
она утверждает, что расщеплённый тор с большим $\gamma (L) = \pi R^2$ 
нельзя поместить гамильтоновым диффеоморфизмом в $B^2 (r) \times \RR^{2n - 2}$ при $r<R$.
\end{ex}

\begin{ex}[Цилиндрическая симплектическая ёмкость.]{}\label{3.2.F}
Пусть $A \subset \RR^{2n}$ — ограниченное подмножество.
Положим \index[symb]{$c(A)$}
\[c(A) = \inf \set{\pi r^2}{ \exists\  g\colon \RR^{2n} \to \RR^{2n}},\]
где $g$ — такой симплектоморфизм, что $g (A) \subset B^2 (r) \times \RR^{2n - 2}$.
Эта функция, определённая на подмножествах $\RR^{2n}$, называется \rindex{цилиндрическая ёмкость}\emph{цилиндрической симплектической ёмкостью}.
На этом языке теорема \ref{3.2.B} читается так:
для замкнутого рационального лагранжева подмногообразия $c(L) \ge \gamma (L)$.
Эта ёмкость является симплектическим инвариантом и удовлетворяет следующему свойству монотонности:
$c (A) \z\le c (B)$ если $A \subset B$
(сравните это с похожим свойством монотонности энергии смещения, раздел \ref{sec:2.4}).
\end{ex}

\begin{ex}[Некоторые обобщения.]{}\label{3.2.G}
Пусть $(M, \Omega)$ — симплектическое многообразие.
Если $M$ открыто, то мы полагаем, что оно «хорошо» ведёт себя на бесконечности (этот класс включает, например, любое кокасательное расслоение со стандартной симплектической структурой, а также произведение кокасательного расслоения с любым замкнутым симплектическим многообразием). 
Возьмём лагранжево подмногообразие $L \subset M$ и рассмотрим
гомоморфизм $\lambda_L\colon \pi_2 (M, L) \z\to \RR$, переводящий любой
диск $\Sigma$ в $M$, граница которого лежит на $L$, в его
симплектическую площадь $\int_\Sigma \Omega$. 
Точно так же, как в случае $M = \RR^{2n}$, мы говорим, что $L$
рационально, если образ $\lambda_L$ дискретен; для рационального $L$
определим $\gamma (L)$ как положительную образующую образа
$\lambda_L$. 
Если $\lambda_L = 0$, то положим $\gamma (L) = + \infty$.
Наше доказательство теоремы \ref{3.2.C} без существенных изменений
распространяется на эту более общую постановку (см. \cite{P1}). 
А именно, мы получаем, что $\e (L) \ge \tfrac12 \gamma (L)$.
В качестве следствия%
\footnote{Оно упущено в \cite[с. 359]{P1}.}
можно получить следующее важное утверждение, доказанное Громовым в \cite{G1}:
$\e (L) = + \infty$ при $\lambda_L = 0$.
Его можно интерпретировать как свойство лагранжева пересечения: если $\lambda_L = 0$, то  для любого гамильтонова диффеоморфизма $\phi$ образ $\phi (L)$ пересекает $L$.
В главе \ref{chap:6} мы обсудим применение этого результата к хоферовской геометрии.

Отметим также, что эти оценки были значительно улучшены в
\rindex{Чеканов}\cite{Ch} при помощи гомологий Флоера (см. также
\cite{O3}). 
В частности, было показано, что каждое (не обязательно рациональное) замкнутое лагранжево подмногообразие $L \subset M$ имеет положительную энергию смещения.
\end{ex}

\begin{ex}[Изопериметрическое неравенство.]{}\rindex{изопериметрическое неравенство}
Мы завершаем этот раздел формулировкой следующего замечательного результата, принадлежащего \rindex{Витербо}Витербо \cite{V2}.
Пусть $L \subset \RR^{2n}$ — замкнутое лагранжево подмногообразие.
Тогда 
\[\e (L)^n \le 2^{\frac{n(n-3)}2} n^n V^2,\]
где $V$ обозначает $n$-мерный евклидов объем $L$.
Точную константу в этом неравенстве ещё предстоит найти. 
\end{ex}

\section{Оценка энергии смещения}\label{3.3}

В этом разделе мы выводим теорему \ref{3.2.C} из теоремы \ref{3.2.B}, используя элементарную геометрию.

\parbf{Шаг 1.}
Пусть $L$ — замкнутое рациональное лаграново подмногообразие и $h_t$, $t \in [0;1]$ — путь гамильтоновых диффеоморфизмов такой, что $h_0 = \1$ и $h_1 (L) \cap L = \emptyset$.
Выберем $\epsilon> 0$.
Не умаляя общности можно считать, что $h_t = \1$ при $t \in [0;\epsilon]$ и $h_t = h_1$ при $t \in [1 - \epsilon;1]$.
Этого можно добиться подходящей репараметризацией потока, сохраняющей его длину (используйте упражнение \ref{1.4.A}).
Пусть $H (x, t)$ — соответствующий гамильтониан.
Положим
\[l
=
\length\{h_t\} 
=
\int_0^1 (\max_x H_t - \min_x H_t)\,dt.\]
Нам нужно доказать, что $l \ge \tfrac12 \gamma (L)$.
Генеральный план в том, чтобы закодировать движение $L$ в потоке как замкнутое лагранжево подмногообразие в $\RR^{2n+2}$, а затем применить \ref{3.2.B}.
Мы воспользуемся лагранжевой надстройкой, описанной в п. \ref{3.1.E}.
Для этого нам понадобится петля гамильтоновых диффеоморфизмов.

\begin{figure}[ht!]
\vskip-0mm
\centering
\includegraphics{mppics/pic-2}
\caption{}\label{pic-2}
\vskip0mm
\end{figure}

\parbf{Шаг 2.}
Рассмотрим следующую петлю гамильтоновых диффеоморфизмов при $t \in [0;2]$
\begin{align*}
g_t
&=
\begin{cases}
\quad h_t&\text{при}\ t\in [0;1],
\\
\quad h_{2-t}&\text{при}\ t\in [1;2]
\end{cases}
\intertext{с гамильтонианом}
G (x, t)&=
\begin{cases}
\quad \phantom{-}H(x,t)&\text{при}\ t\in [0;1],
\\
\quad -H(x,2-t)&\text{при}\ t\in [1;2].
\end{cases}
\end{align*}

\begin{ex*}{Упражнение}
Покажите, что $\int_0^2 G(g_t(x),t)\,dt = 0$ при всех $x$.
\end{ex*}

Взяв лагранжеву надстройку петли $\{g_t\}$ (см. п. \ref{3.1.E}),
получаем новое лагранжево подмногообразие $L' \subset \RR^{2n} \z\times \T^{\ast} S^1$ как образ $L \times S^1$ при отображении 
\[(x, t) \mapsto \Big(g_t (x), -G \big(g_t (x), t\big), t\Big).\]
Напомним, что здесь $S^1 = \RR / 2\ZZ$.
Определим две функции 
\begin{align*}a_+ (t) &= - \min_x G (x, t) + \epsilon
&&\text{и}&
a_- (t) &= - \max_x G (x, t) - \epsilon.
\end{align*}
Ясно, что $L' \subset \RR^{2n} \times C \subset \RR^{2n} \times \T^\ast S^1$, где $C$ обозначает кольцо 
$\{a_- (t) \z<r \z<a_+ (t)\}$
(см. рис. 2).

\parbf{Шаг 3.}
Теперь мы хотим перейти от $\RR^{2n} \times \T^{\ast} S^1$ к $\RR^{2n} \times \RR^2$.
Рассмотрим специальное симплектическое погружение $\theta\: C \to \RR^2$
(этот трюк известен как громовская восьмерка, см. \cite{G1,AL}).

\begin{figure}[ht!]
\vskip-0mm
\centering
\includegraphics{mppics/pic-3}
\caption{}\label{pic-3}
\vskip0mm
\end{figure}

\begin{ex*}{Упражнение}
(см. рис. 3).
Покажите, что существует симплектическое погружение $\theta: C \to \RR^2 (p, q)$ со следующими свойствами: 
\begin{itemize}
\item $\theta$ переводит нулевое сечение $\{r = 0\}$ в восьмерку с равновеликими ушами, и, таким образом, замкнутая форма $\theta^\ast pdq - rdt$ точна (она замкнута, поскольку $\theta$ — симплектоморфизм).
\item $\theta$ — это вложение за пределами пары тонких горловин и оно склеивает эти горловины вместе.
\item площадь внутренних ушей произвольно мала, скажем, по $\epsilon$ у каждого. 
\end{itemize}
\end{ex*}

Заметим, что 
\begin{align*}
\area (C) 
&= \int_0^2(a_+ (t) - a_- (t))\,dt = 
\\
&=2 \int_0^1(\max_x H_t - \min H_t) \,dt + 4\epsilon =
\\
&=2l + 4\epsilon.
\end{align*}
Таким образом, образ $\theta (C)$ можно поместить в диск $B$
площади $2l \z+ 10\epsilon$ (мы добавили $10\epsilon$, чтобы остался зазор).

\parbf{Шаг 4.}
Теперь рассмотрим симплектическое погружение 
\[\theta' = \1 \times \theta: \RR^{2n} \times C \to \RR^{2n} \times \RR^2.\]
Очевидно, что $\theta' (L')$ — погруженное лагранжево подмногообразие лежащее в $\RR^{2n} \times B$. 
Покажем, что $L'' = \theta' (L')$ вложено.

Единственное место, где могут возникнуть двойные точки, — это тонкие горловины.
Но $g_t (L) = L$ при $t \in [-\epsilon;\epsilon]$ 
и $g_t (L) = h_1 (L)$ при $t \in [1 - \epsilon;1 + \epsilon]$, а по предположению $h_1 (L) \cap L = \emptyset$, поэтому двойных точек нет и $\theta'$ — вложение.

\parbf{Шаг 5.}
Нам осталось показать, что $L''$ рационально.
После этого мы сможем применить \ref{3.2.B}.
Докажем, что $\gamma (L) = \gamma (L'')$.
Пусть $\phi$ — композиция лагранжевой надстройки и $\theta'$,
то есть отображение 
\[\phi\colon 
L \times S^1
\to
\RR^{2n} \times \T^\ast S^1
\to
\RR^{2n} (p_1,\dots, p_n, q_1,\dots, q_n) \times \RR^2 (p, q)\]
отправляющее $(x, t)$ в $(g_t (x), \theta (-G (g_t (x), t), t))$.
Тогда $L''=\phi(L \times S^1)$.
Группа $H_1 (L'')$ порождается циклами вида $\phi (b)$, где либо $b \subset L \times {0}$, либо $b = {x_0} \times S^1$ при $x_0 \in L$.
В первом случае $\phi (b) = b \times \{\theta (0, 0)\}$, поэтому симплектические площади циклов $b$ и $\phi (b)$ совпадают.
Осталось рассмотреть второй случай.
Обозначим через $\alpha$ орбиту $\{g_t x_0\}$, $t \in [0; 2]$.
Тогда
\begin{align*}
\int_b\phi^\ast (p_1 dq_1 +\!\cdots\!
+ p_n dq_n + pdq)
&= \int_\alpha (p_1 dq_1 +\!\cdots\!
+ p_n dq_n) + \int_{\Gamma}\theta^\ast pdq =
\\
&= 0 + \int_{\Gamma} r\,dt =
\\
&= - \int_0^2G (g_t (x_0), t)\,dt 
= 0.
\end{align*}
где $\Gamma\subset \T^\ast S^{1}$  — график функции $t\mapsto -G_{t}(g_{t}x_{0})$, то есть проекция цикла $\phi(x_{0}\times S^{1})$ на $\T^\ast S^{1}$.
Это завершает доказательство того, что $L''$ — рациональное лагранжево подмногообразие с $\gamma (L) = \gamma (L'')$.
Напомним, что $L''$ содержится в $\RR^{2n} \times B$.
Принимая во внимание \ref{3.2.B}, получим, что
\[\gamma (L'') \le \area (B) = 2l + 10\epsilon\]
при всех $\epsilon> 0$.
Следовательно, $\e(L) \ge \tfrac12 \gamma (L)$.
\qeds

\chapter
[\texorpdfstring{$\bm{\bar\partial}$-уравнение}{∂-уравнение}]
{Лагранжевы граничные условия на $\bm{\bar\partial}$-уравнение}

В этой главе доказывается теорема~\ref{3.2.B}, которая утверждает, что $\gamma (L) \le \pi r^2$ для любого замкнутого рационального лагранжева подмногообразия $L \subset B^2 (r) \times \RR^{2n-2}$.
Доказательство основано на громовской технике псевдоголоморфных дисков.

\section[\texorpdfstring{Знакомство с $\bar\partial$-оператором}{Знакомство с ∂-оператором}]{Знакомство с $\bm{\bar\partial}$-оператором}\label{sec:4.1}

Отождествим $\RR^{2n} (p_1,q_1,\dots, p_n, q_n)$ с комплексным пространством
\[\CC^n (p_1 + iq_1,\dots, p_n + iq_n)
=
\CC^n (w_1,\dots, w_n)\]
и обозначим через $\langle\ ,\  \rangle$ евклидово скалярное произведение.
У нас появились три геометрические структуры: евклидова, симплектическая и комплексная.
Они связаны следующим образом
\[\langle \xi, \eta\rangle = \omega (\xi, i\eta).\]
Ограничимся проверкой этой формулы при $n = 1$.
Если $\xi \z= (p', q')$ и $\eta = (p'', q'')$, то
\[dp\wedge dq(\xi, i\eta)
=
dp\wedge dq\left(\binom{p'}{q'},\binom{-q''}{p''}\right)
=
p' p'' + q' q''
=
\langle\xi, \eta\rangle.\]
Далее мы измеряем площади и длины с помощью евклидовой метрики.
Рассмотрим единичный круг $D^2 \subset \CC$ с координатой $z \z= x \z+ iy$.
Пусть $f\: D^2 \to \CC^n$ гладкое отображение.
Определим $\bar\partial$-оператор $\bar\partial\: C^\infty (D^2, \CC^n) \to C^\infty  (D^2, \CC^n)$ как 
\[\bar\partial f=\frac12\left(\frac{\partial f}{\partial x} + i \frac{\partial f}{\partial y}\right).\]

\begin{ex*}{Пример}
Пусть $f\: \CC \to \CC$, $z \mapsto \bar z$.
Тогда $f(x,y)=x-iy$ и $\bar\partial f \z= \tfrac12(1+1)=1$.
Заметим, что $\bar\partial f  = \frac{\partial f}{\partial \bar z}$.
\end{ex*}

Введём пару полезных геометрических характеристик отображения $f\: D^2 \to \CC^n$:
\rindex{симплектическая площадь}\emph{симплектическая площадь} $f$, задаваемую как
\[\omega(f) =\int_{D^2}f^\ast\omega\]
и \rindex{евклидова площадь}\emph{евклидову площадь} $f$ определяемую как
\[\area(f)
=\int_{D^2}
\sqrt{
\left\langle\frac{\partial f}{\partial x},\frac{\partial f}{\partial x}\right\rangle
\left\langle \frac{\partial f}{\partial y},\frac{\partial f}{\partial y}\right\rangle
-
\left\langle\frac{\partial f}{\partial x},\frac{\partial f}{\partial y}\right\rangle^2
}
dxdy.
\]

\begin{thm}{Предложение}\label{4.1.A}
\begin{enumerate}[i)]
\item\label{4.1.A.i} $\displaystyle{\area(f)\le
  2\int_{D^2}|\bar\partial f|^2dxdy+\omega(f)}$
\item\label{4.1.A.ii} $\displaystyle{\area(f)\ge |\omega(f)|}$
\end{enumerate}
\end{thm}

\parit{Доказательство.}
Для $\xi, \eta \in \CC^n$ справедливо неравенство
\[
\sqrt{|\xi|^2 | \eta |^2 - \langle\xi, \eta\rangle^2}
\le
|\xi|{\cdot}|\eta| 
\le
\tfrac12(|\xi|^2 +|\eta|^2).
\]
Но
\[\tfrac12|\xi + i\eta|^2 + \omega (\xi, \eta)
=
\tfrac12(|\xi|^2 + |\eta|^2) + \langle\xi, i\eta\rangle + \langle\xi, -i\eta\rangle
=
\tfrac12(| \xi |^2 + | \eta |^2 ),\]
и поэтому 
\[\sqrt{| \xi |^2 | \eta |^2 - \langle\xi, \eta\rangle^2}
\le
\tfrac12| \xi + i\eta |^2 + \omega (\xi, \eta).
\]

Интегрируя полученное поточечное неравенство, получаем 
\ref{4.1.A}.\ref{4.1.A.i}.

Чтобы доказать \ref{4.1.A}.\ref{4.1.A.ii}, нам нужно снова проверить поточечное неравенство.
Предположим, что $\eta \ne 0$, и заметим, что $\langle\eta, i\eta\rangle = 0$.
Проецируя $\xi$ на $\eta$ и $i\eta$, получаем 
\[\left\langle\xi, \frac{\eta}{|\eta|}\right\rangle^2
+
\left\langle\xi, \frac{i\eta}{|i\eta|} \right\rangle^2 \le | \xi |^2.\]
Поскольку $| \eta | = | i\eta |$, это неравенство читается как
\[\langle\xi, \eta\rangle^2 + \omega (\xi, \eta)^2 \le | \xi |^2 | \eta |^2,\]
следовательно, 
\[
|\omega(\xi,\eta)|
\le \sqrt{ | \xi |^2 | \eta |^2 - \langle\xi, \eta\rangle^2}.
\qedsin
\]

\section{Краевая задача}\label{sec:4.2}

Пусть $L \subset \CC^n$ — замкнутое лагранжево подмногообразие и $g\: D^2 \z\times \CC^n \to \CC^n$ — гладкое отображение, ограниченное вместе со всеми своими производными.
Выберем класс $\alpha \in H_2 (\CC^n, L)$ и рассмотрим следующую краевую задачу.

\emph{Найти гладкое отображение $f\: (D^2, \partial D^2) \to (\CC^n, L)$ такое, что} 
\[
\begin{cases}
\quad\bar\partial f(z) = g (z, f (z)),
\\
\quad[f] = \alpha.
\end{cases}
\eqno{(P \big(\alpha, g)\big)}
\]

\begin{ex*}{Пример}
Если $g = 0$ и $\alpha = 0$, то пространство решений уравнения $P (0, 0)$ состоит из постоянных отображений $f (z) \equiv w$ при $w \in L$.
Чтобы в этом убедиться, заметим, что $\omega (f) = 0$.
В самом деле, поскольку $\alpha = 0$ и $L$ лагранжево, кривая $f (\partial D^2)$ ограничивает 2-цепь в $L$ с нулевой симплектической площадью.
Эта цепь вместе с $f (D^2)$ образует замкнутую поверхность в $\CC^n$.
В силу точности формы $\omega$, симплектическая площадь этой поверхности равна нулю.
Значит $\omega (f) = 0$.
Далее, поскольку $g = 0$, получаем $\bar\partial f=0$.
Итак, первая часть \ref{4.1.A} влечёт, что $\area (f) = 0$ и, значит, $\tfrac{\partial f}{\partial x}$ и $\tfrac{\partial f}{\partial y}$ параллельны.
С другой стороны, $\tfrac{\partial f}{\partial x}=-i\tfrac{\partial f}{\partial y}$,
следовательно, $\tfrac{\partial f}{\partial x}\perp\tfrac{\partial f}{\partial y}$.
Отсюда $\tfrac{\partial f}{\partial x}=\tfrac{\partial f}{\partial y}=0$.
Итак, $f$ является постоянным отображением.
Учтя граничные условия, получаем, что образ $f$ лежит в $L$.
\end{ex*}

Предположим теперь, что у нас есть последовательность функций
$\{g_n\}$, которая $C^\infty$-сходится к некоторой функции $g$. 
Пусть $f_n$ — решения соответствующих задач $P(\alpha, g_n)$.
Знаменитая \rindex{Громов}\rindex{теорема о компактности}\emph{теорема Громова о компактности} (см. \cite{G1,AL})
утверждает, что либо $\{f_n\}$ содержит подпоследовательность,
сходящуюся к решению $P (\alpha, g)$, либо происходит выдувание.
Чтобы объяснить, что такое выдувание, мы вводим понятие составного решения. 

\begin{ex*}{Определение}
Рассмотрим следующие данные:
\begin{itemize}
\item Разложение $\alpha = \alpha' + \beta_1 +\dots + \beta_k$, где $\beta_j \ne 0$, $j = 1,\dots,k$;
\item Решение $f$ уравнения $P (\alpha', g)$;
\item Решения $h_j$ уравнения  $P (\beta_j, 0)$, это так называемые псевдоголоморфные диски.%
\footnote{В нашем случае $g=0$ и комплексная структура стандартная, а значит, диски будут голоморфными. — \textit{Прим. ред.}}
\end{itemize}
Этот объект называется \rindex{составное решение}\emph{составным
  решением} уравнения $P(\alpha,g)$, и 
$f(D^2)\z\cup h_1(D^2) \z\cup\dots\cup h_k (D^2)$ называется его образом.
\end{ex*}

Мы говорим, что просходит выдувание, если существует подпоследовательность $\{f_n\}$ (которую мы снова обозначим через $\{f_n\}$), которая сходится к составному решению $P(\alpha,g)$.
Нам будет важно одно свойство этой сходимости — непрерывность
евклидовой площади
\[\area (f_n)
\to
\area(f)
+\sum_{j=1}^k\area (h_j).\]
Mы не приводим точное определение сходимости, оно довольно сложное (см. \cite{G1,AL}).
Иллюстративный пример дан в разделе~\ref{sec:4.4}.

Используя теорему компактности, М. Громов установил следующий важный результат \cite{G1}.

\begin{thm*}{Принцип продолжения}\rindex{принцип продолжения}
Рассмотрим семейство «общего положения» $g_s (z, w)$, $s \in [0;1]$ с
$g_0 = 0$. 
Тогда либо $P (0, g_s)$ имеет решение при всех $s$, либо при некотором $s_\infty \le 1$ происходит выдувание, то есть существует
 такая подпоследовательность  $s_j \to s_\infty$, что
последовательность решений $P(0,g_{s_j})$ сходится к составному решению $P
(0, g_{s_\infty})$. 
\end{thm*}

Понятие «общего положения» следует толковать следующим образом.
Пространство всех семейств $g_s$ можно наделить подходящей структурой
банахова многообразия. 
Семейства общего положения образуют плотное G-дельта-множество (то есть счётное пересечение открытых и плотных подмножеств) в этом пространстве.
В частности, каждое семейство $g_s$ переходит в общее положение после
сколь угодно малого возмущения. 
Мы отсылаем к \cite{G1,AL} за подробностями.

\section{Приложение класса Лиувилля}

Мы приведём доказательство теоремы \ref{3.2.B}, данное \rindex{Сикорав}Сикоравым \cite{S1}.
Предположим, что $L \subset B^2 (r) \times \CC^{n-1}$ — замкнутое
лагранжево подмногообразие. 
Возьмём $g (z, w) = (\sigma, 0 ,\dots, 0) \in \CC^n$ при некотором $\sigma \in \CC$.

\begin{thm}{Лемма}\label{4.3.A}
Если $| \sigma | > r$, то $P (0, g)$ не имеет решений.
\end{thm}

\parit{Доказательство.}
Предположим, что $f$ — решение.
Обозначим через $\phi$ его первую (комплексную) координату.
Тогда 
\[\frac{\partial\phi}{\partial x}+i\frac{\partial\phi}{\partial y} = 2\sigma.\]
Поскольку $L \subset B^2 (r) \times \CC^{n-1}$, получаем $\bigl|\phi|_{\partial D^2}\bigr|\le r$.
Далее
\begin{align*}
2\pi\sigma &= \int_{D^2}\left(\frac{\partial\phi}{\partial x}+i\frac{\partial\phi}{\partial y}\right)dxdy =
\\
&=\int_{D^2} d (\phi dy - i\phi dx)  = 
\\
&=\int_{S^1}\phi dy - i\phi dx.
\end{align*}
Далее $x + iy \z= e^{2\pi i t}$ и $dx + idy \z= 2\pi i e^{2\pi it} dt$, поэтому $dy - idx \z= 2\pi e^{2\pi i t} dt$.
Следовательно,
\[2\pi | \sigma | 
= 
2\pi\left|\int_0^1  e^{\pi i t} \phi(e^{2\pi i t})\,dt\right|
\le
2\pi r \]
и, следовательно, $| \sigma | \le r$.
\qeds

Возьмём теперь любое $\sigma$ такое, что $|\sigma|>r$, и применим принцип продолжения к семейству $g_s = (s\sigma, 0 ,\dots, 0)$, $s \in [0;1]$.
Предыдущая лемма говорит нам, что не существует решения при $s = 1$, поэтому для слегка возмущённого $g_s$ происходит выдувание.
Для простоты будем считать, что выдувание происходит на самом семействе $g_s$.
В общем случае рассуждения остаются без изменений (только надо оценивать с точностью до $\epsilon$), проверка предоставляeтся читателю.

Итак, у нас есть последовательность $s_n \to s_\infty \le 1$ и разложение $0 = \alpha + \beta_1 +\dots+ \beta_k$, $\beta_j \ne 0$.
Пусть $f_n$ — решения $P (0, g_{s_n})$, а $f_\infty$ — решение $P
(\alpha, g_{s_\infty})$ с голоморфными дисками $h_1,\dots,h_k$ такими, что $[h_j] \z= \beta_j$, и
\[\area (f_n)
\to 
\area (f_\infty) + \sum_{j=1}^k\area (h_j).\]
Применяя обе части \ref{4.1.A} и используя, что диски $h_j$
голоморфны, получаем
\[\area (h_j) = \omega (h_j) \ge \gamma (L).\]
Неравенство следует из того, что $[h_j] = \beta_j \ne 0$.
Согласно \ref{4.1.A}.\ref{4.1.A.ii},
\[\area (f_\infty)
\ge
|\omega(f_\infty)|
=
\left|\sum\omega (h_j)\right|
\ge
\gamma(L).\]
Таким образом, $\area (f_\infty) + \sum \area (h_j) \ge 2\gamma (L)$.
С другой стороны, \ref{4.1.A}.\ref{4.1.A.i} влечёт, что 
\[\area(f_n)
\le
2\pi s^2_n|\sigma|^2
\le
2\pi|\sigma|^2.
\]
Здесь мы пользуемся тем, что $\omega (f_n) = 0$ (поскольку $[f_n] = 0$) и $\bar\partial f_n\z=g_{s_n}$.
Вместе эти два неравенства дают $2\pi | \sigma | \ge 2\gamma (L)$, что выполняется при всех $\sigma$ с $| \sigma | > r$.
Поэтому $\pi r^2 \ge \gamma (L)$, и теорема доказана.
\qeds

\parit{Доказательство \ref{3.2.A}.}
Рассмотрим замкнутое лагранжево подмногообразие $L \subset B^2 (r) \times \CC^{n-1}$.
Согласно \ref{4.3.A} задача
\[
\begin{cases}
\quad\bar\partial f(z)=(s\sigma,0,\dots,0),&|\sigma|>r,
\\
\quad[f]=0
\end{cases}
\]
не имеет решения при $s = 1$.
По принципу продолжения произошло выдувание.
Значит, существует ненулевой класс $\beta_1$, который представлен голоморфным диском $h_1$.
Поскольку $h_1 \z\ne \const$, получаем $\omega (h_1)> 0$.
Так мы нашли диск в $\CC^n$, натянутый на $h_1 (\partial D^2)$, который имеет ненулевую симплектическую площадь.
Отсюда $\lambda_L \ne 0$.
\qeds

\section{Пример}\label{sec:4.4}

В этом разделе мы разберём конкретный пример, иллюстрирующий процессы в предыдущем доказательстве.
Пусть $L = \partial D^2 \subset \CC$, и пусть $\sigma = 1$.
Требуется найти все отображения $f\: D^2 \to \CC$ такие, что $f
(\partial D^2) \subset \partial D^2$ и 
\begin{equation}
\begin{cases}
\quad\bar\partial f(z,\bar z)=s,
\\
\quad[f|_{\partial D^2}]=0.
\end{cases}
\label{eq:4.4.A}
\end{equation}
Так как $\tfrac{\partial f}{\partial \bar z} = s$ получаем, что $f (z,
\bar z) = s\bar z + u (z)$ для некоторой голоморфной функции $u$ на
$D^2$. 
Мы утверждаем, что $s+ zu (z)$ голоморфная функция, отображающая
$\partial D^2$ в $\partial D^2$, и что $(s+ zu (z))|_{\partial D^2}$
имеет степень 1. 
Действительно, 
\[zf (z,\bar z) = s | z |^2 + zu (z),\]
так что
\[|z|\cdot| f (z, \bar z) | = \left| s \cdot\left| z \right|^{2} + zu (z)\right|.\]
При $|z|=1$ это выражение можно переписать как
$|f(z,\bar z)|=|s\z+zu(z)|$.
Поскольку $|f(z,\bar z)|=1$, получаем, что $s+zu(z)$ — голоморфная функция, перводящая $\partial D^2$ в $\partial D^2$.
Заметим, что $\deg f = 0$ и $\deg z = 1$, так что $\deg zf= 1 $ и, следовательно, $\deg (s + zu (z)) = 1$.
Все такие голоморфные функции описывают изометрию гиперболической метрики в круге.
Их можно записать как
\[e^{i\theta}\frac{1 - \bar\alpha z}{z-\alpha}\]
при $\theta \in \RR$ и $| \alpha | > 1$.
Таким образом, $s + zu (z) = e^{i\theta}\frac{1 - \bar\alpha z}{z-\alpha}$ и 
\[zu (z)
=
\frac{e^{i\theta} + \alpha s - z (s + e^{i\theta} \bar\alpha)}{z-\alpha}.\]
Поскольку $u$ голоморфна, у неё нет полюсов, так что $e^{i\theta} +
\alpha s = 0$ и, следовательно, $\alpha =-\frac{e^{i\theta}}{s}$.
Теперь  мы видим, что $1 <| \alpha | = | \tfrac1s |$, это влечёт, что $s<1$ и что у уравнения \ref{eq:4.4.A} нет решений при $s \ge 1$.
Значит, при $s = 1$ происходит выдувание.
Для простоты положим $\theta = 0$.
Тогда 
\[u(z)
=
-\frac{s-\frac1s}{z+\frac1s}
=
\frac{1-s^2}{sz+1}\]
и $f_s(z,\bar z)=s\bar z+\frac{1-s^2}{sz+1}$.
Для любого $z \ne -1$ имеем $f_s (z, \bar z) \to \bar z$ при $s\to1$, 
более того, эта сходимость равномерна вне произвольной окрестности $-1$.
Рассмотрим графики $f_s$ в $D^2 \times \CC \subset \CC \times
  \CC$.
Положим $w \z= f_s (z, \bar z)$, так что 
\[(w - s\bar z) (sz + 1) = 1 - s^2.\]
При $s \to 1$ это уравнение переходит в $(w - \bar z) (z + 1) = 0$, а его график приближается к объединению двух кривых 
\[
\left[
\begin{aligned}
\quad w&=\bar z,
\\
\quad z&=-1.
\end{aligned}
\right.
\]

\begin{figure}[th!]
\vskip-0mm
\centering
\includegraphics{mppics/pic-4}
\caption{}\label{pic-4}
\vskip0mm
\end{figure}

Здесь $w = \bar z$ является графиком $f_\infty$ и $\{z = -1\}$ соответствует голоморфному диску с краем на $\{-1\} \times L$.
Спроецировав предельную кривую на $w$-координату, получим объединение двух дисков
которые \?{выдулись}{Яше не понравилось. По-моему, нормально.} из нашего семейства.
Действительно, $f_\infty (z) = \bar z$ является решением $P (-a, 1)$, где $a = [S^1]$ и голоморфный диск $h (z) = z$ является решением $P (a, 0)$.

Увидеть выдувание можно посмотрев на вещественную часть уравнений.
Если рассмотреть графики соответствующих функций
$f_s(x)\z=sx+\frac{1-s^2}{sx+1}$ при $x \in [-1;1]$, 
то получим картинку на рис.~\ref{pic-4}. 
Графики $f_s$ сходятся к объединению двух кривых: график вещественной
части $f_\infty$ и отрезок $I = [-1;1]$, который является
вещественной частью голоморфного диска $\{-1\} \times D^2$.  

\chapter{Линеаризация хоферовской геометрии}

В этой главе мы будем думать про $\rho (\1, \phi)$ как про расстояние между точкой и подмножеством в линейном нормированном пространстве.
Позже это позволит нам получить оценки снизу на $\rho (\1, \phi)$ в некоторых интересных случаях.

\section[Периодические гамильтонианы]{Пространство периодических\\гамильтонианов}

Обозначим через \index[symb]{$\F$}$\F$ пространство всех гладких нормализованных
гамильтонинанов $F\: M \times \RR \to \RR$, 1-периодических по
времени: $F (x, t \z+ 1) = F (x, t)$ при всех $x \in M$ и $t\in\RR$.
Такие $F$ мы часто будем рассматривать как функции на $M \times S^1$,
где $S^1 = \RR / \ZZ$. 
Для гамильтониана $F \in \F$ обозначим через $\phi_F$ отображение
$f_1$ соответствующего гамильтонова потока $\{f_t\}$. 
Заметим, что каждый гамильтонов диффеоморфизм $\phi$ представим таким образом. 
В самом деле, пусть $\{g_t\}$ — любой поток с $g_1 = \phi$.
Выберем функцию $a\: [0; 1] \to [0; 1]$ такую, что $a \equiv 0$ в окрестности нуля и $a \equiv 1$ в окрестности единицы.
Рассмотрим новый поток $f_t = g_{a(t)}$ и 
продолжим его на всё $\RR$ по формуле $f_{t+1} = f_t f_1$.
Ясно, что поток получится гладким.
Утверждение следует из следующего упражнения.

\begin{ex}{Упражнение}\label{5.1.A}
Докажите, что гамильтонов поток $\{f_t\}$, $t \z\in \RR$,
порождается гамильтонианом из $\F$ тогда и только тогда, когда
$f_{t+1} = f_t f_1$ при всех $t$. 
\end{ex}

Рассмотрим подмножество \index[symb]{$\H$}$\H \subset \F$, определенное как 
\[\H = \set{H \in \F}{\phi_H = \1}.\]
Иными словами, гамильтонианы из $\H$ порождают петли гамильтоновых диффеоморфизмов (или \rindex{гамильтонова петля}\emph{гамильтоновы петли}).
Определим норму на $\F$ \index[symb]{$\VERT F \VERT$}
\[\VERT F \VERT = \max_{t} \| F_t \| = \max_{t} (\max_{x} F (x, t) -
\min_{x} F (x, t)).\]
Теперь мы можем сформулировать основную теорему этой главы.

\begin{thm}{Теорема}\label{5.1.B}
\[\rho (\1, \phi_F) = \inf_{H\in\H} \VERT F - H \VERT\]
при любом $F \in \F$ .
\end{thm}

\begin{figure}[ht!]
\vskip0mm
\centering
\includegraphics{mppics/pic-5}
\caption{}\label{pic-5}
\vskip0mm
\end{figure}

Обратите внимание, что правая часть — это просто расстояние от $F$ до $\H$ в смысле нашей нормы (см. рис.~5).
Таким образом, множество $\H$ многое помнит о хоферовской геометрии.
В следующих главах мы установим некоторые интересные свойства $\H$ и внимательно посмотрим на гамильтоновы петли.

Теорема~\ref{5.1.B} — простое следствие следующего утверждения.

\begin{thm}{Лемма}\label{5.1.C}
При любом $\phi \in \Ham (M)$ выполнено равенство
\[\rho (\1, \phi) = \inf \VERT F \VERT,\]
где нижняя грань берётся по всем гамильтонианам $F\in\F$, порождающим $\phi$.
\end{thm}

Используя терминологию некоторых работ автора, можно сказать, что «грубая» хоферовская норма совпадает с обычной.

\parit{Доказательство \ref{5.1.C}.}
Для $\phi \in \Ham (M, \Omega)$ положим $r (\1, \phi) \z= \inf \VERT F
\VERT$ где $F$ пробегает все гамильтонианы $F \in \F$,
порождающие~$\phi$. 
Ясно, что $r (\1, \phi) \ge \rho (\1, \phi)$.
Остаётся доказать обратное неравенство. 
Зафиксируем положительное число $\epsilon$.
Выберем путь $\{f_t\}$, $t \in [0; 1]$ гамильтоновых диффеоморфизмов
таких, что $f_0 = \1$, $f_1 = \phi$ и  $\int_0^1 m (t)\,dt \le \rho
(\1, \phi) + \epsilon$, где $m (t) = \| F_t \|$.

Не умаляя общности, можно считать, что $F \in \F$ и $m (t)> 0$ при всех $t$.
В самом деле, чтобы гарантировать периодичность, можно провести
репараметризацию времени, как же как в начале раздела. 
Обоснование
предположения о положительности функции $m$
даётся в следующем разделе.
Обозначим через $\C$ пространство всех $C^1$-гладких диффеоморфизмов окружности
$S^1$, сохраняющих ориентацию и фиксирующих $0$. 
Отметим, что для $a \in \C$ путь $f_a = \{f_{a(t)}\}$ порождается
нормализованным гамильтонианом $F^a (x, t) = a' (t) F (x, a(t))$, где
$a'$ обозначает производную по $t$ (см. \ref{1.4.A}). 
Пусть $a(t)$ --- обратная функция к 
\[b(t)
=
\frac{\int\limits_0^t m(s)ds}{\int\limits_0^1 m(s)ds}.\]
Обратите внимание, что 
\[\VERT F \VERT = \max a' (t) m (a (t)) = \max (m (t) / b'(t)) = \int_0^1m (t)\,dt.\]
Заключаем, что $\VERT F^a \VERT \le \rho (\1, \phi) + \epsilon$.
Приближая $a$ в $C^1$-топологии гладким диффеоморфизмом из $\C$, мы видим, что можно найти гладкий нормализованный гамильтониан, скажем $\tilde F$, который порождает $\phi$ и удовлетворяет условию $\VERT \tilde F \VERT \le \rho (\1, \phi) + 2\epsilon$.
Поскольку это можно сделать для произвольного $\epsilon$, заключаем, что $r (\1, \phi) \le \rho (\1, \phi)$, что завершает доказательство.
\qeds

\parit{Доказательство \ref{5.1.B}.}
Будем обозначать через $\{f_t\}$ гамильтонов поток, порожденный $F$.
Пусть $\{g_t\}$ — любой другой гамильтонов поток, порожденный $G \in \F$ с $g_1 = \phi_F$.
Разложим $g_t$ как $h_t \circ f_t$.
Как следует из \ref{5.1.A}, $\{h_t\}$ является петлей гамильтоновых диффеоморфизмов, то есть его нормализованный гамильтониан $H$ принадлежит $\F$ и $h_0 = h_1 = \1$.
Наоборот, для каждой петли $\{h_t\}$ поток $\{h_t \circ f_t\}$ порождается гамильтонианом из $\F$,%
\footnote{Вообще говоря, поток $\{f_t \circ h_t\}$ (порядок важен) {}\emph{не} порождается периодическим гамильтонианом!}
и в единичное временя он равен~$\phi_F$.
Далее, 
\[G (x, t) = H (x, t) + F (h^{-1}_t x, t).\]
Пусть $H' (x, t) = -H (h_t x, t)$.
Обратите внимание, что $H'$ порождает петлю $\{h^{-1}_t\}$ и, следовательно, $H'\in\H$.
С другой стороны, из приведённого выше выражения для $G$ следует, что $\VERT G \VERT = \VERT F - H' \VERT$.
Ввиду вышесказанного каждому $G$ соответствует единственный $H'$ и наоборот.
Таким образом, требуемое утверждение немедленно следует из леммы \ref{5.1.C}.
\qeds

\section{Регуляризация}\label{5.2}

В этом разделе мы обоснуем предположение $m(t)>0$ в доказательстве леммы \ref{5.1.C}.
Поток $\{f_t\}$ называется регулярным, если при любом $t$ его нормализованный гамильтониан $F_t$ не обращается тождественно в ноль.
Иными словами, если вектор, касательный к пути $\{f_t\}$, не
обращается в ноль ни при каком $t$.

\begin{thm}{Предложение}\label{5.2.A}
Пусть $\{f_t\}$ — поток, порожденный гамильтонианом из $\F$.
Тогда существует произвольно малая (в $C^\infty$-смысле) петля $\{h_t\}$ такая, что поток $\{h^{-1}_t f_t\}$ регулярен.
\end{thm}

Доказательство разбито на несколько шагов.

1) Сначала поймём, что собственно надо доказывать.
Предположим, что $\{h_t\}$ — петля, порожденная гамильтонианом $H \in \H$.
Тогда гамильтониан потока $\{h^{-1}_t f_t\}$ задаётся выражением $-H (h_t x, t) \z+ F (h_t x, t)$.
Нам надо доказать, что при любом $t$ это выражение не обращается в нуль тождественно.
Иначе говоря,
\begin{equation}
F (x, t) - H (x, t) \not\equiv 0\label{eq:5.2.B}
\end{equation}
при всех $t$.
Итак, нам надо построить сколь угодно малый гамильтониан $H \in \H$, удовлетворяющий \ref{eq:5.2.B}

2) Введём ещё одно полезное понятие:
\emph{$k$-мерной вариации постоянной петли} — это гладкое семейство петель $\{h_t (\epsilon)\}$, где $\epsilon$ принадлежит окрестности $0$ в $\RR^k$ и $h_t (0) = \1$ при всех $t$.
Если $M$ открыто, то мы дополнительно требуем, чтобы носители всех $h_t (\epsilon)$ лежали в некотором компактном подмножестве многообразия $M$.

Вот удобный способ создания вариаций.
Начнём с однопараметрического случая.
Выберем гамильтониан $G \in \F$ такой, что 
\begin{equation}
\int_0^1 G (x, t)\,dt = 0 
\label{eq:5.2.C}
\end{equation}
при любом $x \in M$.
Затем определим $h_t (\epsilon) \in \Ham (M, \Omega)$ как гамильтонов поток в момент $\epsilon$, порожденный независящим от времени гамильтонианом $\int_0^t G(x,s) ds$.

\begin{thm}{Упражнение}\label{5.2.D}
Пусть $H (x, t, \epsilon)$ — нормализованный гамильтониан петли $\{h_t (\epsilon)\}$.
Докажите, что 
\[\frac{\partial}{\partial \epsilon}|_{\epsilon=0} H (x, t, \epsilon) = G (x, t).\]
\end{thm}

Для построения $k$-мерной вариации естественно взять композицию одномерных: 
\[
h_t (\epsilon_1 ,\dots, \epsilon_k)
=
h_t^{(1)} (\epsilon_1) \circ\dots
\circ h_t^{(k)} (\epsilon_k).
\]

Каждый $h^{(j)}$ строится по функции $G^{(j)}$, как указано выше.
Упражнение~\ref{5.2.D} говорит, что частная производная гамильтониана $H (x, t, \epsilon)$ по $\epsilon_j$ при $\epsilon = 0$ равна $G^{(j)}$.

3) Выберем точку $y\in M$ и рассмотрим $2n$-мерное линейное пространство $E = \T_y^\ast M$.
Выберем $2n$ гладких замкнутых кривых $\alpha_1 (t),\z\dots, \alpha_{2n} (t)$ (где $t \in S^1$), удовлетворяющие следующим условиям:
\begin{itemize}
\item $\int_0^1 \alpha_j (t)\,dt = 0$ при всех $j = 1,\dots, 2n$; 
\item векторы $\alpha_1 (t),\dots, \alpha_{2n} (t)$ линейно независимы при любом $t$.
\end{itemize}
Вот конструкция такой системы кривых.
Выбираем базис $u_1, v_1,\z\dots, u_n, v_n$ в $E$ и возьмём кривые вида $u_j\cos 2\pi t + v_j\sin 2\pi t$ и $-u_j\sin 2\pi t \z+ v_j\cos 2\pi t$.

4) Теперь выберем функции $G_1(x,t),\dots, G_{2n}(x,t)$ в $\F$, удовлетворяет условию \ref{eq:5.2.C} и такие, что $d_y G_t^{(j)} = \alpha_j (t)$.
Рассмотрим соответствующую $2n$-мерную вариацию $\{h_t (\epsilon)\}$ постоянной петли, как на шаге 2.
Рассмотрим отображение $\Phi\: S^1 \times \RR^{2n} (\epsilon_1 ,\dots, \epsilon_{2n}) \to E$, определенное как 
\[(t, \epsilon) \to d_y \big(F_t - H_t (\epsilon)\big).\]
Заметим, что $\Phi$ — субмерсия в некоторой окрестности $U$ окружности $\{\epsilon = 0\}$.
Действительно, наша конструкция вместе с обсуждением в шаге 2 влечёт, что 
\[\tfrac{\partial}{\partial\epsilon}|_{\epsilon = 0} \Phi (t, \epsilon) = \alpha_j (t),\]
а эти векторы порождают всё $E$.
Обозначим через $\Psi$ сужение $\Phi$ на $S^1 \times U$.
Поскольку $\Psi$ является субмерсией, множество $\Psi^{-1} (0)$ является одномерным подмногообразием в $S^1 \times U$, поэтому его проекция на $U$ нигде не плотна.
Таким образом, существуют произвольные малые значения параметра $\epsilon$, что $d_y (F_t - H_t (\epsilon)) \ne 0$ при всех~$t$.
Следовательно, для каждого $t$ выполняется условие \ref{eq:5.2.B} — конец доказательства.
\qeds

\section{Пути в данном гомотопическом классе}\label{5.3}

Мы будем понимать \rindex{гомотопия}\emph{гомотопию} как гладкое однопараметрическое семейство путей.
Если не сказано обратное, то мы рассматриваем гомотопии незамкнутых путей с фиксированными конечными точками, а также гомотопии петель с базовой точкой $\1$.
В случае открытого многообразия $M$ обычно предполагается, что носители всех диффеоморфизмов, входящих в двупараметрические семейства, содержатся в компактном подмножестве объемлющего многообразия.

Выберем гамильтониан $F \in \F$ и обозначим через $\{f_t\}$ соответствующий гамильтонов поток.
Рассмотрим величину \index[symb]{$\l(F)$}
\[\l(F) = \inf \length \{g_t\},\]
где нижняя грань берётся по всем гамильтоновым путям $\{g_t\}$, $t
\z\in [0; 1]$ с $g_0 = \1$, $g_1 = \phi_F$, которые гомотопны
$\{f_t\}$ с фиксированными концами. 
Эту величину полезно интерпретировать следующим образом.
Рассмотрим универсальное накрытие $\widetilde\Ham(M, \Omega)$
пространства $(\Ham (M, \Omega), \1)$. 
Оно определяется обычным способом с небольшим исключением —
рассматриваются только гладкие пути и гладкие гомотопии. 
Финслерова структура на $\Ham (M, \Omega)$ канонически поднимается в
универсальное накрытие. 
Таким образом возникает понятие длины гладкой кривой, а значит, и
метрика $\tilde\rho$ на  $\widetilde\Ham(M, \Omega)$. 

Обозначим через $\tilde\1$ каноническое поднятие $\1$ на
$\widetilde\Ham(M, \Omega)$, а через 
$\tilde\phi_F$ — поднятие $\phi_F$, определяемое путём $\{f_t\}$, $t
\in [0; 1]$. 
На этом языке $l (F) = \tilde\rho (\tilde\1, \tilde\phi_F)$, то есть
эта величина отвечает за геометрию универсального накрытия. 

Обозначим через $\H_c$ множество всех гамильтонианов из $\H$,
порождающих стягиваемые петли. 
Иными словами, $\H_c$ — компонента линейной связности нуля в $\H$.

\begin{thm}{Теорема}\label{5.3.A}
Для любого $F \in \F$ 
\[l (F) = \inf_{H\in\H_c} \VERT F - H \VERT.\]
\end{thm}

Доказательство то же, что в \ref{5.1.B}.
В ходе доказательства нужно дополнительно учитывать следующие простые наблюдения:
\begin{itemize}
\item Репараметризация времени, а также процедура регуляризации \ref{5.2} не меняют гомотопический класс пути с фиксированными концами.
Таким образом, минимизировать $\VERT G \VERT$ надо по всем $G \z\in \F$ с $\phi_G = \phi_F$ и таких, что гамильтонов поток $\{g_t\}$ гомотопен $\{f_t\}$ (см. \ref{5.1.C}).
\item Если $g_t = h_t \circ f_t$, где $\{f_t\}$ и $\{g_t\}$ гомотопны с фиксированными концами, то петля $h_t$ стягиваема (сравните с концом доказательства \ref{5.1.B}).
\end{itemize}
Сборка доказательства предоставляется читателю.

\chapter{Лагранжевы пересечения}\label{chap:6}

Теория лагранжевых пересечений изучает одно из самых удивительных явлений симплектической топологии.
В этой главе мы рассмотрим некоторые результаты этой теории, которые в сочетании с идеей линеаризации, изложенной выше, дают довольно мощный инструмент для исследования геометрии группы гамильтоновых диффеоморфизмов.

\section{Точные лагранжевы изотопии}
Пусть $(V^{2n}, \omega)$ — симплектическое многообразие, а $N^n$ —
замкнутое многообразие. 
Рассмотрим лагранжеву изотопию
\[\Phi\: N \times [0;1] \to V,\]
то есть гладкое семейство лагранжевых вложений.
Заметим, что $\Phi^\ast \omega$ должно иметь вид $\alpha_s \wedge \d
s$, где $\{\alpha_s\}$ — семейство 1-форм на $N$ (поскольку
$\Phi^\ast \omega$ обнуляетрся на слоях $N \times \{\point\}$). 
Кроме того, заметим, что $\d \Phi^\ast \omega = \d \alpha_s \wedge \d
s = 0$, откуда следует, что форма $\alpha_s$ замкнута при всех $s$. 

\begin{ex*}{Определение}
Лагранжева изотопия $\Phi$ \rindex{точная изотопия}\emph{точна}, если $\alpha_s$ точна при всех $s$.
\end{ex*}

\begin{ex}{Упражнение}\label{6.1.A}
Покажите, что лагранжева изотопия точна тогда и только тогда, когда
она расширяется до объемлющей гамильтоновой изотопии $V$. 
\emph{Подсказка:} Заметим, что $\alpha_s \z= \d  H_s$ на $N$ и продолжим $H_s
\circ \Phi^{-1}_s$ до нормализованного гамильтониана на $V$. 
\end{ex}

\begin{ex*}{Пример}
Пусть $V$ — поверхность и $N = S^1$.
Лагранжева изотопия $\Phi$ точна тогда и только тогда когда
ориентированная площадь между $\Phi (N \times {0})$ и $\Phi (N \times
{s})$ равна нулю при всех $s$. 
Случай $V = S^2$ показан на рис.~\ref{pic-6}, а  случай $V = \T^\ast S^1 \z= \RR \times S^1$ на рис.~\ref{pic-7}.
Заметим, что в случае цилиндра можно найти симплектическую изотопию, описывающую правую картину.
\end{ex*}

\begin{figure}[ht!]
\begin{minipage}{.48\textwidth}
\centering
\includegraphics{mppics/pic-6}
\end{minipage}\hfill
\begin{minipage}{.48\textwidth}
\centering
\includegraphics{mppics/pic-7}
\end{minipage}

\medskip

\begin{minipage}{.48\textwidth}
\centering
\caption{}\label{pic-6}
\end{minipage}\hfill
\begin{minipage}{.48\textwidth}
\centering
\caption{}\label{pic-7}
\end{minipage}
\vskip-4mm
\end{figure}

Следующий результат играет важную роль в дальнейшем исследовании хоферовской геометрии.
Предположим, что $\{h_t\}$ — петля гамильтоновых диффеоморфизмов, порождённая гамильтонианом $H \in \H$ на $(M, \Omega)$.
Пусть $L \subset M$ — замкнутое лагранжево подмногообразие.
Рассмотрим лагранжеву надстройку (см. п. \ref{3.1.E})
\[L \times S^1 \to (M \times \T^\ast S^1, \Omega + \d r \wedge \d t),\]
\[(x, t) \mapsto (h_t x, -H (h_t x, t), t).\]
Наша цель — исследовать поведение этого лагранжева вложения при однопараметрической деформации.
Пусть $\{h_{t, s}\}$, $s \in [0; 1]$ — гладкое семейство гамильтоновых петель.
Обозначим через $\Phi\: L \z\times S^1 \z\times [0; 1] \z\to M \z\times \T^\ast S^1$ соответствующее семейство лагранжевых надстроек.

\begin{thm}{Теорема}\label{6.1.B}
Лагранжева изотопия $\Phi$ точна.
\end{thm}

Иными словами, гомотопия гамильтоновых петель порождает точную
лагранжевую изотопию лагранжевых надстроек. 
Доказательство теоремы основано на следующем вспомогательном результате.
В его формулировке $H (x, t, s)$ обозначает нормализованный
гамильтониан петли $\{h_{t, s}\}$.
(Тут $t$ и $s$ выполняют разную роль, $t$ — временн\'{а}я переменная, а $s$ — параметр гомотопии.)

\begin{thm}{Предложение}\label{6.1.C}
Равенство 
\[\int_0^1 \frac{\partial H}{\partial s} (h_{t, s} x, t, s)\,\d t = 0\]
выполняется при любых $x \in M$ и $s \in [0; 1]$.
\end{thm}

\parit{Доказательство предложения.}
Доказательство основано на следующей формуле, которая выполняется для произвольного двупараметрического семейства диффеоморфизмов на многообразии.
Проверка формулы предоставляется читателю (см. также \cite{B1}).
Рассмотрим векторные поля $X_{t, s}$ и $Y_{t, s}$ на $M$ определяемые как 
\[\tfrac{\d }{\d t} h_{t, s} x = X_{t, s} (h_{t, s} x)
\quad\text{и}\quad
\tfrac{\d }{\d s} h_{t, s} x = Y_{t, s} (h_{t, s} x).
\]
Тогда 
\[\tfrac{\partial}{\partial s}  X_{t, s}
=
\tfrac{\partial}{\partial t}Y_{t, s} + [X_{t, s}, Y_{t, s}].\]
Обратите внимание, что $X_{t, s}$ и $Y_{t, s}$ являются гамильтоновыми векторными полями при любых $t$ и $s$.
Конечно же, $X_{t, s} = \sgrad H_{t, s}$ для определённого выше $H$.
Напишем $Y_{t, s} = \sgrad F_{t, s}$.
Напомним, что 
\[[\sgrad H, \sgrad F] = -\sgrad  \{H, F\}.\]
Таким образом, получаем, что 
\[\frac{\partial H_{t, s}}{\partial s}
= \frac{\partial F_{t, s}}{\partial t} - \{H_{t, s}, F_{t, s}\}
= \frac{\partial F_{t, s}}{\partial t}+\d F_{t,s}(\sgrad H_{t,s}).
\]
Но последнее выражение, вычисленное в точке $h_{t, s} x$, равно
\[\frac{\d }{\d t} F (h_{t, s} x, t),\]
и мы заключаем, что 
\[\frac{\partial H_{t, s}}{\partial s} (h_{t, s} x)\]
равна полной производной периодической функции.
В частности её интеграл по периоду равен нулю — предложение следует.
\qeds

\parit{Доказательство \ref{6.1.B}.}
Запишем $\Phi^\ast (\Omega + \d r \wedge \d t)$ как $\alpha_s \wedge \d s$.
Требуется проверить, что $\alpha_s$ точна.
Форму $\alpha_s$ можно вычислить явно.

\begin{ex*}[(ср. с п. \ref{3.1.E})]{Упражнение}
Докажите, что равенство
\[\alpha_s (\xi) = \Omega (h_{t, s\ast} \xi, \frac{\partial h_{t, s}}{\partial s}x)\] 
выполняется при всех
$x \in L$, $\xi \in \T_x L$ и 
\[\alpha_s (\tfrac{\partial}{\partial t}) = \frac{\partial H}{\partial s}(h_{t, s}x, t, s).\]
\end{ex*}

Заметим, что первая группа гомологий $H_1 (L \times S^1;\ZZ)$ порождается расщеплёнными циклами вида $C = \beta \times \{0\}$ и $D = \{y\} \times S^1$.
Здесь $\beta$ — цикл на $L$, а $y$ — точка подмногообразия $L$.
Чтобы доказать точность формы $\alpha_s$, достаточно проверить, что её интегралы по всем 1-циклам обнуляются.
Для циклов вида $C$ это следует из того, что $h_{0, s} \equiv \1$ при всех~$s$.
Таким образом, из приведённого выше упражнения следует, что~$\alpha_s$ обращается в нуль на всех векторах, касающихся $L \times \{0\}$.
Далее, 
\[\int_D \alpha_s
= \int_0^1 \frac{\partial H}{\partial s} (h_{t, s} y, t, s)\,\d t. 
\]
По предложению \ref{6.1.C} это выражение равно нулю, что завершает доказательство.
\qeds

\section{Лагранжевы пересечения}

Мы говорим, что лагранжево подмногообразие $N \subset V$ обладает \rindex{свойство лагранжева пересечения}\emph{свойством лагранжева пересечения}, если $N$ пересекает свой образ при любой точной лагранжевой изотопии.
Согласно упражнению~\ref{6.1.A} это можно переформулировать так: $N \cap \phi (N) \ne \emptyset$ при всех $\phi \in \Ham (V, \omega)$, или, иными словами, энергия смещения $N$ бесконечна: $e (N) = + \infty$.

\subsection*{Примеры} 

\begin{ex}[Задача об инфинитезимальном лагранжевом пересечении.]{}\\
Пусть $F$ — автономный гамильтониан на $V$, а $\xi = \sgrad F$ — его гамильтоново
векторное поле. 
Тогда $\xi$ касается $N$ в критических точках $F|_N$ и только в них
(упражнение). 
Поскольку $N$ замкнуто, $F|_N$ должно иметь критические точки, и,
следовательно, нельзя сдвинуть $N$ с себя бесконечно малой
гамильтоновой изотопией. 
\end{ex}

\begin{wrapfigure}{r}{25 mm}
\vskip-0mm
\centering
\includegraphics{mppics/pic-8}
\caption{}\label{pic-8}
\vskip0mm
\end{wrapfigure}

\begin{ex}[Теорема Громова.]{}\label{6.2.B}\rindex{Громов}
В случае если $\pi_2 (V, N) \z= 0$ и $V$ имеет «хорошее» поведение на
бесконечности (скажем, $V$ является произведением замкнутого
многообразия и кокасательного расслоения), Громов \cite{G1} (см. также
статью \rindex{Флоер}Флоера \cite{F}) показал, что $N$ обладает свойством лагранжева
пересечения. 
В частности, это относится к окружности $\{r = 0\}$ в $\T^\ast S^1$
(конечно же, это можно доказать элементарным подсчётом площадей,
см. рис.~7 и обсуждение выше). 
В более общем смысле это справедливо для любой нестягиваемой кривой на
ориентированной поверхности. 
Набросок доказательства теоремы Громова представлен в п. \ref{3.2.G}.
Полное доказательство дано в \cite[Chap. X]{AL}. 

При $\pi_2 (V, N) \ne 0$ свойство лагранжева пересечения может нарушаться.
Возьмём, например, крошечную окружность $N$ на $V \z= S^2$ и сместим её, см. рис.~8.

Тем не менее свойство лагранжева пересечения, очевидно, выполняется для экватора (используйте то, что экватор делит сферу на равновеликие диски).
\end{ex}

\begin{ex}{Определение}\label{6.2.C}
Пусть $L$ — замкнутое лагранжево подмногообразие симплектического многообразия $(M, \Omega)$.
Мы говорим, что $L$ обладает свойством \rindex{устойчивое лагранжево пересечение}\emph{устойчивого лагранжева пересечения}, если $L \times \{r = 0\}$ обладает свойством лагранжева пересечения в $(M \z\times \T^\ast S^1, \Omega + \d r \wedge \d t)$.
\end{ex}

Пара важных примеров приведена в следующих главах.

\begin{ex}[Торы в $\T^\ast \TT^n$.]{}\label{6.2.D}
Рассмотрим лагранжев тор в кокасательном расслоении $\T^\ast \TT^n$ со стандартной симплектической структурой (см. п. \ref{3.1.C}).
Предположим, что он гомологичен нулевому сечению.
Легко проверить, что топологическое предположение теоремы \ref{6.2.B} выполняется.
Следовательно, такие торы обладают свойством устойчивого лагранжева пересечения.
\end{ex}

\begin{ex}[Экватор на $S^2$.]{}\label{6.2.E}
Свойство устойчивого лагранжева пересечения также выполняется для экваторов в $S^2$.
Это следует из сложной теоремы О \cite{O1,O2}, основанной на хитроумной версии гомологий \rindex{Флоер}Флоера.
\end{ex}

Мне не известен пример замкнутого связного лагранжева подмногообразия, обладающего свойством лагранжева пересечения и при этом \emph{не} обладающее его устойчивой версией.

\section{Приложение к гамильтоновым петлям}

Пусть $(M, \Omega)$ — симплектическое многообразие.
Предположим, что $L \subset M$ — замкнутое лагранжево подмногообразие, обладающее свойством устойчивого лагранжева пересечения.  
Пусть $\{g_t\}$ — петля гамильтоновых диффеоморфизмов, порождённая
гамильтонианом $G \in \H$. 
Предположим дополнительно, что
\begin{itemize}
\item $g_t (L) = L$ при всех $t \in S^1$; 
\item $G (x, t) = 0$ при всех $x \in L$, $t \in S^1$.
\end{itemize}
Очевидным примером такой петли является постоянная петля $g_t \z\equiv \1$.
Менее простой пример приведён в \ref{6.3.C}.

\begin{thm}{Теорема}\label{6.3.A}
Пусть $\{h_t\}$ — любая другая петля гамильтоновых диффеоморфизмов, гомотопная описанной выше петле $\{g_t\}$.
Пусть $H \in \H$ — её гамильтониан.
Тогда существуют такие $x \in L$ и $t \in S^1$, что $H (x, t) = 0$.
\end{thm}

Как мы увидим в следующей главе, этот результат приводит к нетривиальным нижним оценкам на хоферовские расстояния.

\parit{Доказательство.} 
Применим дважды конструкцию лагранжевой надстройки к $L$, сначала для $\{g_t\}$, а затем для $\{h_t\}$.
Обозначим через $N_G$ и $N_H$ соответствующие лагранжевы подмногообразия в $M \times \T^\ast S^1$.
По формуле для лагранжевой надстройки мы знаем, что $N_G = L \z\times \{r = 0\}$.
Теорема \ref{6.1.B} говорит, что существует точная лагранжева
изотопия из $N_H$ в $N_G$. 
Таким образом, из свойства устойчивого лагранжева пересечения следует,
что $N_H \cap N_G \ne \emptyset$. 
Пусть $(x, 0, t)$, $x \in L$ — точка пересечения.
Поскольку она лежит в $N_H$, имеем $x = h_t y$ и $0 = -H (h_t y, t, 0)$
при некотором $y \in L$. 
Мы заключаем, что $H (x, t) = 0$.
\qeds

Напомним, что $\H_c$ обозначает пространство всех 1-периодических
гамильтонианов, порождающих стягиваемые петли гамильтоновых
диффеоморфизмов. 
Как непосредственное следствие приведённой выше теоремы мы получаем
следующий результат. 

\begin{thm}{Следствие}\label{6.3.B}
Пусть $L \subset M$ — замкнутое лагранжево подмногообразие,
обладающее свойством устойчивого лагранжева пересечения. 
Тогда для любого $H \in \H_c$ существуют такие $x \in L$ и $t \in S^1$, что $H (x, t) = 0$.
\end{thm}

Например, это верно, когда $M = S^2$, а $L$ — экватор в $S^2$.
Заметим, что утверждение следствия в общем случае становится неверным,
если мы не предполагаем свойства устойчивого лагранжева пересечения.

\begin{ex}{Пример}\label{6.3.C}
Рассмотрим евклидово пространство $\RR^3 (x_1, x_2, x_3)$.
Пусть $M = S^2$ — единичная сфера в пространстве с
индуцированной формой площади. 
Полный оборот вокруг оси $x_3$ — это гамильтонова петля, порождённая
нормализованной функцией Гамильтона $F_1 (x) = 2\pi x_3$. 
(В этом можно убедиться, используя вычисления в \ref{1.4.H}.)
Таким образом, $F_k (x) = 2\pi k x_3$ порождает петлю из $k$ оборотов.
Поскольку гамильтониан $F_k$ тождественно обращается в нуль на экваторе $L = \{x_3
= 0\}$, он должен иметь нуль на \emph{каждой} простой
замкнутой кривой на $S^2$, разделяющей сферу на равновеликие части. 
Обратите внимание, что если $k$ чётно, то $k$ оборотов сферы
$S^2$ представляет собой стягиваемую петлю в $\SO (3)$ и,
следовательно, в $\Ham (S^2)$. 
С другой стороны, $F_k \equiv 2\pi k\epsilon$ на окружности
$C_\epsilon = \{x_3 = \epsilon\}$. 
Отсюда выводим, что явление, описанное в \ref{6.3.B}, очень
жёсткое.%
\footnote{
Это же рассуждение даёт, что при малых  $\epsilon > 0$ подмногообразие $L:=C_{-\epsilon}\cup C_{\epsilon}$ хотя и смещается с себя, тем не менее является стабильно смещяемым.
Недавно Мак и Смит~\cite{MYS} показали, что такое смещение требует некоторого простора.
А именно, для фиксированного достаточно малого $\epsilon>0$ произведение $L$ с экватором сферы достаточно маленькой площади уже несмещаемо.
Это явление было обобщено и использовано для изучения Хоферовской геометрии сферы в~\cite{PS21}. 
  \dpp}
Оно полностью исчезает, если рассматривать окружности, которые делят
сферу на произвольно близкие, но неравные по площади части. 
В самом деле, для любого положительного $\epsilon$ можно выбрать $k$ так,
чтобы $F_k$, сколь угодно велико на $C_\epsilon$!
Решающим, конечно же, является то, что окружность
$C_\epsilon$ не обладает свойством лагранжева пересечения. 
Её можно сместить с себя элементом $\SO (3)$ (ср. с рис.~\ref{pic-8}). 
\end{ex}

\chapter{Диаметр}

В этой главе доказывается, что группа гамильтоновых диффеоморфизмов
замкнутой ориентированной поверхности имеет бесконечный диаметр
относительно хоферовской метрики. 

\section{Начальная оценка}
Пусть $(M, \Omega)$ — симплектическое многообразие и $L \subset M$
— замкнутое лагранжево подмногообразие со свойством  устойчивого
лагранжева пересечения. 
Пусть $F \subset \F$ — такой гамильтониан, что $| F (x, t) | \z\ge C$
при всех $x \in L$ и $t \in S^1$. 
Здесь $C$ — положительная константа.
Следующее предложение даёт оценку снизу на величину $l (F) \z=
\tilde\rho (\tilde\1, \tilde\phi_F)$, введенной в \ref{5.3}. 

\begin{thm}{Предложение}\label{7.1.A}
В данных предположениях $l (F) \ge C$.
\end{thm}

\parit{Доказательство.}
Предложение непосредственно следует из \ref{5.3.A} и \ref{6.3.B}.
Действительно, \ref{6.3.B} утверждает, что каждая функция $H \z\in \H_c$ равна нулю в некоторой точке $(y, \tau)$, где $y \in L$ и $\tau \in S^1$.
Таким образом, $| F (y, \tau) - H (y, \tau) | \ge C$ и, следовательно, $\VERT F - H \VERT \z\ge C$.
Это верно для любого $H \in \H_c$.
Значит \ref{5.3.A} влечёт, что $l (F) \z\ge C$.
\qeds

Мы хотим распространить полученную оценку на $\rho (\1, \phi_F).$
Если группа $\Ham (M, \Omega)$ односвязна относительно
$C^\infty$-топологии (сильной топологии Уитни), то все пути с общими
концами гомотопны и, следовательно, $l (F) = \rho (\1, \phi_F)$. 
Однако если группа $\pi_1 (\Ham (M, \Omega))$
нетривиальна, то нет надежды на обобщение оценки \ref{7.1.A} без
дополнительных предположений. 
Действительно, может случиться, что существует более короткий путь, соединяющий $\1$ с $\phi_F$, который,
конечно, не гомотопен потоку $\{f_t\}$, $t \in [0; 1]$ гамильтониана
$F$. 
Тем не менее, оказывается, что в некоторых интересных случаях эту трудность можно обойти.
Для этого, нужно более внимательно посмотреть на фундаментальную группу $\Ham (M, \Omega)$.

\section{Фундаментальная группа}

О группе $\pi_1 (\Ham (M, \Omega))$ известно немного.
Имеется полная картина для поверхностей (на основе классических
методов) и для некоторых четырёхмерных многообразий \rindex{Громов}
\rindex{Абреу}\cite{G1,A,AM} (на основе теории псевдоголоморфных кривых
Громова). 
В высших размерностях доступны лишь некоторые частичные результаты.
На самом деле я не знаю ни одного симплектического многообразия $M$
размерности $\ge 6$, для которого можно было бы полностью описать группу $\pi_1(\Ham (M, \Omega))$. 
Например, известно, что $\Ham (\RR^{2n})$ односвязно (на самом деле
стягиваемо) при $n = 1, 2$, а при $n = 3$ уже ничего не известно. 
Нам потребуются следующие утверждения о $\pi_1 (\Ham (M, \Omega))$,
где $M$ — замкнутая ориентируемая поверхность. 

\begin{ex}[Сфера (ср. \ref{1.4.H} и \ref{6.3.C}).]{}\label{7.2.A}
Включение $\SO(3) \z\to \Ham(S^2)$ индуцирует изоморфизм фундаментальных групп.
В частности, $\pi_1 (\Ham (S^2))\z=\ZZ_2$.
Нетривиальный элемент получается вращением сферы вокруг вертикальной
оси на один оборот. 
\end{ex}

\begin{ex}[Поверхности рода $\ge 1$.]{}
\label{7.2.B}
В этом случае можно показать, что группа гамильтоновых диффеоморфизмов
односвязна. 
\end{ex}

Эти утверждения хорошо известны специалистам, но, на сколько мне известно,
они не публиковались.
Поэтому я сделаю набросок доказателеьства и добавлю несколько ссылок,
надеясь, что это позволит читателю восстановить полное доказательство.
Мы будем использовать $\Diff_0 (M)$ (соответственно $\Symp_0 (M)$) для
обозначения компоненты единицы группы всех диффеоморфизмов
(соответственно симплектоморфизмов) поверхности $M$. 

\parbf{Набросок доказательства:}

1) \textit{Включение $\pi_1 (\Symp_0 (M)) \to \pi_1 (\Diff_0 (M))$ является
изоморфизмом.}
Чтобы убедиться в этом, рассмотрим пространство $X$ всех форм площади
на $M$ с общей площадью равной $1$. 
Зафиксируем форму площади $\Omega \in X$.
Рассмотрим отображение $\Diff_0 (M) \to X$, переводящее диффеоморфизм
$f$ в $f^\ast \Omega$. 
Можно воспользоваться трюком \rindex{Мозер}Мозера \cite[с. 94--95]{MS}, чтобы
показать, что это отображение является расслоением Серра. 
Заметим, что его слой это не что иное, как $\Symp_0 (M, \Omega)$, а
база $X$ стягиваема. 
После этого желаемое следует из точной гомотопической
последовательности расслоения. 

2) \textit{Топология $\Diff_0 (M)$ известна} (см. \cite{EE}).
В частности, эта группа стягиваема для поверхностей рода $\ge 2$.
Кроме того, если $M \z= S^2$, то она содержит $\SO (3)$ как сильный
деформационный ретракт. 
В случае $M = \TT^2$ она содержит $\TT^2$ как сильный деформационный
ретракт (здесь тор действует на себе сдвигами). 

3) \textit{Включение $j\: \pi_1 (\Ham (M, \Omega)) \to \pi_1 (\Symp_0 (M,
\Omega))$ инъективно} (см. \cite[10.18 iii]{MS}). 
На самом деле это верно для всех замкнутых симплектических многообразий.

4) Собрав все эти утверждения вместе, получаем \ref{7.2.B} для
поверхностей рода $\ge 2$. 
Учитывая, что $\Symp_0 (S^2) = \Ham (S^2)$ (см. \ref{1.4.C}) получаем
\ref{7.2.A}. 
Осталось разобраться со случаем тора $\TT^2$.

5) Выберем точку $y \in \TT^2$ и рассмотрим отображение подстановки $\Diff_0 (\TT^2) \to \TT^2$,
переводящее диффеоморфизм $f$ в точку $f(y)$. 
Оно индуцирует отображение $e_D: \pi_1 (\Diff_0 (\TT^2)) \to \pi_1 (\TT^2)$.
Из шага 2 легко увидеть, что $e_D$ — изоморфизм.
Рассмотрим теперь сужение отображения подстановки на $\Ham (\TT^2)$ и
$\Symp_0 (\TT^2)$ и обозначим через $e_H$ и $e_S$ соответственно
индуцированные ими гомоморфизмы фундаментальных групп. 
Из шага 1 следует, что $e_S$ — изоморфизм.
Используя шаг 3, получаем, что $e_H = e_S \circ j$, где $j$ инъективно.
Из теоремы \rindex{Флоер}Флоера следует, что $e_H$ обращается в ноль (см. \cite{LMP1}).
Таким образом, $\pi_1 (\Ham (\TT^2)) = 0$. 
Доказательство завершено. 
\qeds

\begin{thm}[(\cite{P5})]{Теорема}\label{7.2.C}
Предположим, что 
\begin{itemize}
\item либо $M = S^2$ и $L \subset S^2$ — экватор, 
\item либо $M$ — замкнутая ориентируемая поверхность рода $\ge 1$ и
  $L$ — нестягиваемая простая замкнутая кривая. 
\end{itemize}
Пусть $F \in \F$ — такой гамильтониан, что $|F(x,t)| \ge C$ при всех $x \z\in L$ и $t \in S^1$, где $C$ — положительная постоянная.
Тогда $\rho (\1, \phi_F) \z\ge C$.
\end{thm}

\parit{Доказательство.}
Если $M$ имеет род $\ge 1$, то $l (F) = \rho (\1, \phi_F)$, поскольку группа $\Ham (M, \Omega)$ односвязна (тут мы пользовались утверждением примера \ref{7.2.B}).
Таким образом, результат следует из \ref{7.1.A}.
Если $M = S^2$, то нетривиальный элемент группы $\pi_1 (\Ham (M, \Omega))$ представлен поворотом на один оборот (см. пример \ref{7.2.A}).
Его гамильтониан обращается в нуль на $L$ (см. пример \ref{6.3.C}).
Таким образом, теорема \ref{6.3.A} означает, что каждая функция из $\H$ должна обнуляется в некоторой точке произведения $L \times S^1$.
Необходимое утверждение следует теперь из \ref{5.1.B}.
\qeds

\begin{thm}{Следствие}\label{7.2.D} Группа гамильтоновых диффеоморфизмов замкнутой поверхности имеет бесконечный диаметр относительно хоферовской метрики.
\end{thm}

В самом деле, пусть $M$ и $L$ те же, что в \ref{7.2.C}, и пусть $B \subset M$ — открытый диск, не пересекающийся с $L$.
Возьмём не зависящий от времени гамильтониан $F \in \F$, тождественно равный $C$ вне $B$.
По теореме выше $\rho (\1, \phi_F) \ge C$.
Выбирая $C$ произвольно большим, получаем, что диаметр бесконечен.
Отметим также, что в этом примере носитель функции $\phi_F$ содержится в $B$.
Таким образом, сжимая $B$ и одновременно увеличивая $C$, мы получаем
последовательность гамильтоновых диффеоморфизмов, которая сходится к
тождественному в $C^0$-смысле, но расходится в хоферовской метрике. 

Для замкнутой поверхности $M$ рода $\ge 1$ существует по крайней мере
два других доказательства того, что диаметр $\Ham (M, \Omega)$
бесконечен. 
Одно из них выглядит следующим образом. 

\begin{ex}[(см. \cite{LM2}).]{Упражнение} \label{7.2.E}
Пусть $F$ — нормализованный гамильтониан на $M$.
Предположим, что некоторое регулярное множество уровня $F$ содержит
нестягиваемую замкнутую кривую.  
Рассмотрим поднятие $\tilde f_t$ соответствующего гамильтонова потока
$f_t$ в универсальное накрытие $\tilde M$ поверхности $M$.  
Докажите, что существует $\epsilon> 0$ и семейство дисков $D_t \subset
\tilde M$ площади $\epsilon t$, что $\tilde f_t D_t \cap D_t =
\emptyset$ при всех достаточно больших $t$.  
Выведите из теоремы \ref{3.2.C}, что $\rho (\1, f_t) \to + \infty$ при $t \to + \infty$.  
\end{ex}

Другое доказательство (см. \rindex{Шварц}\cite{Sch3}) основано на
гомологиях Флоера (см. Главу 13 о приложениях гомологий Флоера к
хоферовской геометрии).
Наконец, то, что $\diam \Ham (M, \Omega) = + \infty$, было доказано
для некоторых многомерных многообразий
(см. \rindex{Лалонд}\cite{LM2,Sch3,P5}).
Общий случай остаётся открытым.%
\footnote{
Я. Островер доказал, что универсальное накрытие группы $\Ham(M,\Omega)$ с
поднятием хоферовской метрики имеет бесконечный диаметр для всех
замкнутых симплектических многообразий $(M,\Omega)$, см.~\cite{O03,McD09}.\dpp}

\section{Спектр длин}\label{sec:7.3}

Как было показано, наш метод даёт нижнюю оценку на $\rho (\1, \phi_F)$, при наличии точных данных о фундаментальной группе $\Ham (M, \Omega)$.
Напомним, что в размерностях $\ge 6$ ничего такого нет.
Сейчас мы немного изменим ход рассуждений, расширив при этом класс многообразий, к которым применим метод.
Следующее понятие является одним из главных героев нашего повествования.

\begin{ex}{Определение}\label{7.3.A}
{}\emph{Норму} элемента $\gamma \in \pi_1 (\Ham (M, \Omega))$ определим как \index[symb]{$\nu(\gamma)$}
\[\nu(\gamma)=\inf\length\{h_t\},\]
где нижняя грань берётся по всем гамильтоновым петлям $\{h_t\}$, представляющим $\gamma$.
Множество 
\[\set{\nu (\gamma)}{\gamma \z\in \pi_1 (\Ham (M, \Omega))}\]
называется \rindex{спектр длин}\emph{спектром длин} группы $\Ham (M, \Omega)$.
\end{ex}

\begin{ex*}{Упражнение}

\begin{enumerate}[(i)]
 \item Покажите, что группа $\pi_1 (\Ham (M, \Omega))$ является абелевой. (Используйте то же рассуждение, что и для конечномерных групп Ли.)
Таким образом, мы можем обозначать групповую операцию через $+$, а нейтральный элемент через $0$.
 \item Докажите, что $\nu (\gamma) = \nu (-\gamma)$ и $\nu (\gamma + \gamma') \le \nu (\gamma) + \nu (\gamma')$.
\end{enumerate}

\end{ex*}

На сегодня нет общего утверждения о невырожденности $\nu$ (и
меня не удивит пример, в котором $\gamma \ne 0$, но $\nu (\gamma) = 0$, сравните с
примером \ref{7.3.B}). 
То есть было бы честнее называть $\nu$ псевдонормой, но мы будем
называть её нормой. 
В главе \ref{chap:9} мы опишем метод, который даёт нижние оценки
на $\nu (\gamma)$ и, в частности, позволит вычислить спектр
длин~$S^2$. 

\begin{ex}[Многообразия Лиувилля.]{Пример}\label{7.3.B}\rindex{многообразие Лиувилля}
Мы говорим, что открытое симплектическое многообразие $(M, \Omega)$
обладает \rindex{свойство Лиувиля}\emph{свойством Лиувиля}, если
существует гладкое семейство диффеоморфизмов 
\[D_c: M \to M,\quad c \in (0; + \infty),\]
что $D_1 = \1$ и $D_c^\ast \Omega = c\Omega$ при всех $c$.
Такие диффеоморфизмы $D_c$ конечно, не могут иметь компактных
носителей при $c \ne 1$. 
Важный класс примеров это кокасательные расслоения со стандартной
симплектической структурой (см. пример \ref{3.1.C}).
Диффеоморфизм $D_c$ в этом случае задаётся послойной гомотетией $(p,
q) \z\mapsto (cp, q)$. 
Мы утверждаем что спектр длин $\Ham (M, \Omega)$ равен $\{0\}$ при
условии, что $(M, \Omega)$ 
обладает свойством Лиувиля.
Доказательство этого утверждения основано на следующем простом наблюдении.
\end{ex}

\begin{ex*}{Упражнение}
Пусть $\{f_t\}$ — гамильтонов поток, порождённый нормализованным гамильтонианом $F (x, t)$.
Тогда при любом $c> 0$ поток $\{D_c f_t D_c^{-1}\}$ снова является
гамильтоновым и его нормализованный гамильтониан равен $cF (D_c^{-1} x, t)$.
\end{ex*}

Пусть $\{h_t\}$ — петля гамильтоновых диффеоморфизмов.
Рассмотрим семейство гомотопных петель $\{D_c h_t D_c ^{-1}\}$.
Из упражнения выше следует, что длины петель стремятся к нулю при $c \to 0$, поэтому каждую петлю можно продеформировать в петлю произвольной малой длины.
Мы заключаем, что спектр длин равен $\{0\}$ (без каких-либо сведений о фундаментальной группе).
В следующей главе мы разберём этот пример в контексте классической механики.

\section{Уточнение оценки}

\begin{thm}{Теорема}\label{7.4.A}
Пусть $(M, \Omega)$ — симплектическое многообразие и $L \subset M$
— замкнутое лагранжево подмногообразие со свойством устойчивого
лагранжева пересечения. 
Предположим, что спектр длин группы $\Ham (M, \Omega)$ ограничен
сверху некоторым $K \ge 0$. 
Далее предположим, что гамильтониан $F \in \F$ такой, что $| F (x,
t) | \z\ge C$ при всех $x \in L$ и $t \in S^1$. 
Тогда 
\[\rho (\1, \phi_F) \ge C - K.\]
\end{thm}

\parit{Доказательство.}
Выберем произвольное $\epsilon> 0$. 
Пусть $\{g_t\}$ — путь гамильтоновых диффеоморфизмов, соединяющий $\1$ с $\phi_F$.
Рассмотрим петлю $\{f_t \circ g_t ^{-1}\}$.
По нашему предположению эта петля гомотопна петле $\{h_t\}$, длина которой не превосходит $K + \epsilon$.
Путь $\{f_t\}$ гомотопен с фиксированными концами с композицией $\{h_t \circ g_t\}$.
Следовательно
\[l (F) \le \length \{h_t\} + \length \{g_t\}.\]
Поскольку $l (F) \ge C$, то в силу \ref{7.1.A} получаем, что $\length \{g_t\} \ge C \z- K - \epsilon$.
Таким образом, $\rho (\1, \phi_F) \ge C - K$.
\qeds

\chapter{Рост и динамика}

В этой главе обсуждается асимптотическое геометрическое поведение
однопараметрических подгрупп в группе гамильтоновых диффеоморфизмов и
описывается связь между геометрией и инвариантными торами в
классической механике. 

\section[Инвариантные торы]{Инвариантные торы в классической механике}

Инвариантные лагранжевы торы гамильтоновых динамических систем играют
важную роль в классической механике. 
Начнём с препятствия к существованию инвариантных торов, которое
происходит из геометрии группы гамильтоновых диффеоморфизмов
(см. \ref{8.1.C}). 

Рассмотрим $n$-мерный тор $\TT^n$ с евклидовой метрикой $ds^2 \z=
\sum^n_{j = 1} dq_j^2$. 
Евклидов геодезический поток описывается гамильтоновой системой на
кокасательном расслоении $\T^\ast \TT^n$ со стандартной симплектической
формой $\Omega = dp \wedge dq$. 
Функция Гамильтона задается формулой $F (p, q) = \tfrac12 | p |^2$.
Решая гамильтонову систему 
\[
\begin{cases}
\quad\dot p &= 0,\\
\quad\dot q &= p,
\end{cases}
\]
находим гамильтонов поток $f_t (p, q) = (p, q + pt)$.
Опишем динамику этого потока.
Каждый тор $\{p = a\}$ инвариантен относительно $\{f_t\}$. 
Более того, ограничение потока на каждый такой тор представляет собой
обычное (квази)-периодическое движение $q \mapsto q + at$. 
Все эти торы гомологичны нулевому сечению кокасательного расслоения.
Геометрически они соответствуют семействам параллельных евклидовых
прямых на $\TT^n$. 

Приведенный выше поток принадлежит важному классу независящих от
времени гамильтоновых систем, называемых \rindex{интегрируемая система}\emph{интегрируемыми} системами
(см. \cite{Ar}). 
Интегрируемые системы можно охарактеризовать тем, что их
энергетические уровни расслоены (в дополнении множеств меры ноль)
инвариантными торами средней размерности с (квази-)периодической
динамикой. 
Традиционно они считаются простейшей динамикой в классической механике.
Возникает естественный вопрос: что происходит с инвариантными торами
при возмущении системы. 

Теория \rindex{Колмогоров}\rindex{Арнольд}\rindex{Мозер}Колмогорова — Арнольда — Мозера (или \rindex{КАМ-теория}КАМ-теория, см. \cite{Ar})
говорит нам, что при малых возмущениях выживает б\'{о}льшая часть
инвариантных торов. 
Однако необходимо потребовать, чтобы вектор вращения был «достаточно
иррациональным». 
В частности, это означает, что если в некоторых угловых координатах
$\theta$ на инвариантном торе динамика задается формулой $\dot\theta =
a$, где $a = (a_1 ,\dots, a_n) \in \RR^n$, то $a_1 ,\dots, a_n$ линейно независимы над $\QQ$. 
При б\'{о}льших возмущениях «топологически существенные» торы могут исчезнуть.
Можно, например, продеформировать евклидову метрику на $\TT^n$ до римановой метрики, геодезический поток которой не имеет инвариантных торов, несущих квазипериодическое движение и гомологичных нулевому сечению (см. \cite{AL}). 

Заметим, что в нашем начальном примере инвариантные торы $\{p = a\}$
лагранжевы. 
Это общее явление, которое сейчас прояснится.

\begin{ex}{Упражнение}\label{8.1.A}
Пусть $F\: M \to \RR$ — независящий от времени гамильтониан на симплектическом многообразии~$M$.
Рассмотрим замкнутое подмногообразие $L \subset \{H = c\}$.
Покажите, что если $L$ лагранжево, то $L$ инвариантно относительно гамильтонова потока~$F$.
\emph{Подсказка:} используйте линейную алгебру, чтобы доказать, что $\sgrad F$ касается~$L$.
\end{ex}

В некоторых интересных случаях приведенное выше утверждение можно обратить.

\begin{thm}[(\cite{He})]{Предложение}\label{8.1.B}\rindex{Эрман}
Пусть $F\: \T^\ast \TT^n \to \RR$ — гамильтониан, независящий от времени.
Рассмотрим инвариантный тор $L \subset \{H = c\}$, несущий квазипериодическое движение $\dot\theta = a$, где координаты вектора $a$, как и раньше, линейно независимы над $\QQ$.
Тогда $L$ лагранжев.
\end{thm}

\parit{Доказательство.}
Возьмём точку $x \in L$.
Предположим, что $\Omega|_{\T_x L}$ имеет вид $\sum b_{ij} d\theta_i
\wedge d\theta_j$. 
Поскольку динамика — это просто сдвиг, имеем $\theta (t) = \theta
(0) + at$, и, значит, $\Omega|_{\T_y L} = \sum b_{ij} d\theta_i \wedge d\theta_j$ для каждой точки $y$, лежащей на траектории $x$. 
Заметим, что каждая траектория плотна на торе.
Поэтому $\Omega = \sum b_{ij} d\theta_i \wedge d\theta_j$ всюду на~$L$.
Но $\Omega|_{\T L}$ — точная 2-форма.
Следовательно, $b_{ij} = 0$ при всех $i$, $j$, и, значит, $\Omega|_{\T L} = 0$. 
Мы показали, что $L$ лагранжево.
\qeds

Пусть $F$ — независящий от времени гамильтониан на $(\T^\ast \TT^n, \Omega)$ с компактным
носителем. 
Определим число \index[symb]{$E(F)$}$E(F)=\sup|E|$, где точная верхняя грань берётся
по всем таким $E$, что уровень энергии $\{F \z= E\}$ \textit{содержит
лагранжев тор, гомологичный нулевому сечению}. 
Как найти нетривиальную оценку сверху на $E (F)$?
Такого сорта вопросы изучаются в рамках обратной КАМ-теории.%
\footnote{Цель КАМ-теории — доказать существование инвариантных
  лагранжевых торов, в то время как обратная КАМ-теория изучает
  препятствия к их существованию, см., например \cite{Mac} и
  приведённые там ссылки.} 
Обозначим через $\{f_t\}$ гамильтонов поток, порожденный $F$.

\begin{thm}[(ср. \cite{BP2,P8})
]{Теорема}\label{8.1.C}
  $\rho (\1, f_t) \ge E (F) t$ при всех $t \z\in \RR$.
\end{thm}

\parit{Доказательство.}
Отметим, что достаточно доказать неравенство при $t = 1$ (время можно
репараметризовать).  
Пусть $L \z\subset \{H = E\}$ — лагранжев тор, гомологичный нулевому сечению.
Тогда $L$ обладает свойством устойчивого лагранжева пересечения
(см. пример \ref{6.2.D}). 
Поскольку $(\T^\ast \TT^n, \Omega)$ — многообразие Лиувилля, спектр
длин $\Ham (M, \Omega)$ равен $\{0\}$ (см. \ref{7.3.B}). 
Следовательно, все условия теоремы \ref{7.4.A} выпонены. 
Из этой теоремы следует, что $\rho (\1, f_1) \ge E$. 
Взяв верхнюю грань по всем таким $E$, получаем нужную оценку.
\qeds

Геометрическое содержание этой оценки следует разуметь в более общем контексте роста однопараметрических подгрупп гамильтоновых диффеоморфизмов.

\section{Рост однопараметрических подгрупп}\label{sec:8.2}

Пусть $(M, \Omega)$ — симплектическое многообразие и $\{f_t\}$ —
однопараметрическая подгруппа в $\Ham (M, \Omega)$, порожденная
нормализованным гамильтонианом $F \in \A$. 
Одна из центральных задач хоферовской геометрии — изучить
взаимосвязь между функцией $\rho (\1, f_t)$ и динамикой потока
$\{f_t\}$.
Например, теорема \ref{8.1.C} утверждает, что инвариантные торы
автономного гамильтонова потока на $\T^\ast \TT^n$,  гомологичные
нулевому сечению и с квазипериодической динамикой, вносят вклад в
линейный рост функции $\rho (\1, f_t)$. 
Существует еще одна, чисто геометрическая причина интереса к этой функции.
Она исходит из теории геодезических хоферовской метрики (см. главу~\ref{chap:12}). 
Гамильтонов путь $\{f_t\}$ называется
\rindex{кратчайшая}\emph{кратчайшей}, если каждый из его отрезков
минимизирует длину между своими концами.
Предположительно (см. \ref{12.6.A}) все {}\emph{достаточно короткие}
отрезки любой однопараметрической подгруппы являются кратчайшими
(иными словами, каждая однопараметрическая подгруппа локально
минимизирует длину). 
Однако, как мы увидим в \ref{8.2.H}, {}\emph{длинные отрезки} могут перестать
быть кратчайшими. 
Нарушение минимальности — любопытное явление, которое всё ещё не изучено.
В настоящее время к нему известно несколько подходов.
Один из них основан на теории сопряженных точек в хоферовской
геометрии, он обсуждается в главе \ref{chap:12}. 
Здесь мы обсудим подход, основанный на понятии \rindex{асимптотический рост}\emph{асимптотического роста}
однопараметрической подгруппы (см. \rindex{Бялый}\cite{BP2}). 
Асимптотический рост определяется следующим образом: \index[symb]{$\mu(F)$}
\[\mu(F)=\lim_{t \to + \infty}\frac{\rho (\1, f_t)}{ t \| F \|}.\]

\begin{ex*}{Упражнение}
Покажите, что указанный выше предел существует.
\emph{Подсказка:} используйте субаддитивность функции $\rho (\1, f_t)$, то есть $\rho (\1, f_{t + s}) \le \rho (\1, f_t) + \rho (\1, f_s)$. 
\end{ex*}

Ясно, что $\mu (F)$ лежит в интервале $[0; 1]$.
Если $\mu (F) <1$, то путь $\{f_t\}$ не является кратчайшей.

Рассмотрим несколько примеров поведения функции $\rho (\1, f_t)$.
\rindex{Хофер}Х. Хофер \cite{H2} доказал, что каждая однопараметрическая подгруппа в
$\Ham (\RR^{2n})$ локально кратчайшая. 
С другой стороны, \rindex{Сикорав}Ж.-К. Сикорав \cite{S2} обнаружил поразительное явление --- каждая такая однопараметрическая подгруппа лежит на конечном расстоянии от тождественного отображения (в частности $\mu=0$).
Таким образом, вся подгруппа не может быть кратчайшей. 

\begin{thm}{Теорема}\label{8.2.A}
Пусть $\{f_t\}$ — однопараметрическая подгруппа в $\Ham (\RR^{2n})$,
порожденная гамильтоновой функцией $F$ с компактным носителем. 
Предположим, что носитель $F$ содержится в евклидовом шаре радиуса
$r$. 
Тогда функция $\rho (\1, f_t)$ ограничена: $\rho (\1, f_t) \z\le 16\pi
r^2$. 
\end{thm}

Мы отсылаем читателя к \cite[с. 177]{HZ} за полным доказательством (см. также обсуждение в предложении \ref{12.6.E}). 

Вернёмся теперь к однопараметрическим подгруппам в $\Ham (\T^\ast \TT^n)$.
Прежде всего доказано, что все они локально кратчайшие \cite{LM2}.
Иными словами, $\rho (\1, f_t) = t \| F \|$ при условии, что $t$
достаточно мало. 
Конечно же, это означает, что в общем случае оценка в \ref{8.1.C} не является
точной при малых $t$. 
Действительно, в общем случае $E (F)$ строго меньше, чем $\| F \|$.
Тем не менее в случае $n = 1$ и $F \ge 0$ оценка \ref{8.1.C}
асимптотически точна! 
Заметим, что, поскольку каждая замкнутая кривая на цилиндре $\T^\ast
\TT^1$ является лагранжевой, величина $E (F)$ равна верхней грани тех
вещественных чисел $E$, для которых уровень $\{F = E\}$ содержит
нестягиваемую вложенную окружность. 

\begin{thm}[(\cite{PS})]{Теорема}\label{8.2.B}\rindex{Зибург}
  Пусть F — неотрицательный гамильтониан с компактным носителем на
  цилиндре $\T^\ast \TT^1$ с $\| F \| = 1$. 
  Тогда обратный КАМ-параметр $E (F)$ совпадает с асимптотическим
  ростом $\mu (F)$: $E (F) = \mu (F)$. 
\end{thm}

\parit{Доказательство.}
Нам достаточно доказать, что $\mu (F) \le E (F)$.
Тогда, примениив \ref{8.1.C}, мы получим желаемый результат.

Если $E (F) \z= \max F = 1$, то \ref{8.1.C} влечёт $\mu (F) = E (F)$.
Предположим теперь, что $E (F) <1$.
Идея состоит в разложении потока $\{f_t\}$ на два
коммутирующих потока с простой асимптотикой. 
Выберем $\epsilon> 0$ достаточно малым и рассмотрим гладкую
неубывающую функцию $u: [0; + \infty) \to [0; + \infty)$ со следующими
    свойствами:  
\begin{itemize}
\item $u (s) = s$ при $s \le E (F) + \epsilon$;
\item $u (s) = E (F) + 2\epsilon$ при $s \ge E (F) + 3\epsilon$;
\item $u (s) \le s$ при всех $s$.
\end{itemize}
Рассмотрим новые гамильтонианы $G = u \circ F$ и $H = F - G$ и
обозначим через $\{g_t\}$ и $\{h_t\}$ соответствующие гамильтоновы
потоки. 
Эти потоки коммутируют и $f_t = g_t h_t$. 
Таким образом, 
\begin{equation}\rho (\1, f_t) \le \rho (\1, g_t) + \rho (\1, h_t).
\label{eq:8.2.C}
\end{equation}
Обратите внимание, что $\| G \| \le E (F) + 2\epsilon$.
Следовательно, 
\begin{equation}
 \rho (\1, g_t) \le t (E (F) + 2\epsilon).
\label{eq:8.2.D}
\end{equation}
Далее, носитель $H$ содержится в подмножестве $D_\epsilon = \{F \ge E
(F) + \epsilon\}$. 
Для достаточно малого $\epsilon$ общего положения множество
$D_\epsilon$ является областью, граница которой состоит из {}\emph{стягиваемых}
замкнутых кривых (тут используется определение величины $E(F)$). 
Предположим теперь, что носитель гамильтониана $F$ содержится в кольце 
\[A = \set{(p, q) \in \T^\ast\TT^1}{| q | \le a / 2}\]
для некоторого $a> 0$. 
Отметим, что $\partial D_\epsilon \subset A$.
Следовательно, множество $D_\epsilon$ содержится в некотором множестве
$D' \subset A$, которое является конечным объединением попарно
непересекающихся замкнутых дисков с общей площадью не более $a$. 
Поскольку цилиндр имеет бесконечную площадь, из теоремы
Дакорогны — \rindex{Мозер}Мозера \cite[1.6]{HZ} легко вытекает существование симплектического вложения $i\: \RR^2 \to\T^\ast \TT^1$ и конечного объединения $D'' \subset \RR^2$ евклидовых дисков таких, что $i$
отображает $D''$ диффеоморфно на $D'$. 
Ясно, что $i$ индуцирует естественный гомоморфизм 
\[i_\ast: \Ham (\RR^2) \to \Ham (\T^\ast \TT^1).\]
Важно отметить, что $i_\ast$ не увеличивает хоферовские расстояния.
Наш поток $h_t$ лежит в образе $i_\ast$, то есть $h_t = i_\ast (e_{t})$,
где $e_{t}$ — однопараметрическая подгруппа в $\Ham (\RR^2)$,
гамильтониан которой имеет носитель в $D''$. 
Таким образом, теорема \ref{8.2.A} влечёт, что 
\[
\rho (\1, h_t) \le 16a.
\]
Комбинируя это неравенство с (\ref{eq:8.2.D}) и (\ref{eq:8.2.C}) получаем, что
\begin{equation}
  \rho (\1, f_t) \le t (E (F) + 2\epsilon) + 16a 
  \label{eq:8.2.E} 
\end{equation}
при всех $t> 0$.
Поделив на $t$ и перейдя к пределу при $t \to +\infty$, получаем, что
$\mu (F) \le E (F) + 2\epsilon$. 
Теорема следует поскольку $\epsilon$ произвольно мало.
\qeds

\begin{ex}{Замечание}\label{8.2.F}
То же доказательство показывает, что если $E (F) \z= 0$, то функция
$\rho (\1, f_t)$ ограничена. 
Действительно, поскольку (\ref{eq:8.2.E}) выполняется при всех
$\epsilon> 0$, получаем, что $\rho (\1, f_t) \z\le 16a$. 
Принимая во внимание $E (F) \le \mu (F)$, получаем следующее
утверждение типа «жёсткости»: \textit{если $\mu (F) = 0$, то функция $\rho
(\1, f_t)$ ограничена} (дальнейшее обсуждение в \ref{sec:8.4}). 
\end{ex}

Теорема \ref{8.2.B} и замечание \ref{8.2.F} верны для всех открытых
поверхностей бесконечной площади (см. \cite{PS}). 
Более того, можно изменить определение величины $E (F)$ и распространить эти
утверждения на произвольные (не обязательно неотрицательные)
гамильтонианы $F$. 
Пока что нет обобщений \ref{8.2.B} на высшие размерности.
Однако следующее рассуждение показывает, что оценка \ref{8.1.C} бывает точной,
по крайней мере в следующем очень частном случае. 

Пусть $F\: \T^\ast \TT^n \to \RR$ — гамильтониан с компактным носителем, независящий от времени, который удовлетворяет следующим условиям: 
\begin{enumerate}[(i)]
\item $F \ge 0 $
\item $\max F = 1$ 
\item множество максимума $\Sigma =\{F = 1\}$ является гладким сечением кокасательного расслоения. 
\end{enumerate}

Оказывается, геометрия соответствующего потока $\{f_t\}$ сильно
зависит от того, является ли $\Sigma$ лагранжевым или нет! 

Предположим, что $\Sigma$ — лагранжево подмногообразие.
Тогда по определению $E (F) = 1$.
Следовательно, теорема \ref{8.1.C} влечёт, что $\{f_t\}$ —
кратчайшая и, в частности, $E (F) = \mu (F)$. 

Предположим теперь, что $\Sigma$ нелагранжево, то есть $\Omega$ не
обращается в нуль хотя бы на одном касательном пространстве к
$\Sigma$. 
Тогда неизвестно, равны ли $\mu (F)$ и $E (F)$ между собой.
Однако мы утверждаем, что оценка из \ref{8.1.C} по крайней мере
нетривиальна, а именно  
\begin{equation}
E (F) \le \mu (F) <1.\label{8.2.G}
\end{equation}

Чтобы объяснить это неравенство, нам понадобится следующий довольно общий результат.
В некоторых интересных случаях он позволяет доказать, что $\mu(F)$ строго меньше $1$. 

\begin{thm}{Теорема}\label{8.2.H}
Пусть $F$ — независящий от времени нормализованный гамильтониан на
симплектическом многообразии $ (M, \Omega)$. 
Пусть $\Sigma_+$  и $\Sigma_-$ 
  множества минимума и максимума $F$ соответственно.
Предположим, что существует такой диффеоморфизм $\phi \in \Ham (M, \Omega)$, что
либо $\phi (\Sigma_+) \cap \Sigma_+ = \emptyset$, либо $\phi
(\Sigma_-) \cap \Sigma_- = \emptyset$. 
Тогда $\mu (F) <1$ и, в частности, гамильтонов поток
гамильтониана $F$ не образует
кратчайшую.  
\end{thm}

Доказательство теоремы приведено в разделе \ref{sec:8.3}.

Приведём набросок доказательства \ref{8.2.G}.
Поскольку $\Sigma$ нелагранжево, можно показать, что существует такой гамильтонов диффеоморфизм $\phi$, что
$\phi(\Sigma) \cap \Sigma = \emptyset$. 
Таким образом, комбинируя теоремы \ref{8.1.C} и \ref{8.2.H} получаем,
что $E (F) \le \mu (F) <1$. 
Доказательство существования гамильтонова диффеоморфизма $\phi$,
смещающего нелагранжево подмногообразие $\Sigma$, сложное. 
Оно основано на $h$-принципе \rindex{Громов}Громова для отношений в частных производных.
На самом деле можно доказать большее, а именно что энергия смещения
$\Sigma$ обращается в нуль. 
Мы отсылаем читателя к \rindex{Лауденбах}\rindex{Сикорав}\cite{P2,LS} за доказательствами этого
результата и его обобщений. 

Отметим также, что \rindex{Зибург}К. Ф. Зибург \cite{Si2} обобщил неравенство \ref{8.1.C}
на {}\emph{неавтономные} гамильтоновы потоки. 
Обратите внимание, что в случае, зависящем от времени, уже
нетривиально определить обратный КАМ-параметр из-за отсутствия
сохранения энергии. 
Определение Зибурга основано на идее минимального действия,
разработанной \rindex{Мазер}Дж. Мазером. 

\section[Вырямление кривых]{Вырямление кривых}\label{sec:8.3}

Пусть $(M, \Omega)$ — симплектическое многообразие.
Для любой функции $F \in \A$, $F \not\equiv 0$, положим 
\[\delta(F)=\inf_\phi \frac{\|F+ F \circ \phi\|}{2\|F\|},\]
где точная нижняя грань берётся по всем $\phi \in \Ham (M, \Omega)$.

\begin{thm}[(\cite{BP2})]{Теорема}\label{8.3.A}
\[\mu (F) \le \delta (F).\]
\end{thm}

Теорема \ref{8.2.H} является непосредственным следствием теоремы выше. 
В самом деле, если $F$ удовлетворяет предположениям \ref{8.2.H}, то
$\delta (F) <1$. 

\parit{Доказательство \ref{8.3.A}.}
Не умаляя общности, можно считать, что $\|F \| = 1$.
Выберем $\phi \in \Ham (M, \Omega)$ и $T> 0$.
Запишем 
\[f_{2T}= (f_T \circ \phi \circ f_T \circ \phi^{-1}) \circ (\phi \circ
f_T^{-1} \circ \phi^{-1} \circ f_T) = A \circ B.\] 
Заметим, что $B$ — коммутатор, поэтому $\rho (\1, B) \le 2\rho (\1, \phi)$.
Диффеоморфизм $A$ порождается путем $g_t = f_t \phi f_t \phi ^{-1}$
при $t \in [0;T]$, и его гамильтониан равен  
\[G (x, t) = F (x) + F (\phi^{-1} f_t^{-1} x).\]
Таким образом, имеем 
\begin{align*}
\|G_t \| &= \|F + F \circ \phi ^{-1} \circ f_t ^{-1} \| =
\\
&=\|F \circ f_t + F \circ \phi ^{-1} \| =
\\
&=\|F + F \circ \phi ^{-1} \|,\end{align*}
поскольку $F \circ f_t = F$ (по закону сохранения энергии).
Итак, 
\[\rho (\1, f_{2T}) \le T \| F + F \circ \phi ^{-1} \| + 2\rho (\1, \phi)\]
и, следовательно, 
\[\frac{\rho (\1, f_{2T})}{2T}
\le\frac12 \|F + F \circ \phi ^{-1} \| + \frac{\rho (\1, \phi)}{T}\]
при всех $\phi \in \Ham (M , \Omega)$.
Переходя к пределу при $T \to \infty$, получаем $\mu (F) \le \delta (F)$.
\qeds

Отсылаем читателя к \cite{LM2} и \cite{P9}, где описаны различные
процедуры выпрямления кривых в хоферовской геометрии. 
Мы вернемся к этому вопросу в главе \ref{chap:11}.
Величина $\delta (F)$ допускает следующее естественное обобщение.
Положим 
\[
\delta_N(F)
=
\inf_{\phi_1 ,\dots, \phi_{N-1}}
\frac1N
\frac{\|\sum_{j=0}^{N-1} F \circ \phi_j\|}{\|F\|}
\] 
где $\phi_0 = \1$, а нижняя грань берётся по всем последовательностям $\{\phi_j\}$, $j = 1,\dots, N-1 $ гамильтоновых диффеоморфизмов. 
На этом языке $\delta = \delta_2$. 
Можно показать \cite{P9}, что на замкнутом симплектическом многообразии $\delta_N (F) \to 0$ при $N \to + \infty$ для любой функции $F \in \A$.
Заметим, что при $N>2$ неравенство $\mu(F)\le\delta_N(F)$, напоминающее теорему~\ref{8.3.A},
становится вообще говоря неверным.
Мне неизвестна скорость убывания последовательности $\delta_N (F)$. 

\section{А что если рост нулевой?}\label{sec:8.4}

Как мы видели в \ref{8.2.F}, на цилиндре 
(и, в более общем смысле, на открытых поверхностях бесконечной площади)
каждая однопараметрическая подгруппа с нулевым асимптотическим ростом является ограниченной. 
А что происходит на других симплектических многообразиях?

\begin{ex}{Задача}\label{8.4.A}
Существует ли симплектическое многообразие $(M, \Omega)$ и
однопараметрическая подгруппа $\{f_t\}$ в $\Ham (M, \Omega)$ такие,
что функция $\rho (\1, f_t)$ имеет промежуточный рост
(например, растет как $\sqrt{t}$)? 
\end{ex}

Эта задача открыта даже для такого простого симплектического
многообразия, как двумерный тор $\TT^2$. 
В оправдание невежества, заметим, что, как показывает следующий
результат, гамильтониан $F \z\in \A (\TT^2)$ с $\mu (F) = 0$, не может находится в общем положении.%
\footnote{
Для любого замкнутого симплектического многообразия и
$C^{\infty}$-общей функции $F$ на нём выполняется $\mu(F)>0$,
см.~\cite[Section 6.3.1]{PR14}.\dpp}

\begin{thm}{Теорема}\label{8.4.B}
Если 0 — регулярное значение $F \in \A (\TT^2)$, то $\mu (F)> 0$.
\end{thm}

\parit{Доказательство.}
Множество $D = \{F = 0\}$ состоит из конечного числа попарно
непересекающихся вложенных окружностей. 
Значит, существует нестягиваемая простая замкнутая кривая $L \z\subset
\TT^2$ такая, что $L\cap D = \emptyset$. 
Таким образом, $|F (x)| \z> C$ при фиксированном $C> 0$ и всех $x \in L$.
Поскольку $\pi_1 (\Ham (\TT^2)) = 0$ (см. пример \ref{7.2.B}), из
\ref{7.4.A} следует, что $\rho (\1, f_t) \ge Ct$ при всех $t$ и, в
частности, $\mu (F)> 0$. 
\qeds

\chapter{Спектр длин}\label{chap:9}

В этой главе мы опишем метод вычисления спектра длин в хоферовской геометрии.
Наш подход основан на теории симплектических расслоений над двумерной сферой.

\section{Положительная и отрицательная части хоферовской нормы}

Пусть $(M, \Omega)$ — замкнутое симплектическое многообразие.
Для $\gamma \in \pi_1 (\Ham(M, \Omega))$ положим \index[symb]{$\nu_\pm$}
\begin{align*}
\nu_+ (\gamma) &= \inf_F \int_0^1 \max_x F (x, t)\,dt = \inf_F \max_{x,t}F (x, t),
\\ 
\nu_- (\gamma) &= \inf_F \int_0^1 \big(-\min_x F (x, t)\big)\,dt = \inf_F \min_{x,t} \big(-F (x, t)\big),
\end{align*}
где точная нижняя грань берётся по всем нормализованным периодическим
гамильтонианам $F \in \H$, порождающим петлю в классе~$\gamma$. 
Равенства в определениях доказываются точно так же, как и лемма \ref{5.1.C}.

\begin{ex*}{Упражнение}
Докажите, что $\nu_+ (\gamma) \z= \nu_- (-\gamma)$ и $\nu(\gamma) \z\ge
\nu_- (\gamma) \z+ \nu_+ (\gamma)$ 
(сравните с открытой задачей в \ref{sec:2.4}).
\end{ex*}

Вычисление $\nu_+ (\gamma)$ оказывается нетривиальным даже в следующем
простейшем случае. 
Рассмотрим симплектическое многообразие $(S^2, \Omega)$ нормализованное так, что $\int_{S^2} \Omega = 1$.
Пусть $\gamma$ — класс одного оборота $\{f_t\}$, и пусть $F \in \H$
— гамильтониан, порождающий $\{f_t\}$. 
Как видно из \ref{6.3.C}, $\max F \z= - \min F = \frac12$.
(Множитель $2\pi$ из \ref{6.3.C} исчез ввиду приведенной выше
нормировки формы $\Omega$.) 
Таким образом, $\nu_+ (\gamma) \le \frac12$.
Но на самом деле выполняется равенство!

\begin{thm}{Теорема}\label{9.1.A}
$\nu_+ (\gamma) = \frac12$.
\end{thm}

Гамильтонова петля $\{f_t \}$, представляющая класс $\gamma \ne 0$,
называется \rindex{замкнутая кратчайшая}\emph{замкнутой кратчайшей}, если $\length\{f_t \} =
\nu(\gamma)$.
Обратите внимание, что замкнутая кратчайшая не является кратчайшей.

\begin{thm}[(\cite{LM2})]{Следствие} $\nu(\gamma) = 1$ и $\{f_t\}$ —
  замкнутая кратчайшая.
\end{thm}

\parit{Доказательство.}
Прежде всего, $\length\{f_t\} =\frac12 - (-\frac12 ) = 1$ и,
следовательно, $\nu(\gamma) \le 1$. 
Поскольку $2\gamma = 0$
\[\nu_- (\gamma) = \nu_+ (-\gamma) = \nu_+ (\gamma) =\tfrac12.\]
Кроме того, $\nu(\gamma) \ge \nu_- (\gamma) + \nu_+ (\gamma) = 1$.
Отсюда $\nu(\gamma) = 1$ и $\{f_t\}$ — замкнутая кратчайшая.
\qeds

Теорема \ref{9.1.A} доказана в разделе~\ref{sec:9.4}.
Обобщение на $\CP^n$ при $n \ge 2$ дано в \cite{P3}.

\section[\texorpdfstring{Симплектические расслоения над
    $S^2$}{Симплектические расслоения над S²}]{Симплектические
  расслоения над $\bm{S^2}$} 
\label{sec:9.2}

Пусть $(M, \Omega)$ — замкнутое симплектическое многообразие.
Далее будем считать, что $H^1 (M;\RR) = 0$.%
\footnote{От этого предположения легко избавиться, см. теорию \emph{гамильтоновых} симплектических расслоений в \cite{MS} и \cite{P4}.} 
Как следствие, $\Ham(M, \Omega)$ совпадает со связной компонентой $\1$
в $\Symp(M, \Omega)$. 
Пусть $p \: P \to S^2$ — гладкое расслоение со слоем $M$ со
следующей послойной симплектической структурой.
Для каждого $x \in S^2$ на $p^{-1} (x)$ задана такая симплектическая форма
$\Omega_x$, что $\Omega_x$ гладко зависит от $x$ и $(p^{-1} (x),
\Omega_x)$  симплектоморфно $(M, \Omega)$.
Кроме того, мы всегда выбираем ориентацию на $S^2$ как часть данных
(поэтому $P$ также ориентировано).
Назовём $p\: P\to S^2$ \rindex{симплектическое расслоение}\emph{симплектическим расслоением} (подробнее
см. \cite{MS}). 

Каждая петля $\{f_t\}$ гамильтоновых диффеоморфизмов $M$ порождает
симплектическое расслоение.
Возьмём две копии единичного $2$-диска $D_+^2$, $D_-^2$ такие, что
$D_+^2$ имеет положительную ориентацию, а $D_-^2$ — обратную.
Определим многообразие
\[P =  M  \times D_-^2 \cup_\psi M \times D_+^2,\]
где отображение склейки $\psi$ задаётся следующим образом
\[\psi \: M \times S^1 \to M \times S^1,\quad (z, t) \mapsto (f_t z, t).\]
Ясно, что $P$ имеет естественную структуру симплектического расслоения
над $S^2$, ведь $f_t$ — симплектоморфизмы (вдобавок $S^2$
получает ориентацию по построению).

Гомотопные петли приводят к изоморфным симплектическим расслоениям, то
есть существуют гладкие изоморфизмы, сохраняющие послойную
симплектическую структуру и ориентацию.
Кроме того, это построение можно обратить — по заданному расслоению
$P \to S^2$ с выбранной тривиализацией над одной точкой можно
восстановить гомотопический класс $\gamma$.
Заметим также, что класс $\gamma = 0$ соответствует тривиальному
расслоению $S^2 \times (M, \Omega)$.
Оставляем доказательство этих утверждений читателю.
Далее будем писать $P = P(\gamma)$, где $\gamma$ — гомотопический класс
$\{f_t\}$.

Давайте посмотрим на расслоение $P(\gamma)$ для одного оборота $S^2$.
Заметим, что база, и слой этого расслоения $S^2$.
Будет полезно отождествить $S^2$ с комплексной проективной прямой $\CP^1$ и перейти к комплексной точке зрения.

Прежде чем продолжить, сделаем отступление о симплектической геометрии
комплексных проективных пространств. 
Пусть $E$ — $2n$-мерное вещественное векторное пространство с
комплексной структурой $j$ 
(здесь $j$ — линейное преобразование $E\to E$ такое, что $j^2 = -\1$),
скалярным произведением $g$ и симплектической формой~$\omega$, такими что
\[g(\xi, \eta ) = \omega(\xi, j\eta)\]
для всех $\xi, \eta \in E$ (ср. с \ref{sec:4.1}).
Конечно, используя $j$, мы можем рассматривать $E$ как комплексное
векторное пространство. 
Такая пара $(\omega, g)$ называется эрмитовой структурой на
комплексном пространстве $(E, j)$.
Рассмотрим единичную сферу
\[S = \set{\xi \in E}{g(\xi, \xi) = 1}\]
и действие окружности на $S$, определённое формулой 
\[\xi \mapsto e^{2\pi jt} \xi,\quad t \in \RR/\ZZ.\]
Орбитами этого действия являются множества вида $S \cap l$, где $l$ — комплексная прямая в $E$. 
Таким образом, пространство орбит $S/S^1$ можно канонически отождествить с \rindex{проективное пространство}\emph{комплексным проективным пространством} $\PP(E)$.
Действие сохраняет сужение $\omega$ на $\T S$, и, кроме того, касательные пространства к орбитам действия лежат в ядре этого сужения.
Обозначим через $\Omega$ проекцию $\tfrac1\pi \omega$ на $\PP(E)$.
Эта форма замкнута.
Далее, воспользовавшись элементарной линейной алгеброй, получаем, что $\Omega$ невырождена, и, таким образом, мы получаем симплектическую форму на $\PP(E)$. 

Форма $\Omega$ называется стандартной симплектической формой (или формой Фубини — Штуди) на $\PP(E)$, ассоциированной с эрмитовой структурой $(\omega, g)$ на комплексном пространстве $(E, j)$.
Приведённое выше построение является частным случаем редукции \rindex{Марсден}Марсдена — \rindex{Вайнштейн}Вайнштейна \cite{MS}, играющей ключевую роль в теории действий
групп на симплектических многообразиях.
Фактор $\tfrac1\pi$ выше выбран по следующей причине.

\begin{ex}{Упражнение}\label{9.2.A}
  Покажите, что интеграл от $\Omega$ по проективной прямой в $\PP(E)$
  равен $1$.
\end{ex}

Из теоремы \rindex{Мозер}Мозера \cite{MS} следует, что различные эрмитовы структуры
на $(E, j)$ порождают симплектические формы на
$\PP(E)$ эквивалентные с точностью до диффеоморфизма. 

\begin{ex}{Упражнение}\label{9.2.B}
  Рассмотрим пространство $\CC^n$ со стандартной эрмитовой структурой
  (см. \ref{sec:4.1}).
  Покажите, что стандартная симплектическая форма на $\PP(\CC^n ) =
  \CP^{n-1}$ инвариантна относительно группы $\PU (n) = \U (n)/S^1$ и
  что эта группа действует на $\CP^{n-1}$ гамильтоновыми
  диффеоморфизмами.
  Покажите, что 
  \[\PU (2) = \SU (2)/\{\1, -\1\}.\]
\end{ex}
 
\begin{ex}{Упражнение}\label{9.2.C}
  Рассмотрим петлю проективных унитарных преобразований $\CP^1$,
  которая в однородных координатах $(z_1 \z: z_2)$ на $\CP^1$ задается формулой
  \[(z_1 : z_2 ) \mapsto (e^{-2\pi it} z_1 : z_2 ),\quad t \in [0; 1].\]
  Покажите, что эта петля представляет собой нетривиальный элемент
  $\pi_1 (\Ham(\CP^1 ))$.
  \emph{Подсказка:} используйте канонический изоморфизм
  $\SU (2)/\{\1; -\1\}\to \SO(3)$ (см. \cite{DFN}).  
\end{ex}

У приведенного выше построения существует естественный
«параметрический» вариант, который даёт симплектические расслоения со
слоем $\CP^{n-1}$, оснащённым стандартной симплектической формой. 
Пусть $E \to S^2$ — комплексное векторное расслоение ранга $n$, и
пусть $\PP(E)$ — его проективизация. 
Каждая эрмитова структура на $E$ порождает послойную симплектическую
форму $\Omega_x$ на $\PP(E)$. 
Здесь $\Omega_x$ равна стандартной симплектической форме на слое
$\PP(E_x)$.
Различные эрмитовы структуры приводят к изоморфным симплектическим
расслоениям. 

Вернёмся теперь к симплектическому расслоению $P(\gamma)$, где
$\gamma$ — нетривиальный элемент $\pi_1 (\Ham(S^2 ))$. 
Пусть $T \to \CP^1$ — тавтологическое расслоение,
то есть его слой над комплексной прямой в $\CC^2$ это сама прямая. 
Пусть $C = \CC \times \CP^1$ — тривиальное расслоение.
 
\begin{ex}{Упражнение}\label{9.2.D}
Докажите, что симплектическое расслоение $P(\gamma)$ изоморфно $\PP(T
\oplus C)$. 
\emph{Подсказка:} представьте базу $\CP^1$ в виде объединения двух дисков 
\[D^2_- = \set{(x_0 : x_1 ) \in \CP^1}{ |x_0 /x_1 | \le 1}\]
и
\[D^2_+ = \set{(x_0 : x_1 ) \in \CP^1}{|x_1 /x_0 | \le 1}.\]
Рассмотрите координату $t$ на окружности $S^1 = \partial D^2_+ =
\partial D^2_-$ такую, что $x_1 /x_0 = e^{2\pi it}$. 
Расслоение $T \oplus C$ можно тривиализовать над $D^2_-$ и $D^2_+$ так,
что функция перехода $S^1 \times \CC^2 \to S^1 \times \CC^2$ имеет вид  
\[(t, z_1, z_2 ) \mapsto (t, e^{-2\pi it} z_1, z_2 )\]
(подробное описание тавтологического линейного расслоения дано в
\cite{GH}). 
Теперь результат следует из \ref{9.2.C}.
\end{ex}

\begin{ex}{Упражнение}\label{9.2.E}
Покажите, что $\PP(T \oplus C)$ биголоморфно эквивалентно комплексному
раздутию $\CP^2$ в одной точке. 
Расслоение получается собственным прообразом пучка прямых, проходящих
через точку раздутия. 
\end{ex}

\section{Симплектические связности}

{

\begin{wrapfigure}[12]{o}{40 mm}
\vskip-0mm
\centering
\includegraphics{mppics/pic-9}
\caption{}\label{pic-9}
\vskip0mm
\end{wrapfigure}

Пусть $p\: P\to S^2$ — симплектическое расслоение со слоем $(M,\Omega)$.
Связность $\sigma$ на $P$ (то есть поле двумерных подпространств,
трансверсальных слоям, см. рис. \ref{pic-9}) называется
симплектической, если её параллельный перенос сохраняет послойную
симплектическую структуру. 
Можно показать, что всякое симплектическое расслоение допускает
симплектическую связность, см. \cite{GLS,MS}. 
\rindex{симплектическая связность}
\rindex{кривизна симплектической связности}

\begin{ex*}{Пример}
Пусть $E \to S^2$ — комплексное векторное расслоение с эрмитовой метрикой.
Тогда каждая эрмитова связность на $E$ индуцирует симлектическую связность на проективизированном расслоении $\PP(E)\z\to S^2$.
\end{ex*}

Давайте вспомним определение кривизны связности.
Для данных $x\in S^2$ и $\xi$, $\eta \z\in \T_x S^2$, продолжим $\xi$
и $\eta$ до полей в окрестности точки $x$, 
и пусть $\tilde\xi$ и $\tilde\eta$ — горизонтальные поднятия
полученных векторных полей. 
По определению кривизна $\rho^\sigma$ в точке $x$ это 2-форма на
базе расслоения,
принимающая значения в алгебре Ли векторных полей на слое
$p^{-1}(x)$. 
Она определяется как $\rho^\sigma (\xi, \eta) \z= ([\tilde\xi,
  \tilde\eta])^\vert$, где «$\vert$» обозначает проекцию $[\tilde\xi,
  \tilde\eta]$ вдоль связности на слой $p^{-1} (x)$. 
Таким образом, если $\sigma$ — симплектическая связность, то
$\rho^\sigma (\xi, \eta)$ лежит в алгебре Ли группы
$\Symp(p^{-1}(x))$, и, поскольку $H^1 (M;\RR) = 0$, кривизна
$\rho^\sigma (\xi, \eta)$ является гамильтоновым векторным полем на
слое. 
Отождествив гамильтоновы векторные поля с нормализованными
гамильтонианами, можно думать, что $\rho^\sigma (\xi, \eta)$ есть
нормализованный гамильтониан 
\[p^{-1} (x) \to \RR.\]
Зафиксируем форму площади $\tau$ на $S^2$ такую, что $\int_{S^2} \tau
= 1$ (тут нам потребовалась ориентация $S^2$). 
Поскольку каждая 2-форма на $S^2$ в любой точке
пропорциональна $\tau$, мы можем написать
\[\rho^\sigma = L^\sigma \tau,\]
где $L^\sigma$ — функция на $P$.

}

Теория симплектических связностей была развита в работах
В. Гиймена,
Е. Лермана
и С. Штернберга \cite{GLS,MS}. 
Основными объектами этой теории являются форма сцепления
симплектической связности и класс сцепления симплектического
расслоения.  
Касательное пространство в точке $(x, z) \in P$ можно разложить по симплектической
связности $\sigma$: 
\[\T_{(x,z)} P = \T_z p^{-1} (x) \oplus \T_x S^2.\]
Определим \rindex{форма сцепления}\emph{форму сцепления} $\delta^\sigma$ связности $\sigma$
как 2-форму на $P$, заданную формулой  
\[\delta^\sigma (v \oplus \xi, w \oplus \eta) = \Omega_x (v, w) -
\rho^\sigma (\xi, \eta)(z).\] 
Здесь $z \in p^{-1} (x)$, $v$, $w \in \T_z p^{-1} (x)$ и $\xi$, $\eta
\in \T_x S^2$, а $\rho^\sigma (\xi, \eta)$ рассматривается как функция
на слое $p^{-1} (x)$. 

Оказывается, форма сцепления замкнута.
Обозначим через $c$ её класс когомологий в $H^2 (P;\RR)$.
Очевидно, что сужение $c$ на любой слой $p^{-1}(x)$ совпадает с классом $[\Omega_x]$
симплектической формы.  
Кроме того, можно доказать, что $c^{n+1} = 0$, где $2n = \dim M$.
Следующий результат показывает, что класс $c$ однозначно определяется
этими двумя свойствами. 

\begin{thm}[(См. \cite{GLS,MS})]{Теорема}\label{9.3.A}\rindex{Штернберг}
  Класс $c$ — это единственный класс когомологий в $H^2 (P;\RR)$,
  такой что $c|_{\fiber} = [\Omega_x]$ и $c^{n+1} \z= 0$.
\end{thm}

В частности, $c$ — инвариант симплектического расслоения $P$, независящий от выбора связности $\sigma$. 
Мы называем класс $c$ \rindex{класс сцепления}\emph{классом сцепления}~$P$.
Подробные доказательства всех этих результатов даны в \rindex{Гиймен}\cite{GLS} и \cite{MS}.

Следующее построение играет важную роль в нашем подходе к спектру длин.

\parbf{Построение слабого сцепления}
(\cite{GLS,MS})\rindex{Лерман}
Для достаточно малого $\epsilon > 0$ существует гладкое семейство
замкнутых 2-форм $\omega_t$ на $P$ при $t \in [0;\epsilon)$ такое, что 
\begin{itemize}
\item $\omega_0 = p^\ast \tau$
\item $[\omega_t] = tc + p^\ast [\tau]$
\item $\omega_t|_{\fiber} = t\Omega_x$
\item форма $\omega_t$ симплектическая при всех $t > 0$.
\end{itemize}

Положим \index[symb]{$\epsilon(P)$}$\epsilon(P) = \sup \{\epsilon\}$,
где точная верхняя грань берётся по всем таким деформациям.  
Эта величина измеряет, насколько сильно слабое сцепление.
Заметим, что $\epsilon(P) = +\infty$ для тривиального расслоения $P = M \times S^2$.
Таким образом, в некотором смысле $\epsilon(P)$ измеряет нетривиальность расслоения.

Есть ещё один способ измерения нетривиальности расслоения, который
работает и для расслоений с другими структурными группами. 
Его предложил \rindex{Громов}Громов \cite{G2} для унитарных векторных расслоений.
Далее мы обсудим симплектическую версию этой конструкции. 
Идея состоит в том, чтобы измерить минимально возможную норму кривизны
симплектической связности на $P$. 
Введём следующее понятие.

\begin{ex*}{Определение}
\[\chi_+ (P) = \sup_\sigma \frac1{\max_P L^\sigma}\]  — (положительная часть) \rindex{$K$-площадь}\emph{симплектической $K$-площади} $P$.
Здесь точная верхняя грань берётся по всем симплектическим связностям на $P$, а $L^\sigma$ определяется равенством $\rho^\sigma = L^\sigma \tau$.
\end{ex*}

\begin{ex*}{Упражнение}
Докажите, что обе введенные выше величины, $\epsilon(P)$ и $\chi_+
(P)$, не зависят от выбора формы площади $\tau$ на $S^2$ с $\int_{S^2}
\tau = 1$. 
\emph{Подсказка:} сначала применим теорему Мозера \cite{MS} о том, что для
любых двух таких форм, скажем $\tau_1$ и $\tau_2$, существует
диффеоморфизм $a\:S^2\to S^2$, изотопный единице и такой, что $a^\ast
\tau_2 = \tau_1$.
Затем поднимем $a$ до послойно симплектического диффеоморфизма $A$ расслоения $P$.
Это означает, что $p(A(z)) = a(p(z))$ для всех $z \in P$ и сужение $A_x$ диффеоморфизма $A$ на любой слой $p^{-1} (x)$ удовлетворяет равенству $(A_x)^\ast \Omega_{a(x)} = \Omega_x$.
Такой подъём можно построить с помощью любой симплектической связности
на $P$. 
Обратите внимание, что $A$ переводит любую деформацию слабого
сцепления для $\tau_2$ в деформацию слабого сцепления для
$\tau_1$.
Это доказывает, что $\epsilon(P)$ не зависит от выбора формы площади.
Далее, $A$ действует на пространстве симплектических связностей на $P$
по правилу $\sigma \mapsto A_\ast \sigma$.
Покажите, что 
\[\rho^{A_\ast \sigma} (a_\ast \xi, a_\ast \eta)(Az) = \rho^\sigma (\xi, \eta)(z)\]
для каждой точки $z \in P$ и любой пары векторов $\xi$, $\eta \in
\T_{p(z)} S^2$.
То, что $\chi_+(P)$ не зависит от выбора формы площади, является
простым следствием из этой формулы.
\end{ex*}

\begin{thm}[(\cite{P4})]{Теорема}\label{9.3.B}
  Пусть $P = P(\gamma)$.
  Тогда $\epsilon(P) = \chi_+ (P) \z= \frac1{\nu_+(\gamma)}.$
\end{thm}

Мы докажем более слабое утверждение, а именно что $\epsilon(P) \z\ge
\chi_+ (P)$, но его достаточно для доказательства теоремы
\ref{9.1.A}.

\parit{Доказательство неравенства $\epsilon(P) \ge \chi_+ (P)$.}
Пусть $\sigma$ — симплектическая связность.
Рассмотрим формы
\[\omega_t = p^\ast \tau + t\delta^\sigma,\]
где $\delta^\sigma$ — форма сцепления.
В точке $(x, z) \in P$ имеем 
\[\omega_{t,(x,z)}
=
t\Omega_x \oplus (-tL^\sigma (x, z)\tau) + p^{*}\tau
=
t\Omega_x \oplus (1 - tL^\sigma (x, z))\tau.\]
Ясно, что $\omega_t$ удовлетворяет первым трём свойствам в построении
слабого сцепления.
Форма $\omega_t$, является симплектической, пока $1 \z- tL^\sigma (x, z)
> 0$ или, что то же самое, пока $L^\sigma (x, z) < \frac1t$ при всех
$x$,~$z$.
Это условие означает, что 
$\max_P L^\sigma (x, z) < \frac1t$
или 
\[\frac1{\max_P L^\sigma (x, z)} > t.\]
Выберем теперь произвольное $\kappa > 0$ и симплектическую связность
$\sigma$ так, что  
\[\frac{1}{\max_P L^\sigma(x, z)} > \chi_+ (P) - \kappa.\]
Таким образом, существует деформация сцепления $\omega_t$ для $t \z\in[0;\chi_+ (P) \z- \kappa)$.
То есть $\epsilon(P) \ge \chi_+ (P) - \kappa$ для любого $\kappa > 0$,
и мы заключаем, что
\[\epsilon(P) \ge \chi_+ (P).\qedsin\]
\medskip

Заметим, что мы доказали существование деформации слабого сцепления.

\parit{Доказательство неравенства $\chi_+ (P) \ge \frac{1}{\nu_+ (\gamma)}$.}
Доказательство основано на следующем упражнении.

\begin{ex*}{Упражнение}
  Пусть $p\: P \to S^2$ — симплектическое расслоение.
  Пусть $\omega$ — такая замкнутая 2-форма на $P$, что
  $\omega|_{\fiber} = \Omega_x$.
  Положим 
  \[\sigma_{(x,z)} = \set{\xi \in \T_{(x,z)} P }{ i_\xi \omega =
    0\ \text{на}\  \T_z p^{-1} (x)}.\] 
  Покажите, что $\sigma$ определяет симплектическую
  связность на $P$. 
\end{ex*}

Пусть $\{f_t \}$, $t \in [0;1]$, — произвольная петля гамильтоновых
диффеоморфизмов, порождённая нормализованным гамильтонианом $F \z\in
\H$.
Зафиксируем полярные координаты $u \in (0; 1]$ (радиус) и $t \z\in S^1
= \RR/\ZZ$ (нормированный угол) на $D^2$. 
Возьмём такую монотонную функцию срезки $\phi(u)$, что $\phi(u) = 0$
вблизи $u = 0$ и $\phi(u) = 1$ вблизи $u = 1$. 
Положим 
\[P = M \times D_-^2 \cup_\psi M \times D_+^2,\]
где $\psi(z, t) = (f_t z, t)$.
Определим замкнутую 2-форму $\omega$ на $P$ равенством 
\[\omega=
\begin{cases}
\quad\Omega&\text{на}\ M\times D^2_+,
\\
\quad\Omega+d(\phi(u)H_t(z))\wedge dt&\text{на}\ M\times D^2_-,
\end{cases}
\]
где $H_t (z) = F (f_t z , t)$.

\begin{ex*}{Упражнение}
Докажите, что $\omega$ хорошо определена, то есть $\psi^\ast \Omega \z= \Omega + dH_t \wedge dt$. (Это можно проделать прямым вычислением.)
\end{ex*}

Вычислим кривизну $\rho^\sigma$ симплектической связности,
ассоциированной с $\omega$. 
Заметим, что $\rho^\sigma$ обращается в нуль на $D_+^2$ (поскольку
$\Omega$ индуцирует плоскую связность на $D_+^2$), поэтому остается
вычислить $\rho^\sigma$ на $D_-^2$. 
Заметим, что $\rho^\sigma = 0$ вблизи $0$ на $D_-^2$, и это хорошо, ведь мы избегаем сингулярности в нуле.
Чтобы вычислить кривизну, мы должны горизонтально поднять $\tfrac{\partial}{\partial u}$ и $\tfrac{\partial}{\partial t}$ в точках $(x, z) \in M \times D_-^2$.

Пусть горизонтальный подъём векторного поля $\tfrac{\partial}{\partial u}$ имеет вид 
\[\widetilde{\tfrac{\partial}{\partial u}}
=
\tfrac{\partial}{\partial u} + v
\quad\text{для некоторого}\quad
v \in \T_z p^{-1} (x).\]
По определению связности, $\omega\left(\widetilde{\frac{\partial}{\partial u}},w\right)=0$ для всех $w\in \T_zM$, значит 
\[0=\omega\left(\widetilde{\tfrac{\partial}{\partial u}}, w\right) = \omega(v, w) = \Omega(v, w)\]
для всех $w \in \T_z p^{-1} (x)$.
Из невырожденности $\Omega$ следует, что $v = 0$ и, значит, 
\[\widetilde{\tfrac{\partial}{\partial u}}
=
\tfrac{\partial}{\partial u}\]
Полагая, что
\[\widetilde{\tfrac{\partial}{\partial t}}
=
\tfrac{\partial}{\partial t}+v
\quad\text{для другого}\quad
v \in \T_z p^{-1} (x),\]
как и раньше, получаем, что 
\[0
=
\omega\left(\widetilde{\tfrac{\partial}{\partial t}}, w\right)
= 
\omega(\tfrac{\partial}{\partial t}+v, w) = \Omega(v, w)- d(\phi(u)H_t)(w)
\]
для всех $w \in \T_z p^{-1} (x)$.
Итак, 
$i_v \Omega = d(\phi(u)H_t)$.
Следовательно, 
$v \z= - \sgrad \phi(u)H_t \z= -\phi(u) \sgrad H_t$,
и можно сделать вывод, что 
\[\widetilde{\tfrac{\partial}{\partial t}}
=
\tfrac{\partial}{\partial t}-\phi(u) \sgrad H_t.\]
Вычислив кривизну, получим 
\[\rho^\sigma
\left(\tfrac{\partial}{\partial t},\tfrac{\partial}{\partial u}\right)
=
\left[\tfrac{\partial}{\partial t}-\phi(u)\sgrad H_t,\frac{\partial}{\partial u}\right]^\vert
=
\phi'(u)\sgrad H_t.\]
(На самом деле коммутатор уже вертикальное векторное поле.)
Переходя от векторных полей к гамильтонианам в определении кривизны, мы видим, что
\[\rho^\sigma\left(\tfrac{\partial}{\partial t},\tfrac{\partial}{\partial u}\right)
=\phi'(u) H_t.
\]

Зафиксируем $\kappa > 0$ и выберем форму площади на $D_-^2$ вида $(1 \z- \kappa)dt \wedge du$ (напомним, что у $D_-^2$ обратная ориентация).

Продолжим её до формы площади на $D_+^2$ такой, что $\area (D_+^2) = \kappa$.
Полученную форму обозначим $\tau$.
Мы видим, что $\rho^\sigma = L^\sigma (u, t, z)\tau$, где 
\[
L^\sigma(u,t,z)=
\begin{cases}
\quad0&\text{на\ } M\times D^2_+\,,
\\
\quad\frac{\phi'(u)H_t(z)}{1-\kappa}&\text{на\ } M\times D^2_-\,.
\end{cases}
\]
Наконец, выберем $\phi$ так, что $\phi' (u) \le 1 + \kappa$, и
$\{f_t\}$ так, что $\max_z F_t \z= \max_z H_t \le \nu_+ (\gamma) +
\kappa$. 
Получаем
\[\max_P L^\sigma \le  \frac{1+\kappa}{1-\kappa} (\nu_+ (\gamma) + \kappa),\]
значит 
\begin{align*}
\chi_+ (P) &= \sup_\sigma \frac1{\max_PL^\sigma}\ge
\\
&\ge\frac{1-\kappa}{(1 + \kappa)(\nu_+ (\gamma) + \kappa)}. 
\end{align*}
Поскольку $\kappa$ произвольно, мы доказали, что 
\[\chi_+ (P) \ge \frac{1}{\nu_+ (\gamma)}.\]
\qeds

\section{Приложение к спектру длин}\label{sec:9.4}

Последний шаг в доказательстве теоремы \ref{9.1.A} состоит в следующей
оценке на $\epsilon(P)$, которая обсуждается в следующей главе.

\begin{thm}{Теорема}\label{9.4.A}
Пусть $P = P(\gamma)$, где $\gamma$ — поворот $S^2$ на 1 оборот.
Тогда
\[\epsilon(P) \le 2.\]
\end{thm}

\parit{Доказательство \ref{9.1.A}.}
Мы знаем, что $\nu_+ (\gamma) \le \tfrac12$.
С другой стороны, из \ref{9.3.B} и \ref{9.4.A} получаем, что
\[2 \ge \epsilon(P) \ge \chi_+ (P) \ge\frac1{\nu_+(\gamma)},\]
поэтому $\nu_+ (\gamma) = \tfrac12$ и $\epsilon(P) = \chi_+ (P) = 2$.
\qeds

\begin{ex}{Упражнение}
Рассмотрим голоморфные линейные расслоения $C$ и $T$ над $\CP^1$ —
соответственно тривиальное и тавтологическое (см. \ref{sec:9.2}). 
Пусть $\nabla_C$ — естественная плоская связность на~$C$.
Существует единственная связность, скажем $\nabla_T$, на $T$,
     которая сохраняет полученную из $\CC^2$ каноническую эрмитову
     структуру на $T$, и такая, что её $(0, 1)$-часть совпадает с
     оператором $\bar\partial$ (см. \cite{GH}). 
Рассмотрим симплектическое расслоение $P = \PP(T \oplus C)$ и
обозначим через $\sigma$ связность на $P$, полученную из $\nabla_T
\oplus \nabla_C$. 
Очевидно, что $\sigma$ симплектична (на самом деле её параллельный
перенос сохраняет как симплектическую, так и комплексную структуру на
слоях). 
Пусть $\tau$ — форма площади Фубини — Штуди на $\CP^1$,
определенная в \ref{sec:9.2}. 
Докажите, что 
\[\frac1{\max_P L^\sigma}=2,\]
где $\rho^\sigma = L^\sigma \tau$.
В частности, $\sigma$ является связностью с минимально возможной
кривизной, где кривизна «измеряется» относительно $\tau$.  
\end{ex}

\chapter[Деформации симплектических форм]{Деформации симплектических форм и псевдоголоморфные кривые}

В этой главе мы получим верхнюю оценку на параметр сцепления
\ref{9.4.A} и, таким образом, завершим вычисление спектра длин $\Ham(S^2)$.
Эта оценка оказывается частным случаем более общей задачи о
деформациях симплектических форм (см. \cite{P7}). 
В её решении будет использована громовская теория псевдоголоморфных кривых.

\section{Задача деформации}

\begin{wrapfigure}[8]{r}{30 mm}
\vskip-6mm
\centering
\includegraphics{mppics/pic-10}
\caption{}\label{pic-10}
\vskip0mm
\end{wrapfigure}

Пусть $(P, \omega)$ — замкнутое симплектическое многообразие и $l$ —
луч в $H^2 (P;\RR)$ с началом в $[\omega]$, см. рис. \ref{pic-10}. 

\begin{ex*}{Задача}
Как далеко можно деформировать $\omega$ так, чтобы класс когомологий двигался вдоль $l$ и форма оставалась симплектической? 
\end{ex*}

Обратите внимание, что добавление малой замкнутой 2-формы к $\omega$ оставляет её симплектической — вопрос в том, как далеко можно зайти.
В деформации сцепления форма $\omega_0=p^\ast \tau$ вырождена, но форма $\omega_t$ в классе $[p^\ast \tau ] + tc$ становится симплектической при малых $t$.
В этом случае нам удастся пройти по лучу в направлении класса сцепления $c$ ровно на $\epsilon(P)$.

Приведём пример препятствия к существованию бесконечной деформации.
Предположим, что $\dim P = 4$ и что $\Sigma \subset P$ — вложенная 2-сфера такая, что $\omega|_{\T \Sigma}$ — форма площади (иными словами, $\Sigma$ — симплектическое подмногообразие в $P$).
Через $(A, B)$ будет обозначаться индекс пересечения классов гомологий $A$ и $B$.

\begin{ex*}{Определение}
Пусть $\Sigma \subset P^4$ — симплектическая вложенная сфера.
Если $([\Sigma], [\Sigma]) = -1$, то $\Sigma$ называется \rindex{исключительная сфера}\emph{исключительной сферой}.
\end{ex*}

\begin{thm}[(\cite{McD1})]{Теорема}\label{10.1.A}
Пусть $\Sigma \subset (P^4, \omega)$ — исключительная сфера.
Пусть $\omega_t$, $t \in [0;1]$, — симплектическая деформация $\omega$.
Тогда $([\omega_1 ], [\Sigma]) > 0$.
\end{thm}

Иными словами, гиперплоскость $(x, [\Sigma]) = 0$ в $H^2 (P;\RR)$ образует непроницаемую стенку для симплектических деформаций.

Доказательство использует теорию псевдоголоморфных кривых.
Ниже мы дадим набросок доказательства, а также
приложение к доказательству \ref{9.4.A}. 

\section[\texorpdfstring{И снова $\bar\partial$-уравнение}{И снова ∂-уравнение}]{И снова $\bm{\bar\partial}$-уравнение}

\rindex{почти комплексная структура}\emph{Почти комплексная структура} $j$ на многообразии $P$ это поле
эндоморфизмов $\T P \to \T P$ таких, что $j^2 = -\1$. 
Важный класс примеров приходит из комплексной алгебраической геометрии.
Каждое комплексное многообразие (то есть многообразие с атласом у
которого голоморфны отображения склеек) имеет каноническую почти
комплексную структуру $\xi \to \sqrt{-1}\xi$. 
Возникающие таким образом почти комплексные структуры называются
\emph{интегрируемыми}. 
Важно, что интегрируемость эквивалентна
обнулению некоторого тензора, связанного с почти комплексной
структурой (\cite{NN}). 
Как следствие получается, что всякая почти
комплексная структура на (вещественной) поверхности интегрируема.%
\footnote{Поскольку почти комплексная структура на поверхности это конформная структура плюс ориентация,
это утверждение также следует из теоремы об униформизации. — \textit{Прим. ред.}}
Более того, на 2-сфере любые две почти комплексные структуры диффеоморфны (это классическое утверждение — так называемая теорема об униформизации \cite{AS}).

\begin{ex*}{Определение}
Пусть $(P, \omega)$ — симплектическое многообразие.
Почти комплексная структура $j$ называется \rindex{совместимая комплексная структура}\emph{совместимой} с $\omega$, если билинейная форма $g(\xi,
\eta) \z= \omega(\xi, j\eta)$ определяет риманову метрику на~$P$. 
\end{ex*}

\begin{ex}{Упражнение}\label{10.2.A}
Пусть $(P, \omega)$ — вещественная симплектическая
поверхность и пусть $j$ — почти комплексная структура на $P$. 
Тогда либо $j$, либо $-j$ совместима с $\omega$.
\end{ex}

Пусть $(P, \omega)$ — симплектическое многообразие и $j$ —
совместимая интегрируемая почти комплексная структура на $P$.
Тройка $(P, \omega, j)$ называется \rindex{кэлерово
  многообразие}\emph{кэлеровым многообразием}.  
Например, $\CP^n$ со стандартными $\omega$ и $j$ является кэлеровым
многообразием (см. \ref{sec:9.2}). 
Таким образом, каждое комплексное подмногообразие в $\CP^n$ кэлерово
относительно индуцированной структуры. 
И наоборот, если $(P, \omega, j)$ — замкнутое кэлерово многообразие
такое, что класс когомологий $[\omega]$ целочисленен, то $(P, j)$
голоморфно вкладывается в $\CP^n$ для некоторого $n$ (это знаменитая
теорема Кодаиры, см. \cite{GH}). 

Глубокое наблюдение, сделанное \rindex{Громов}Громовым \cite{G1}, состоит в том
что некоторые мощные методы алгебраической геометрии можно обобщить на
\rindex{квазикэлерово многообразие}\emph{квазикэлеровы многообразия}
$(P, \omega, j)$. 
Здесь $j$ — совместимая почти комплексная структура, возможно не интегрируемая. 
Примечательно, что теория голоморфных кривых без существенных
изменений распространяется на неинтегрируемый случай. 
Важность такого обобщения обусловлена тем, что
каждое симплектическое многообразие допускает совместимую почти
комплексную структуру.
При этом существуют симплектические многообразия, которые не допускают
кэлеровой структуры~\cite{MS}. 

\begin{ex}[\cite{MS}]{Упражнение}\label{10.2.B}
Пусть $E$ — чётномерное линейное пространство с невырожденной
кососимметричной билинейной формой~$\omega$. 
Покажите, что пространство комплексных структур $j\: E \z\to E$, $j^2
= -\1$, совместимых с $\omega$ стягиваемо. 
\end{ex}

Таким образом, совместимая почти комплексная структура на
симплектическом многообразии — это сечение расслоения, слои которого
стягиваемы. 
Следовательно, такие структуры существуют и, кроме того, образуют
стягиваемое пространство. 
По этой же причине различные задачи продолжения, связанные с почти
комплексными структурами, допускают положительное решение. 

\begin{ex}{Упражнение}\label{10.2.C}
Пусть $(P, \omega)$ — четырёхмерное симплектическое многообразие, и
пусть $\Sigma \subset P$ — симплектическое подмногообразие. 
Предположим, что $\Sigma$ снабжена почти комплексной структурой $j$,
согласованной с $\omega|_{\T \Sigma}$. 
Покажите, что $j$ продолжается до совместимой почти комплексной
структуры на $P$. 
\end{ex}

Теперь перейдём к теории псевдоголоморфных кривых на квазикэлеровых многообразиях.

\begin{ex*}{Определение}
Отображение $\phi\:(S^2, i) \to (P, j)$ является псевдоголоморфной
(или $j$-голоморфной) кривой, если  
$\phi_\ast \circ i = j \circ \phi_\ast$.
\end{ex*}

\begin{ex}{Упражнение}\label{10.2.D}
Покажите, что на комплексном многообразии приведённое выше определение
эквивалентно обычному уравнению Коши — Римана.
\end{ex}

Определим 
\[\bar\partial\phi=\frac12(\phi_\ast+j\circ\phi_\ast\circ i)\]

\begin{ex}[(ср. \ref{4.1.A})]{Упражнение}\label{10.2.E}
Для заданных $(P, \omega, j, g)$ и $j$-го\-ло\-морф\-ной кривой $\phi
\: (S^2, i) \to (P, j)$ покажите, что $\area_g (\phi(S^2)) \z=
\int_{S^2}\phi^\ast \omega$. 
Более того, если $\phi$ непостоянна, то $\int_{S^2} \phi^\ast \omega
> 0$.  
\end{ex}

В дополнение сказанному, если $\phi$ — вложение, то $\phi(S^2)$ —
симплектическое подмногообразие в $P$. 
Действительно, ограничение симплектической формы совпадает с формой
римановой площади. 

Пусть $(P, \omega)$ — симплектическое многообразие.
Выберем почти комплексную структуру $j$, совместимую с $\omega$.
Тогда касательное расслоение $\T P$ получит структуру комплексного
векторного расслоения. 
Поскольку пространство совместимых $j$ связно, соответствующие
характеристические классы не зависят от выбора $j$. 
Обозначим через \index[symb]{$c_1$}$c_1$ первый класс Черна расслоения $\T P$
относительно любой согласованной почти комплексной структуры. 

\begin{ex}[(формула присоединения)]{Упражнение}\label{10.2.F}
Пусть $(P^4, \omega)$ — симплектическое многообразие с согласованной
почти комплексной структурой $j$. 
Пусть $\Sigma \subset P$ — вложенная $j$-голоморфная сфера.
Покажите, что
\[1 +\tfrac12 (([\Sigma], [\Sigma]) - c_1 (\Sigma)) = 0.\]
\emph{Подсказка:} используйте, что $([\Sigma], [\Sigma])$ —
самопересечение в комплексном нормальном расслоении $\nu_\Sigma$ и
$\T_\Sigma P = \T\Sigma \oplus \nu_\Sigma$. 
\end{ex}

\begin{thm*}{Следствие}
Пусть $\Sigma$ — симплектически вложенная сфера.
Тогда $c_1(\Sigma) = 1$ в том и только в том случае, когда $([\Sigma], [\Sigma]) = -1$. 
\end{thm*}

\parit{Доказательство.}
Поскольку $\Sigma$ симплектически вложена, существует
$\omega$-совместимая почти комплексная структура $j$ такая, что
$\Sigma$ $j$-голоморфна (см. \ref{10.2.C}). 
Теперь утверждение следует из формулы присоединения \ref{10.2.F}. 
\qeds

\section{Приложение к сцеплению}\label{sec:10.3}

В этом разделе мы выводим \ref{9.4.A} из \ref{10.1.A}.
Напомним, что мы изучаем деформацию сцепления расслоения $P(T \oplus C) \to \CP^1$, где $T$ и $C$ — соответственно тавтологическое и тривиальное расслоения. 

\begin{ex*}{Упражнение}
Пусть $E$ — комплексное векторное пространство, и пусть $l \in P(E)$
— прямая в $E$. 
Покажите, что $\T_l P(E)$ канонически изоморфно $\Hom(l, E/l) = l^\ast \otimes E/l$.
\end{ex*}

Прежде всего, мы хотим вычислить кольцо (ко)гомологий расслоения $P \z= P(T \oplus C)$.
Обозначим через $[F]$ гомологический класс слоя, а через $\Sigma$
сечение, соответствующее подрасслоению $0\oplus C$ ранга $1$. 
Ясно, что $([F], [F]) = 0$ и $([F], [\Sigma]) = 1$.
Для вычисления $([\Sigma], [\Sigma])$, обратим внимание на то, что
нормальное расслоение $\nu_\Sigma$ — это просто ограничение на
$\Sigma$ касательного расслоения к слоям. 
Из приведённого выше упражнения следует, что 
\[\nu_\Sigma = \Hom(C, T \oplus C/C) = \Hom(C, T) = C^\ast \otimes T = T.\]
Таким образом, $c_1 (\nu_\Sigma) = -1$. 
Так как $([\Sigma], [\Sigma])$ это самопересечение в нормальном расслоении,
мы заключаем, что $([\Sigma], [\Sigma]) \z= -1$.

Пусть $\omega_t$ — деформация сцепления.
Из теоремы Мозера \cite{MS} легко следует, что для любого достотчно малого
$t$ форма $\omega_t$ симплектоморфна кэлеровой форме относительно
стандартной комплексной структуры на $P(T \oplus C)$. 
Не умаляя общности можно считать, что $\omega_t$ — кэлерова
форма при малых $t$. 
Поскольку $\Sigma$ — голоморфное сечение расслоения $P$, мы получаем, что
$\Sigma$ симплектично и вложено. 
Учитывая, что $([\Sigma], [\Sigma]) = -1$, получаем, что $\Sigma$ —
исключительная сфера. 
Таким образом, из теоремы \ref{10.1.A} следует, что $([\omega_t],
[\Sigma]) > 0$ для всех~$t$. 

Напомним, что $[\omega_t] = p^\ast [\tau] + tc$, где $[\tau]$ —
образующая группы $H^2 (S^2;\ZZ)$, соответствующая ориентации сферы, а $c$ — класс сцепления.  
Классы $p^\ast [\tau]$ и $c$ однозначно определяются следующими соотношениями: 
\begin{align*}
(p^\ast [\tau], [F]) &= 0,
&
(p^\ast [\tau], [\Sigma]) &= 1,
\\
(c, [F]) &= 1,
&
c^2 &= 0.
\end{align*}
Мы будем использовать двойственность Пуанкаре для отождествления гомологий и когомологий.

\begin{ex*}{Упражнение}
Докажите, что $p^\ast [\tau] = [F]$ и $c = [\Sigma] + \tfrac12 [F]$.
\end{ex*}

Таким образом, в силу теоремы \ref{10.1.A}
\[([\omega_t], [\Sigma]) = ([F] + t[\Sigma] + \tfrac t2[F], [\Sigma]) = 1 - t + \tfrac t2= 1 -\tfrac t2>0,\]
и, значит, $t < 2$.
Поскольку это верно для любой деформации сцепления, получаем $\epsilon(P)\le2$,
что завершает доказательство.
\qeds

\section{Псевдоголоморфные кривые}\label{sec:10.4}

Здесь мы кратко изложим громовскую теорию \rindex{псевдоголоморфная кривая}псевдоголоморфных кривых \cite{G1,AL}. 
Пусть $(P^{2n}, \omega)$ — симплектическое многообразие и пусть $A
\in H_2 (P; \ZZ)$ — примитивный класс,
то есть $A$ нельзя представить в виде $kB$, где $k > 1$ — целое
число и $B \in H_2 (P; \ZZ)$. 
В частности, $A \ne 0$.
Пусть $\J$ — пространство всех $\omega$-совместимых почти
комплексных структур на $P$, а $\mathcal{N}$ — пространство всех
гладких отображений $f\: S^2 \to P$ таких, что $[f] = A$. 
Определим $\mathcal{X} \subset \mathcal{N} \times \J$ как 
\[\mathcal{X}
=
\set{(f,j)}{f \in \mathcal{N},\  j \in \J \ \text{и}\  \bar
  \partial_j f = 0}.\] 
Тогда $\mathcal{X}$ — гладкое подмногообразие в $\mathcal{N} \times
\J$ и проекция $\pi\:\mathcal{X}\z\to\J$~—
оператор Фредгольма, то есть 
$\Ker \pi_\ast$ и $\Coker \pi_\ast$ конечномерны.
Индекс оператора $\pi$ удовлетворяет условию 
\[\Index \pi_\ast
\df
\dim (\Ker \pi_\ast) - \dim (\Coker \pi_\ast) 
= 
2(c_1 (A) + n).
\]
Фредгольмовость оператора $\pi$ позволяет применить бесконечномерную версию теоремы Сарда \cite{Sm}.

\begin{figure}[ht!]
\vskip0mm
\centering
\includegraphics{mppics/pic-11}
\caption{}\label{pic-11}
\vskip0mm
\end{figure}

Пусть $j_0$ и $j_1$ — регулярные значения проекции $\pi$ (то есть оператор $\pi_\ast$
сюръективен для всех $x \in \pi^{-1} (j_k)$, $k = 0, 1$). 
Тогда $\pi^{-1} (j_0)$ и $\pi^{-1} (j_1)$ — гладкие подмногообразия.
Для пути $\gamma$ общего положения, соединяющего $j_0$ с $j_1$,
прообраз $\pi^{-1}(\gamma)$ является гладким подмногообразием
размерности $\Index (\pi) + 1$ и $\partial\pi^{-1} (\gamma) = \pi^{-1}
(j_0) \cup \pi^{-1} (j_1)$ (см. рис. \ref{pic-11}). 

Решающую роль играет свойства компактности $\pi^{-1} (\gamma)$.
Прежде всего заметим, что слой $\pi^{-1} (\gamma)$ сам по себе не
компактен, ведь на нём действует некомпактная группа $\PSL(2, \CC)$. 
Здесь $\PSL(2, \CC)$ — группа конформных преобразований $(S^2, i)$ и
она действует на $\pi^{-1} (j_0)$, ведь 
если отображение $h\: S^2 \z\to S^2$ конформно и $(f, j_0) \in \pi^{-1} (j_0)$, то
и $(f \circ h, j_0) \in \pi^{-1} (j_0)$. 

\begin{ex*}{Упражнение}
Покажите, что это действие свободно.
(Воспользуйтесь тем, что класс $A$ примитивен.)
\end{ex*}

\begin{thm*}{Теорема Громова о компактности}\rindex{теорема о компактности}\rindex{Громов}
Либо пространство модулей $\pi^{-1} (\gamma)/\PSL(2, \CC)$ компактно, либо существует такое семейство $(f_k, j_k) \z\in \pi^{-1} (\gamma)$, что $j_k \to j_\infty$ и $f_k$ «сходится» к составной $j_\infty$-голоморфной кривой в классе $A$.
\end{thm*}

Мы не определяем сходимость.
Для нас важно, что если пространство модулей $\pi^{-1} (\gamma)/\PSL(2, \CC)$ не компактно, то существует \rindex{составная кривая}\emph{составная} $j$-голоморфная кривая в классе $A$ для некоторого $j \in \J$.
Составная кривая определяется следующим образом.
Пусть $A \z= A_1 +\dots+ A_d$ при $d > 1$ — такое разложение $A$, что $A_k \ne 0$ для всех $k$, и пусть $\phi_k \: S^2 \to P$ — $j$-голоморфные кривые в классах $A_k$, $k = 1,\dots,d$.
Мы говорим, что эти данные определяют составную кривую в классе~$A$.
Объединение $\phi_k (A_k)$ называется его образом и обычно считается связным.

Рассмотрим теперь следующую ситуацию.
\begin{itemize}
\item Множество $\J'$ тех $j \in \J$, у которых есть составная кривая в классе $A$, имеет коразмерность не менее $2$ в $\J$.

\item точка $j_0 \in \J\setminus \J'$ является регулярным значением проекции $\pi$ и $\pi^{-1} (j_0)/\PSL(2, \CC)$ (которое, таким образом, является компактным многообразием без края) не кобордантно нулю (то есть не ограничивает компактное многообразие).
 
\end{itemize}
Например, точка не кобордантна нулю, ведь она не ограничивает компактное многообразие.

Тогда для любого регулярного значения $j_1 \in \J \setminus \J'$ множество $\pi^{-1} (j_1)/\PSL(2, \CC)$ не пусто.
В самом деле, если бы оно было пусто то, соединив $j_0$ с $j_1$ путём общего положения $\gamma$ в $\J \setminus\J'$, получим
\[\partial(\pi^{-1}(\gamma)/\PSL(2,\CC))
=
\pi^{-1}(j_0)/\PSL(2,\CC),
\]
а это противоречит тому, что $\pi^{-1} (j_0)/\PSL(2, \CC)$ не кобордантно нулю.
Рис. \ref{pic-12} иллюстрирует противоречие.

\begin{figure}[ht!]
\vskip0mm
\centering
\includegraphics{mppics/pic-12}
\caption{}\label{pic-12}
\vskip0mm
\end{figure}

Это явление соответствует \rindex{принцип продолжения}принципу продолжения, с которым мы уже встречались в
случае псевдоголоморфных дисков в разделе \ref{sec:4.2}: 
либо решения $\bar\partial$-уравнения сохраняются, либо происходит выдувание.

\section{Сохранение исключительных сфер}

Приведём набросок доказательства теоремы \ref{10.1.A}.
Мы полагаем, что $(P^4, \omega)$ — симплектическое многообразие и $\Sigma \subset P$ — исключительная сфера с $A = [\Sigma]$.
Пусть $\omega_t$, $t \in [0;1]$, — деформация симплектической формы $\omega = \omega_0$.

\parbf{Шаг 1.}
Выберем такую $\omega$-совместимую почти комплексную структуру $j_0$, что $\Sigma$ является $j_0$-голоморфной, и расширим $j_0$ до однопараметрического семейства $j_t$ такого, что $j_t$ является $\omega_t$-совместимой (ср. \ref{10.2.C}).

\parbf{Шаг 2.}
Теория, описанная в предыдущем разделе, без изменений работает для
множества всех $j$, совместимых с симплектическими структурами в
данном классе деформации. 
И значит
\[\Index\pi = 2(c_1 (A) + n) = 2(1 + 2) = 6.\]
Однако $\dim_\RR\PSL(2, \CC) = 6$, поэтому $\dim \pi^{-1}
(j_0)/\PSL(2, \CC) = 0$ при условии, что значение $j_0$ регулярно. 

\parbf{Шаг 3.}
В размерности 4 мы видим следующее.
Два различных ростка $j$-голоморфных кривых всегда пересекаются с
положительным индексом в общей точке. 
Это утверждение хорошо известно если $j$ интегрируема.
В неинтегрируемом случае оно просто следует из линейной алгебры при
условии, что пересечение трансверсально. 
Однако доказательство требует тонкого локального анализа
нетрансверсальных пересечений. 
Как следствие, в $A$ существует единственная $j_0$-голоморфная кривая
— если бы была другая, скажем $\Sigma'$, то $([\Sigma], [\Sigma'])
= ([\Sigma], [\Sigma]) = -1$, что противоречит тому, что $([\Sigma],
[\Sigma']) \ge 0$. 
Мы заключаем, что $\pi^{-1} (j_0)/\PSL(2, \CC)$ состоит ровно из одной
точки. 

\parbf{Шаг 4.}
Выберем регулярные структуры $j_0$ и $\{j_t\}$.
Мы утверждаем, что в общем случае составные кривые не появляются.
Действительно, положим $A = A_1 +\dots + A_d$, $d > 1$ и $c_1 (A) = 1
= c_1 (A_1) +\dots + c_1 (A_d)$. 
Таким образом, по крайней мере один из классов Черна в правой части неположителен.
Не умаляя общности, можно предположить, что $c_1 (A_1) \le 0$ и что
$A_1$ представляется $j_s$-голоморфной кривой для некоторого $s \in
[0;1]$. 
Обозначим через $\pi_{A_1}$ проекцию для класса $A_1$ (см. начало \ref{sec:10.4}).
Мы видим, что $\pi_{A_1}^{-1} (\gamma)/\PSL(2, \CC)$ непусто.
С другой стороны, 
\[\dim \pi_{A_1}^{-1} (\gamma)/\PSL(2, \CC) = 2(c_1 (A_1) + 2) - 6 + 1 \le -1,\]
что невозможно.
Таким образом, в общем случае выдувания не происходит (это отражает то, что $\codim \J' \ge 2$). 
Следовательно, $A$ представляется $j_1$-голоморфной кривой и, следовательно, $([\omega_1], [A]) > 0$ (см. \ref{10.2.E}).

\chapter[Эргодическая теория]{Приложение к эргодической теории}\label{chap:11}

В настоящей главе мы обсудим асимптотический геометрический инвариант,
связанный с фундаментальной группой группы $\Ham(M, \Omega)$, и применим его
в классической эргодической теории (см. \cite{P9}). 

\section[Гамильтоновы петли и динамика]{Гамильтоновы петли как динамические объекты}\label{sec:11.1}

Пусть $(M,\Omega)$ — замкнутое симплектическое многообразие.
Для заданного иррационального числа $\alpha$ и гладкой петли $h\:S^1\z\to\Ham(M,\Omega)$ можно определить отображение \rindex{косое произведение}\emph{косого произведения} $T_{h,\alpha} \: M \times S^1 \z\to M\times S^1$ как $T_{h,\alpha}(y,t)=(h(t)y,t+\alpha)$.
Наша цель — связать геометрию и топологию гамильтоновых петель с динамикой соответствующих косых произведений.%
\footnote{В этой главе петли свободные, мы не предполагаем, что $h(0) = \1$.}

Приведённое выше определение является частным случаем гораздо более
общего понятия косого произведения  \cite[с. 231]{CFS}, которое
интенсивно изучалось несколько десятилетий. 
Есть по крайней мере две важные причины интереса к этому понятию.
Во-первых, оно служит основой для математических моделей случайной
динамики (см. обзор \cite{Ki}). 
Во-вторых, даёт нетривиальные примеры систем с интересными
динамическими свойствами. 

Нас будет интересовать строгая эргодичность.
Напомним, что гомеоморфизм $T$ компактного топологического пространства $X$ \rindex{строгая эргодичность}\emph{строго эргодичен}, если он имеет ровно одну инвариантную борелевскую вероятностную меру, скажем $m$, которая, кроме того, положительна на непустых открытых подмножествах. 
Строго эргодические гомеоморфизмы эргодичны и обладают рядом
дополнительных замечательных свойств. 
Упомянем одно свойство, которое сыграет решающую роль.
А именно, если $T$ строго эргодичен, то для произвольной непрерывной
функции $F$ на $X$ средние по времени
$\tfrac1N\sum_{i=0}^{N-1}F(T^ix)$ равномерно сходятся к среднему по
пространству $\int_XFdm$ и, в частности, сходятся при всех $x \in X$. 
Заметим, что в общем случае для эргодических преобразований такая сходимость
имеет место только {}\emph{почти везде}. 
Контраст между «везде» и «почти везде» становится совсем чётким, если убедиться в наличии чисто топологических препятствий к строгой эргодичности. 
Например, 2-сфера не допускает строго эргодических гомеоморфизмов.
Действительно, из теоремы Лефшеца следует, что каждый гомеоморфизм $S^2$ имеет либо неподвижную точку, либо орбиту периода $2$, и мы видим, что инвариантная мера, сосредоточенная на такой орбите, противоречит определению строгой эргодичности. 
В главе~\ref{sec:11.2} мы опишем более хитрое препятствие к строгой эргодичности, вытекающее из симплектической топологии. 

Мы говорим, что петля $h\:S^1\to\Ham(M,\Omega)$ \rindex{строго эргодичная петля}\emph{строго эргодична}, если косое произведение $T_{h,\alpha}$ строго эргодично
для некоторого $\alpha$.%
\footnote{Заметим, что каждый класс $T_{h,\alpha}$ сохраняет
  каноническую меру на $M \times S^1$, индуцированную симплектической
  формой.
Таким образом, в нашей постановке строгая эргодичность означает, что
эта мера является (с точностью до множителя) единственной инвариантной
мерой.} 
На этом языке наш центральный вопрос можно сформулировать следующим образом.

\begin{ex*}{Вопрос}
Какие гомотопические классы $S^1 \to \Ham(M, \Omega)$ могут быть
представлены строго эргодическими петлями? 
\end{ex*}

Приведём пример в котором на этот вопрос можно полностью ответить.
Пусть $M_\ast$ — раздутие комплексной проективной плоскости $\CP^2$
в одной точке. 
Выберем кэлерову симплектическую структуру $\Omega_\ast$ на $M_\ast$,
интеграл которой по прямой общего положения равен 1, а интеграл по
исключительному дивизору равен $\tfrac13$. 
Ввиду \ref{9.2.E}, можно также думать, что $M_\ast = \PP(T \oplus
C)$, где $T$ и $C$ — тавтологическое и тривиальное голоморфное
линейные расслоения над $\CP^1$ соответственно. 
На этом языке исключительный дивизор соответствует нашему старому
приятелю — сечению $\Sigma$ из главы \ref{sec:10.3}, а общая прямая гомологична
сумме $\Sigma$ со слоем. 
Периоды симплектической формы выбраны таким образом, чтобы класс её
когомологий был кратен первому классу Черна многообразия $M$
(см. определение \ref{11.3.A}). 
Легко видеть, что $(M_\ast, \Omega_\ast)$ допускает эффективное
гамильтоново действие унитарной группы $\U(2)$, то есть существует
мономорфизм $i\: \U(2) \to \Ham(M_\ast, \Omega_\ast)$. 
Фундаментальная группа $\U(2)$ равна $\ZZ$.
Недавно \rindex{Абреу}Абреу и Макдафф \cite{AM} доказали, что включение $\pi_1(U(2))
\to \pi_1 (\Ham(M_\ast, \Omega_\ast))$ является изоморфизмом и,
следовательно, $\pi_1(\Ham(M_\ast, \Omega_\ast)) = \ZZ$. 
Насколько мне известно, это простейший пример симплектического
многообразия с $\pi_1 (\Ham) = \ZZ$. 

\begin{thm}{Теорема}\label{11.1.A}
Только тривиальный класс $0\z\in\pi_1(\Ham(M_\ast, \Omega_\ast))$ представим строго эргодической
петлёй.
\end{thm}

Доказательство теоремы разбито на две части.
Прежде всего необходимо установить существование стягиваемых строго
эргодических петель. 
Это можно сделать чисто эргодическими методами в достаточно общем случае.
Мы отсылаем читателя к \cite{P9} за подробностями.
Во-вторых, нужно доказать, что каждый класс $\gamma \ne 0$ не может
быть представлен строго эргодической петлёй. 
Препятствие исходит из геометрии группы $\Ham(M, \Omega)$.
Мы обсудим это подробнее.

\section{Асимптотической спектр длин}\label{sec:11.2}
\rindex{асимптотической спектр длин}

Определим \rindex{асимптотическая норма}\emph{асимптотическую норму} элемента
$\gamma\in\pi_1(\Ham(M, \Omega))$ как
\[\nu_\infty(\gamma)=\lim_{k\to+\infty}\frac{\nu(k\gamma)}{k},\]
где $\nu$ — норма, введённая в \ref{sec:7.3}. 
Это понятие похоже на асимптотический рост $\mu$ определённый в \ref{sec:8.2}.
Предел существует, поскольку последовательность $\nu(k\gamma)$ субаддитивна.

\begin{thm}{Теорема}\label{11.2.A}
Пусть $\gamma\in\pi_1(\Ham(M,\Omega))$ — класс, представленный
гладкой строго эргодической петлёй.
Тогда асимптотическая норма $\nu_\infty (\gamma)$ обращается в нуль.
\end{thm}

\parit{Доказательство.}
Доказательство основано на процедуре асимптотического сокращения
кривой в духе \ref{sec:8.3}. 
Пусть $h\: S^1 \to \Ham(M,\Omega)$ — гладкая петля гамильтоновых
диффеоморфизмов, определяющая строго эргодическое косое произведение
$T(y, t) = (h(t)y, t+\alpha)$. 
Пусть $\gamma$ — соответствующий элемент из $\pi_1(\Ham(M, \Omega))$.
Обозначим через $H(x, t)$ нормализованный гамильтониан, порождающий
петлю $h(t)^{-1}$. 
Пусть $h_k(t) = h(t + k\alpha)^{-1}$ и  
\[f_N(t) = h_0(t) \circ \dots \circ h_{N-1}(t).\]
Из \ref{1.4.D} следует, что петля $f_N$ порождается нормализованным гамильтонианом 
\begin{align*}
F_N(y,t)
&=
H(y,t)
+ H(h_0(t)^{-1}y,t+\alpha)
+\dots
\\
&\dots
+
H(h_{N-2} (t)^{-1} \circ \dots\circ h_0(t)^{-1}y, t + (N - 1)\alpha).
\end{align*}
Это выражение можно переписать следующим образом:
\[F_N(y, t) = \sum_{k=0}^{N-1} H\circ T^k(y, t).\]
Поскольку гомеоморфизм $T$ строго эргодичен и функция $F_N$ имеет
нулевое среднее, мы заключаем, что
\[\frac1N\int_0^1\big(\max_{y\in M}F_N(y, t) - \min F_N(y, t)\big)\,dt \to 0\]
при $N \to \infty$.

Но выражение в левой части в точности равно $\frac1N\length\{f_N(t)\}$.
Обратите внимание, что петля $\{f_N(t)\}$ представляет элемент $-N\gamma$.
Поскольку $\nu(N\gamma) = \nu(-N\gamma)$, мы получаем, что
$\nu(N\gamma)/N$ стремится к нулю при $N \to \infty$, 
то есть асимптотическая норма элемента $\gamma$ обращается в нуль.
\qeds

Я не знаю \emph{точного} значения $\nu_\infty (\gamma)$ ни в одном примере,
где эта величина строго положительна (например, для раздутия пространства $\CP^2$ в
\ref{sec:11.1}).
Трудность в том, что во всех известных примерах, где можно точно
вычислить хоферовскую норму $\nu(\gamma)$, существует замкнутая петля
$h(t)$, минимизирующая длину в своём гомотопическом классе (то есть
\emph{замкнутая кратчайшая}).
Однако, оказывается, что всякая непостоянная замкнутая кратчайшая
перестаёт быть таковой после достаточного числа итераций.
Иными словами, петлю $h(Nt)$ можно сократить, если $N$ достаточно велико.
Доказательство основано на следующем обобщении
описанной выше процедуры сокращения.
Пусть $H(y, t)$ — нормализованный гамильтониан петли $h(t)^{-1}$.
Не умаляя общности, можно предположить, что $h(0) = \1$, и что $H(y,0)$
не обращается в нуль тождественно.
Обозначим через $\Gamma$ множество всех точек многообразия $M$, в которых функция
$|H(y, 0)|$ достигает максимального значения.
Поскольку $M \setminus \Gamma$ — непустое открытое подмножество, а
группа гамильтоновых диффеоморфизмов действует транзитивно на $M$,
можно выбрать такую последовательность
\[\1=\phi_0,\phi_1,\dots,\phi_{N-1}\in\Ham(M,\Omega),\]
что 
\[\Gamma\cap\phi_1(\Gamma)\cap\dots\cap\phi_{N-1}(\Gamma)=\emptyset\]
Рассмотрим петлю $f_N(t) = h(t)^{-1} \circ \phi_1h(t)^{-1}\phi_1^{-1}
\circ \dots \circ \phi_{N-1}h(t)^{-1}\phi_{N-1}$.
Покажем, что она короче, чем петля $h(Nt)$.
Действительно, заметим, что её гамильтониан $F_N$ в момент времени $t
= 0$ записывается следующим образом:
\[F_N(y,0) = \sum_{i=0}^{N-1} H(\phi_i^{-1}y, 0).\]
Пусть
\[a(t)
=
\max_{y\in M} F_N(y, t) - \min_{y\in M}  F_N(y, t)
\]
и
\[b(t)
=
N(\max_{y\in M} H(y, t) - \min_{y\in M}  H(y, t)).
\]
Выбор последовательности $\{\phi_i\}$ влечёт, что $a(0) < b(0)$.
Поскольку $a(t) \z\le b(t)$ для всех $t$, мы получаем, что
$\int_0^1a(t)\,dt < \int_0^1 b(t)\,dt$, и это доказывает утверждение.
Мы заключаем, что если \textit{ненулевой класс $\gamma\in
  \pi_1(\Ham(M, \Omega))$ представлен минимальной замкнутой
  геодезической, то $\nu_\infty (\gamma)$ строго меньше, чем
  $\nu(\gamma)$.}

Было бы интересно исследовать дальнейшие ограничения на гомотопические
классы гладких строго эргодических петель в группе гамильтоновых
диффеоморфизмов. 

\section{Алгебра в помощь}

Вернёмся к теореме \ref{11.1.A}.
В этом разделе мы наметим доказательство того, что нестягиваемые петли
в $\Ham(M_\ast, \Omega_\ast)$ не бывают строго эргодичными.

Пусть $(M^{2n},\Omega)$ — замкнутое симплектическое многообразие.
Рассмотрим отображение
\[I\:\pi_1(\Ham(M,\Omega))\to\RR\]
определённое следующим образом.
Пусть $\gamma\in\pi_1(\Ham(M,\Omega))$.
Рассмотрим соответственное симплектическое расслоение $P(\gamma)$.
Обозначим через $u$ первый класс Черна вертикального касательного
расслоения над $P(\gamma)$;
его слой в точке $x\in P(\gamma)$ есть (симплектическое) векторное пространство, касательное к слою, проходящему через $x$. 
Как и раньше, $c$ обозначает класс сцепления.
Определим «характеристическое число» 
\[I(\gamma)= \int_{P(\gamma)} c^n\smallsmile u.\]
Легко видеть, что $I\: \pi_1(\Ham(M, \Omega))\to \RR$ — гомоморфизм
(\rindex{Лалонд}\cite{P6,LMP2}). 

\begin{ex}{Определение}\label{11.3.A} Симплектическое многообразие $(M,\Omega)$ называется \rindex{монотонное многообразие}\emph{монотонным}, если $[\Omega]$ является положительным кратным класса $c_1 (\T M)$.
\end{ex}

\begin{thm}[(\cite{P6})]{Теорема}\label{11.3.B}
Пусть $(M, \Omega)$ — замкнутое монотонное симплектическое многообразие.
Тогда найдётся такая положительная константа $C>0$, что
$\nu(\gamma)\ge C|I(\gamma)|$ для всех $\gamma\z\in\pi_1(\Ham(M,\Omega))$.
\end{thm}

Иными словами, гомоморфизм $I$ калибрует хоферовскую норму на
фундаментальной группе.
Доказательство теоремы основано на теории, описанной в двух
предыдущих главах, в сочетании с результатами из \cite{Se}.
Недавно \rindex{Зейдель}Зейдель получил обобщение этого неравенства на немонотонные
симплектические многообразия.
Так как $I$ — гомоморфизм, оценка в \ref{11.3.B} проходит для
асимптотической хоферовской нормы:
\[\nu_\infty(\gamma)\ge C|I(\gamma)|\]
В отличии от нормы $\nu$, гомоморфизм $I$ можно относительно легко вычислить, и в этом его большое преимущество.
Например, можно показать, что $I(\gamma) \ne 0$, где $\gamma$ —
образующая $\pi_1(\Ham(M_\ast,\Omega_\ast)) = \ZZ$.
Таким образом, из \ref{11.3.B} следует, что асимптотическая норма
каждого нетривиального элемента $\pi_1(\Ham(M_\ast, \Omega_\ast))$
строго положительна.
Из \ref{11.2.A} следует, что такой элемент не представим
строго эргодической петлёй.

 \chapter[Геодезические]{Элементы вариационной теории геодезических}\label{chap:12}

Мы уже обсудили ряд результатов о геодезических группы гамильтоновых
диффеоморфизмов. 
В этой главе мы посмотрим на геодезические с точки зрения
вариационного исчисления. 
Пусть задан гладкий путь гамильтоновых диффеоморфизмов — можно ли
его сократить малой вариацией с фиксированными концами? 
Этот вопрос приходит из классической теории геодезических на римановых многообразиях. 
Интересно, что в хоферовской геометрии, по крайней мере при некоторых
предположениях о невырожденности, у этого вопроса есть довольно точный
ответ с ясным динамическим смыслом \rindex{Устиловский}\cite{U}. 

\section{Что такое геодезическая?}

Риманова интуиция подсказывает, что геодезические надо определять
как критические точки функционала длины. 
Попробуем формализовать это определение.

Пусть $\{f_t\}$, $t\in[a; b]$, — гладкий путь в $\Ham(M,\Omega)$.
\rindex{вариация}\emph{Вариацией} $\{f_t\}$ называется гладкое
семейство путей $\{f_{t,\epsilon}\}$ с $t \in [a; b]$ и $\epsilon \in
(-\epsilon_0; \epsilon_0)$, удовлетворяющих условиям
\[f_{a,\epsilon} = f_a,\quad f_{b,\epsilon} = f_b\quad\text{и}\quad f_{t,0} = f_t\]
для всех $t$ и $\epsilon$.
Мы всегда предполагаем, что общий носитель
$\overline{\bigcup_{t,\epsilon}\supp f_{t,\epsilon}}$ компактен.
При такой вариации рассмотрим длину пути $\{f_{t,\epsilon}\}$ как функцию от $\epsilon$ \index[symb]{$\ell(\epsilon)$}
\begin{align*}
\ell(\epsilon)&=\int_a^b\|F(\cdot,t,\epsilon)\|\,dt=
\\
&=\int_a^b \big(\max_x F(x,t,\epsilon)-\min_x F(x,t,\epsilon)\big)\,dt,
\end{align*} 
где $F(x, t, \epsilon)$ — гамильтониан, порождающий путь $\{f_{t,\epsilon}\}$ для данного~$\epsilon$.

\begin{ex}{Предварительное определение}\label{12.1.A}
Путь $\{f_t\}$ является геодезической, если
\begin{itemize}
\item его скорость постоянна, то есть $\|F(\cdot, t)\|$ не зависит от $t$;
\item для любой гладкой вариации $\{f_t\}$ функция длины $\ell(s)$ имеет критическую точку при $\epsilon = 0$.
\end{itemize}
\end{ex}

Тут мы сталкиваемся с трудностью — даже при гладких вариациях функция $\ell(\epsilon)$ не обязана быть гладкой!
Таким образом, необходимо уточнить понятие критической точки.
Наша первая задача — выяснить структуру функций длины
$\ell(\epsilon)$, связанных с вариациями данного пути.

\begin{ex}{Предложение}\label{12.1.B}
Для каждой вариации её функция длины $\ell(\epsilon)$ выпукла в $0$ с точностью до
второго порядка. 
То есть существуют выпуклая функция $u(\epsilon)$, а также числа
$\delta > 0$ и $C > 0$ такие, что  
\[|\ell(\epsilon) - u(\epsilon)| \le C\epsilon^2\]
для всех $\epsilon \in (-\delta; \delta)$.
\end{ex}

Обратите внимание, что, каким бы ни было определение критической точки,
оно не должно зависеть от членов второго порядка. 
Далее, единственным естественным кандидатом на критическую точку выпуклой
функции является её точка минимума.
Таким образом, мы приходим к следующему понятию.

\begin{ex}{Определение}\label{12.1.C}
Значение $\epsilon = 0$ является критической точкой функции длины $\ell(\epsilon)$ тогда и только тогда, когда она является точкой минимума выпуклой функции $u(\epsilon)$, удовлетворяющей неравенству \ref{12.1.B}
\end{ex}

\begin{ex}{Упражнение}\label{12.1.D}
Проверьте следующие утверждения.
\begin{itemize}

\item Убедитесь, что приведённое выше определение правильное, то есть не зависит от выбора выпуклой функции $u$, удовлетворяющей условию \ref{12.1.B}.

\item Докажите, что если функция $\ell(\epsilon)$ гладкая, то определение \ref{12.1.C} совпадает с обычным.

\item Докажите, что если $\ell$ достигает своего локального минимума или максимума в
  нуле, то ноль является критической точкой $\ell$ в смысле
  \ref{12.1.C}.

\item Выведите из \ref{12.1.B}, что росток функции $1 - |\epsilon|$
  в точке $0$ не может возникнуть как функция длины какой-либо гладкой
  вариации.
\end{itemize}
\end{ex}

Теперь мы готовы ответить на вопрос, поставленный в заголовке раздела:
\rindex{геодезическая}\emph{путь $\{f_t\}$, $t\in[a; b]$, называется
  геодезической, если он удовлетворяет \ref{12.1.A} и \ref{12.1.C}.} 

Как мы увидим в \ref{12.2.A}, ограничение геодезической, определённой
на $[a; b]$, на любой подотрезок снова является геодезической. 
Таким образом, понятие геодезических естественно распространяется на пути, определённые на произвольных интервалах времени. 
А именно, пусть $I\subset \RR$ — интервал и $f\: I\to\Ham(M,\Omega)$ — гладкий путь.
Если ограничение $f$ на любой отрезок $[a; b] \subset I$ —
геодезическая, то и сам путь $f$ называется геодезической. 
Далее мы сосредоточимся на геодезических $\{f_t\}$, которые определены
на единичном интервале $[0; 1]$. 
Более того, поскольку хоферовская метрика биинвариантна, мы всегда
можем сдвинуть геодезическую и считать, что $f_0 = \1$.  

Перейдём к доказательству \ref{12.1.B}.
Нам будет удобно использовать идею линеаризации, представленную в главе 5.
Рассмотрим пространство \index[symb]{$\F_0$}$\F_0$ всех
нормализованных гамильтонианов $M \times [0; 1] \to \RR$ с нормой
\index[symb]{$\VERT F \VERT_0$} 
\[\VERT F \VERT_0 = \int_0^1 \big(\max_x F(x,t) - \min_x F(x,t)\big)\, dt.\]
Обозначим через \index[symb]{$\H_0$}$\H_0 \subset \F_0$ подмножество,
состоящее из всех гамильтонианов, порождающих петли $\{h_t\}$
гамильтоновых диффеоморфизмов: $h_0 = h_1 = \1$.
В отличие от соглашений в главе 5,
мы \textbf{не предполагаем}, что гамильтонианы в $\F_0$ и в $\H_0$
периодичны по времени. 
Обозначим через \index[symb]{$\V$}$\V$ множество всех гладких семейств
$H(x, t, \epsilon)$ функций из $\H_0$ таких, что $H(x, t, 0) \equiv
0$. 

\begin{thm}{Предложение}\label{12.1.E}
Пусть $\{f_t\}$ — гладкий путь гамильтоновых диффеоморфизмов с $t
\in [0; 1]$ и $f_0 = \1$. 
Множество функций длины $\ell(\epsilon)$, связанных с вариациями пути $\{f_t\}$, состоит из всех функций вида 
\[\epsilon\mapsto\VERT F - H(\epsilon)\VERT_0,\]
где $H\in\V$.
\end{thm}

Предложение \ref{12.1.B} является непосредственным следствием этой формулы.
Действительно, положим $u(\epsilon) = \VERT F - \epsilon H'(0)\VERT_0$.
Ясно, что $u$ выпукло и совпадает с $\ell$ с точностью до второго порядка.

\parit{Доказательство \ref{12.1.E}.}
Любая вариация пути $f_t$ записывается в виде $f_{t,\epsilon} = h_{t,\epsilon}^{-1}\circ f_t$, где $h_{1,\epsilon}$ — гладкое семейство петель с $h_{t,0} = \1$.
Гамильтониан $H(x,t,\epsilon)$ петель $\{h_{t,\epsilon}\}$ принадлежит $\V$.
Далее, 
\[F(x, t, \epsilon) = -H(h_{t,\epsilon}x, t, \epsilon) + F(h_{t,\epsilon}x, t),\]
поэтому
\[\ell(\epsilon) = \VERT F - H(\epsilon)\VERT_0.\]
\qeds

Решите следующие упражнения, используя \ref{12.1.E}.

\begin{ex}{Упражнение}\label{12.1.F}
Предположим, что при всяком $t$ функция $F(x, t)$ имеет единственную
точку максимума и минимума соответственно, причём эти точки
невырождены в смысле теории Морса. 
Тогда \textit{для любой} вариации $\{f_t\}$ функция длины $\ell(\epsilon)$ гладкая в окрестности нуля. 
\emph{Подсказка:} используйте теорему о неявной функции.
\end{ex}

\begin{ex}{Упражнение}\label{12.1.G}
Предположим, что множество точек максимума $F$ имеет непустую внутренность.
Постройте вариацию, у которой функция длины негладкая в нуле.%
\footnote{Достаточно чтобы максимум достигался хотя бы в двух точках. — \textit{Прим. ред.}}
\end{ex}

\section{Описание геодезических}\label{sec:12.2}

Будем говорить, что гамильтониан $F\in \F_0$ имеет \rindex{фиксированные экстремумы}\emph{фиксированные экстремумы}, если существуют такие две точки $x_-, x_+\in M$, что $F(x_-, t) \z= \min_x F(x, t)$ и $F(x_+, t) = \max_x F (x, t)$ для всех $t\in[0; 1]$, причём функция $F(x_+, t) - F(x_-, t)$ не зависит от $t$.
Важнее всего то, что экстремальные точки $x_-$ и $x_+$ функции $F(\cdot, t)$ не зависят от времени.

\begin{thm}{Теорема}\label{12.2.A}
Путь $\{f_t\}$ является геодезическим тогда и только тогда, когда соответствующий гамильтониан $F\in\F_0$ имеет фиксированные экстремумы.
\end{thm}

В частности, каждый автономный гамильтонов поток является геодезическим.
Гамильтонианы с фиксированными экстремумами введены в \cite{BP1}, где они были названы \rindex{квазиавтономный поток}\emph{квазиавтономными}.
Теорема \ref{12.2.A} по существу доказана в \rindex{Лалонд}\cite{LM2} (хотя там геодезическая определяется немного иначе).

Для доказательства теоремы \ref{12.2.A} нам предстоит более детально исследовать структуру вариаций.
Обозначим через $\V_1$ касательное пространство к $\H_0$ в точке $0$: 
\[\V_1=\set{\frac{\partial H}{\partial\epsilon}(0)
}{H\in \V}\]

\begin{thm}{Предложение}\label{12.2.B}
Пространство $\V_1$ состоит из всех функций $G\in\F_0$, удовлетворяющих условию
\[\int_0^1G(x,t)\,dt=0\]
при всех $x\in M$.
\end{thm}

Предложение доказывается так же как \ref{5.2.D} и \ref{6.1.C}, где решён случай периодических во времени вариаций.

\parit{Доказательство теоремы \ref{12.2.A}.}
Пусть $F\in\F_0$ — такая функция, что $\|F(\cdot, t)\|$ не зависит от $t$.
Для семейства $H(\epsilon)$ из $\V$ положим $G = H'(0)\in\V_1$.
Определим функции $\ell(\epsilon) = \VERT F - H(\epsilon)\VERT_0$ и $v(\epsilon) = \VERT F - \epsilon G\VERT_0$.
С учётом \ref{12.1.E} гамильтониан $F$ порождает геодезическую тогда и только тогда, когда $\ell(\epsilon)$ имеет критическую точку в смысле \ref{12.1.C} при $\epsilon= 0$ для каждого $H\in\V$.
Поскольку $\ell(\epsilon)$ и $v(\epsilon)$ совпадают с точностью до членов второго порядка, а $v(\epsilon)$ выпукла, это эквивалентно тому, что $v(\epsilon)$ имеет точку минимума при $\epsilon= 0$ для любого $G\in\V_1$.
Поэтому для доказательства теоремы \ref{12.2.A} достаточно проверить эквивалентность следующих условий:
\begin{enumerate}[(i)]
\item\label{12.2.i} $F$ имеет фиксированные экстремумы;
\item\label{12.2.ii} $\VERT F - \epsilon G\VERT_0 \ge \VERT F\VERT_0$ для всех $G\in\V_1$, $\epsilon\in\RR$.
\end{enumerate} 

Предположим, что выполняется (\ref{12.2.i}).
Положим 
$u=F-\epsilon G$ для $G\in \V_1$.
Тогда из \ref{12.2.B} следует, что 
\[\VERT u\VERT_0 \ge \int_0^1 u(x_+,t)-u(x_-,t)\,dt=\VERT F\VERT_0,\]
и, следовательно, мы получаем (\ref{12.2.ii}).

А теперь предположим, что выполняется (\ref{12.2.ii}).
Возьмём 
\[G(x, t) = F(x, t) - \int_0^1F(x,t)\,dt.\]
Неравенство (\ref{12.2.ii}) даёт 
\[
\max_x\!\int_0^1 F(x, t)\,dt- \min_x\!\int_0^1 F(x, t)\,dt
\ge
\int_0^1 \big(\max_xF(x, t)\,dt-\min_x F(x, t)\big)\,dt,
\]
что возможно, только если $F$ имеет фиксированные экстремумы.
\qeds

\section{Устойчивость и сопряжённые точки}

Геодезическая называется \rindex{устойчивая
  геодезическая}\emph{устойчивой}, если при каждой вариации функция
длины $\ell(\epsilon)$ достигает своего минимального значения в нуле. 
Иными словами, устойчивые геодезические нельзя укоротить небольшими
вариациями с фиксированными концами. 
Задача описания устойчивых геодезических в полной общности остаётся открытой.%
\footnote{Думаю, её можно решить существующими методами негладкого анализа.
Заметим также, что первое утверждение теоремы~\ref{12.3.A},
которое даёт достаточное условие устойчивости, выполняется для
{}\emph{любых} геодезических, см. \cite{LM3}.  
Доказательство использует теорию псевдоголоморфных кривых.}
Ниже мы приведём решение для определённого класса
\rindex{невырожденная геодезическая}\emph{невырожденных}
геодезических. 
По определению геодезическая невырождена, если для каждого $t$
соответствующий гамильтониан имеет единственную точку максимума и
минимума соответственно, причём эти точки невырождены в смысле теории
Морса. 
Например, каждый автономный путь, порождённый гамильтонианом с единственными невырожденными максимумом и минимумом, является невырожденной геодезической. 
Напомним, что для любой вариации невырожденной геодезической функция длины $\ell(\epsilon)$ гладкая (см. \ref{12.1.F}). 
Теория, развитая в этом и двух следующих разделах, восходит к
\rindex{Устиловский}Устиловскому \cite{U} (см. также
\rindex{Лалонд}\cite{LM2}). 

Пусть $\{f_t\}$, $t \in [0; 1]$, $f_0 = \1$ — невырожденная
геодезическая, порождённая гамильтонианом $F\in \F_0$. 
Обозначим через $x_-$ и $x_+$ соответственно независимые от времени
точки минимума и максимума функции $F(\cdot, t)$. 
Заметим, что $x_\pm$ в этом случае являются неподвижными точками $\{f_t\}$.
Рассмотрим линеаризованные потоки $f_{t*}$ на $\T_{x_+}M$ и $\T_{x_-}M$.
Мы говорим, что такой поток имеет нетривиальную $T$-периодическую
орбиту с $T\ne0$, если $f_{T*}\xi=\xi$ для некоторого касательного
вектора $\xi\ne0$. 

\begin{thm}[(\cite{U})]{Теорема}\label{12.3.A}

\begin{itemize}
\item Предположим, что линеаризованные потоки не имеют нетривиальных $T$-периодических орбит с $T\in(0;1]$.
Тогда $\{f_t\}$ устойчив.
\item Предположим, что $\{f_t\}$ устойчив.
Тогда линеаризованные потоки не имеют нетривиальных периодических орбит с $T\z\in(0;1)$.
\end{itemize}
\end{thm}

Этот результат можно интерпретировать как описание сопряжённых точек вдоль геодезических хоферовской метрики —  сопряжённые точки соответствуют нетривиальным замкнутым орбитам линеаризованного потока в точках $x_\pm$.
До сопряжённой точки геодезическая устойчива, а после неё теряет устойчивость (то есть геодезическую можно укоротить малой вариацией).
За дополнительными сведениями мы отсылаем читателя к \cite{U}.
Доказательство \ref{12.3.A} дано в \ref{sec:12.5}.
Оно основано на формуле второй вариации, которую мы опишем в следующем разделе.

\begin{ex}{Упражнение}\label{12.3.B}
Выведите из \ref{12.3.A}, что \textit{достаточно короткий} отрезок невырожденной геодезической устойчив.
\end{ex}

\section{Формула второй вариации}
\rindex{формула второй вариации}

Пусть $\{f_t\}$ — невырожденная геодезическая, порождённая
гамильтонианом $F(x, t)$ с точками максимума/минимума $x_\pm$. 
Обозначим через $C_\pm(t)$ оператор линеаризованного уравнения в точке $x_\pm$,
то есть 
\[\frac{d}{dt} f_{t*}(x_\pm)=C_\pm(t)f_{t*}(x_\pm).\]
Рассмотрим пространства 
\[V_\pm = \set{\text{гладкие отображения}\ v\: [0;1] \to
  \T_{x_\pm}M}{v(0) = v(1) = 0}.\] 
Для данной вариации $\{f_{t,\epsilon}\}$ геодезической $\{f_t\}$
определим элементы $v_\pm\in V_\pm$ формулой  
\[v_\pm = \left.\frac d{d\epsilon}\right|_{\epsilon=0} f_{t,\epsilon} x_\pm.\]
Будет удобно рассматривать максимальную и минимальную части функции
длины, связанные с вариацией, по отдельности. 
Положим 
\[\ell_+(\epsilon) =\int_0^1\max_x F(x,t,\epsilon)\,dt\]
и 
\[\ell_{-}(\epsilon)=\int_0^1\min_x F(x,t,\epsilon)\,dt.\]
Ясно, что $\ell(\epsilon) = \ell_+(\epsilon) - \ell_{-}(\epsilon)$.

\begin{thm}[\cite{U}]{Теорема}\label{12.4.A}\rindex{Устиловский}
\[\frac{\d^2\ell_\pm}{d\epsilon^2}(0)=Q_\pm(v_\pm)\]
и, следовательно,
\[\frac{\d^2\ell_\pm}{d\epsilon^2}(0)=Q_+(v_+)-Q_-(v_-),\]
где
\[Q_\pm(v)=-\int_0^1\left(\Omega(C^{-1}_\pm \dot v,\dot v)+\Omega(\dot
v,\dot v)\right)\,dt.\] 

\end{thm}

\begin{ex*}[(изопериметрическое неравенство)]{Пример}
Рассмотрим стандартную симплектическую плоскость $\RR^2(p, q)$ с
симплектической формой $\omega = dp \wedge dq$. 
Пусть $v\: [0;1] \to\RR^2$ — гладкая кривая такая, что $v(0)=v(1)=0$.
Напомним следующие понятия евклидовой геометрии: 
\begin{align*}
\length(v)&=\int_0^1|\dot v|\,\d t,
\\
\energy(v)&=\int_0^1|\dot v|^2\,\d t,
\\
\area(v)&=\int_0^1\langle \tfrac12(p\d q-q\d p),\dot v\rangle\,\d t=
\\
&=\frac12\int_0^1 (p\dot q-q\dot p)\,\d t=
\\
&=\frac12\int_0^1\omega(v,\dot v)\,\d t,
\end{align*}
где $v(t) = (p(t), q(t))$.
Изопериметрическое неравенство гласит, что
\[4\pi\area(v)\le\length(v)^2\le \energy(v).\]
Вернёмся к нашей симплектической задаче.
Мы предполагаем, что
гамильтониан вблизи $x_-$ задаётся формулой $F(x) = \pi \lambda(p^2
\z+ q^2)$ с $\lambda > 0$. 
Таким образом, мы получаем гамильтонову систему 
\[
\begin{cases}
\quad\dot p &= -2\pi\lambda q\;,
\\
\quad\dot q &= 2\pi\lambda p\;.
\end{cases}
\]
Полагая $z = p + iq$, мы получаем линейное уравнение $\dot z = 2\pi
\lambda iz$. 
В этом случае $C_-(t) = 2\pi\lambda i$ и $C_-^{-1}(t) = -\frac i{2\pi\lambda}$.
Учитывая, что $\omega(\xi,i\xi) \z= |\xi|^2$, получаем
\begin{align*}
Q_-(v)&=-\int_0^1\left(\omega(-\tfrac i{2\pi\lambda}\dot v,\dot
v)+\omega(\dot v, v)\right)\,dt= 
\\
&=-\frac 1{2\pi\lambda}\int_0^1|\dot v|^2\,dt-\int_0^1\omega(\dot v,v)\,dt=
\\
&=-\frac 1{2\pi\lambda}\energy(v)+2\area(v)\;.
\end{align*}
Уравнение $\dot z = 2\pi\lambda iz$ имеет решение $z(t) =
e^{2\pi\lambda it}z(0)$ при $t \in [0;1]$ и, следовательно, при
$\lambda < 1$ у него нет орбит периода $1$. 
Применив теоремы \ref{12.3.A} и \ref{12.4.A}, получаем, что $Q_-(v)\le
0$ для всех плоских кривых $v$ удовлетворяющих $v(0) = v(1) = 0$. 
Поэтому 
\[4\pi \lambda \area(v) \le \energy (v)\]
при всех $\lambda < 1$.
При $\lambda\to1$ получаем изопериметрическое неравенство.
\end{ex*}

В доказательстве \ref{12.4.A} нам понадобится следующее вспомогательное утверждение.

\begin{thm}{Лемма}\label{12.4.B}
Пусть $\{f_{t,\epsilon}\}$ — вариация пути $\{f_t\}$, а $F(x, t,\epsilon)$ — её гамильтониан.
Тогда 
\[\int_0^1\frac{\partial F}{\partial \epsilon}(f_{t,\epsilon}x,t,\epsilon)\,\d t=0.\] 

\end{thm}

\parit{Доказательство.}
\begin{enumerate}[1)]
\item Приведённая выше формула верна для любой вариации постоянной петли, то есть когда $f_t = \1$ при всех $t$.
Доказательство то же что в \ref{6.1.C}.
\item Рассмотрим теперь общий случай.
Запишим $f_{t,\epsilon} = f_t \circ h_{1,\epsilon}$, где $h_{1,\epsilon}$ — вариация постоянной петли.
Обозначим через $H(x, t, \epsilon)$ гамильтониан петли $\{h_{1,\epsilon}\}$.
Тогда $F(x,t,\epsilon) = F(x,t) + H(f_t^{-1}x,t,\epsilon)$, и, значит, 
\[\frac{\partial F}{\partial \epsilon}(x,t,\epsilon)=\frac{\partial H}{\partial \epsilon}(f_t^{-1}x,t,\epsilon).\]
Таким образом,
\[\frac{\partial F}{\partial \epsilon}(f_{t,\epsilon}x,t,\epsilon)=\frac{\partial H}{\partial \epsilon}(h_{t,\epsilon}x,t,\epsilon).\]
Требуемое утверждение следует теперь из шага 1.
\qeds
\end{enumerate}

\parit{Доказательство \ref{12.4.A}.}
Давайте вычислим вторую производную $\ell_+(\epsilon)$,
а с $\ell_-(\epsilon)$ поступим точно также.
Будем работать в стандартных симплектических координатах $x = (p, q)$ вблизи $x_+$ и $\langle\xi,\eta\rangle$ обозначает евклидово скалярное произведение.
Обозначим через $i$ комплексную структуру $(p, q) \mapsto (-q, p)$.
В этих обозначениях $\Omega(\xi, i\eta) = \langle\xi,\eta\rangle$ и гамильтониан запишется как
\[\frac{\d}{\d t}f_tx=i\frac{\partial F}{\partial x}(f_tx,t).\]
Значит
\[C_+(t)=i\frac{\partial^2F}{\partial x^2}(x_+,t).\]
Для упрощения формул введём следующие обозначения: 
\begin{align*}
a&=\frac{\partial^2 F}{\partial x\partial\epsilon}(x_+,t,0),
\\
b&=\frac{\partial x_+}{\partial\epsilon}(t,0),
\\
c&=\frac{\partial^2 F}{\partial\epsilon^2}(x_+,t,0),
\\
K&=\frac{\partial^2 F}{\partial x^2}(x_+,t).
\end{align*}

Теорема о неявной функции гарантирует, что $F(\cdot, t, \epsilon)$ имеет единственную точку максимума $x_+(t, \epsilon)$, которая гладко зависит от $t$ и $\epsilon$.
Поскольку 
\[\ell_+(\epsilon)=\int_0^1F(x_+(t,\epsilon),t,\epsilon)\,\d t,\]
получаем, что
\[\frac{\d \ell_+}{\d \epsilon}(\epsilon)
=
\int_0^1
\frac{\partial F}{\partial x}(x_+(t,\epsilon),t,\epsilon)\frac{\partial x_+}{\partial \epsilon}(t,\epsilon)
+
\frac{\partial F}{\partial x}(x_+(t,\epsilon),t,\epsilon)\,\d t.
\]
Поскольку $x_+(t,\epsilon)$ является критической точкой $F(\cdot,t,\epsilon)$, это выражение можно упростить
\[\frac{\d \ell_+}{\d \epsilon}(\epsilon)
=
\int_0^1
\frac{\partial F}{\partial x}(x_+(t,\epsilon),t,\epsilon)\,\d t.
\]
Снова дифференцируя по $\epsilon$ и полагая $\epsilon = 0$, мы получаем 
\begin{equation}
\frac{\d^2 \ell_+}{\d\epsilon^2}(0)
=
\int_0^1 (\langle a, b\rangle +c)\,\d t.
\label{eq:12.4.C}
\end{equation}
Дифференцируя уравнение 
$\frac{\partial F}{\partial x}(x_+ (t, \epsilon), t, \epsilon) = 0$ по
$\epsilon$, получаем 
\[Kb+a = 0.\]
Условие невырожденности гарантирует, что коэффициент $K$ обратим и,
таким образом,  
\begin{equation}
b=-K^{-1}a.
\label{eq:12.4.D}
\end{equation}

Согласно лемме~\ref{12.4.B}, 
\[\int_0^1\frac{\partial F}{\partial
  \epsilon}(f_{t,\epsilon},t,\epsilon)\,\d t=0.\] 
Дифференцируя это равенство по $\epsilon$ при $\epsilon = 0$ и $x =
x_+$, получаем  
\begin{equation}
\int_0^1(\langle a, v_+\rangle+c)\,\d t = 0.
\label{eq:12.4.E}
\end{equation}

Рассмотрим уравнение Гамильтона
\[\frac{\d}{\d t}f_{t,\epsilon}x=i\frac{\partial F}{\partial
  x}(f_{t,\epsilon}x,t,\epsilon).\] 
Продифференцировав его по $\epsilon$ при $\epsilon = 0$ и $x = x_+$,
получим
\[\dot v_+=iKv_++ia,\]
или, эквивалентно,
\begin{equation}
a=-i\dot v_+-Kv_+.
\label{eq:12.4.F}
\end{equation}

Применив \ref{eq:12.4.D} и \ref{eq:12.4.F} вместе, получаем
\[b=-K^{-1}(-i\dot v_+ - K v_+) = K^{-1}i\dot v_+ + v_+.\]
С учётом \ref{eq:12.4.C} получим
\begin{align*}
\frac{\d^2\ell_+}{\d\epsilon}(0)
&=
\int_0^1(\langle a, b\rangle+c)\,\d t=
\\
&=\int_0^1(\langle a, K^{-1}i\dot v_+\rangle+\langle a, v_+\rangle+c)\,\d t=
\\
&\!\!\!\stackrel{\text{\tiny \ref{eq:12.4.E}}}{=}
\int_0^1 \langle a,
K^{-1}i\dot v_+\rangle\,\d t 
\\
&\!\!\!\stackrel{\text{\tiny \ref{eq:12.4.F}}}{=}
-\int_0^1\langle i\dot
v_++Kv_+, K^{-1}i\dot v_+\rangle\,\d t= 
\\
&\stackrel{{*}}{=}
-\int_0^1
(\langle K^{-1}i\dot v_+, i\dot v_+\rangle +\langle v_+,
i\dot v_+\rangle)\,\d t= 
\\
&=
-\int_0^1
(-\langle C_+^{-1}\dot v_+, i\dot v_+\rangle
+\langle v_+, i\dot v_+\rangle)\,\d t=
\\
&=
-\int_0^1
(\Omega(C_+^{-1}\dot v_+,\dot v_+)+\Omega(\dot v_+,v_+))\,\d t.
\end{align*}
Равенство $({*})$ следует из симметрии $K$.
Кроме того, мы воспользовались тем, что $C_+^{-1} = -K^{-1} i$ и
$\Omega(\xi, i\eta) = \langle \xi, \eta\rangle$ 
(или, эквивалентно, $-\Omega(\xi, \eta) = \langle \xi, i\eta\rangle$).
\qeds

\section{Анализ формулы второй вариации}\label{sec:12.5}

\begin{thm}{Предложение}\label{12.5.A}
Для произвольных $a \in V_+$ и $b \in V_-$ существует вариация
$\{f_{t,\epsilon}\}$ пути $\{f_t\}$ такая, что $v_+ = a$ и $v_- = b$. 
\end{thm}

\parit{Доказательство.}
Выберем гамильтониан $G\in\F_0$ такой, что
\[\int_0^1G(x, t)\,\d t = 0\]
для всех $x\in M$ (то есть $G\in\V_1$ в обозначениях \ref{sec:12.2}).
Мы уточним этот выбор позже.
Определим вариацию $\{h_{1,\epsilon}\}$ постоянной петли $h_{1,0} =
\1$ следующим образом: $h_{1,\epsilon}$ — поток гамильтониана
$\int_0^tG(x,s)\,\d s$ с временной координатой $\epsilon$. 
Рассмотрим вариацию $f_{t,\epsilon} = f_th_{1,\epsilon}$ пути $\{f_t\}$.
Тогда
\begin{align*}
v_+(t) &= \left.\frac{\d}{\d \epsilon}\right|_{\epsilon=0} f_{t,\epsilon} x_+ =
\\
&=\left.\frac{\d}{\d \epsilon}\right|_{\epsilon=0}f_th_{t,\epsilon} x_+=
\\
&=
f_{t*} \left.\frac{\d}{\d \epsilon}\right|_{\epsilon=0} h_{1,\epsilon} x_+.
\end{align*}
По построению $\{h_{1,\epsilon}\}$,
\[\left.\frac{\d}{\d \epsilon}\right|_{\epsilon=0}h_{t,\epsilon}x_+
=
\int_0^t\sgrad G(x_+,s)\,\d s.\]
Пусть 
\[c(t) = \frac{\d}{\d t}f^{-1}_{t*} a(t).\]
Выберем теперь $G \in \V_1$ так, чтобы равенство
\[G(x, t) = \Omega(x - x_+, c(t))\]
выполнялось в стандартных симплектических координатах вблизи $x_+$.
Поскольку $a(0) = a(1) = 0$ по определению $V_+$, это условие
согласовано с тем, что $G\in\V_1$. 
Явное выражение для $G$ даёт 
\[\sgrad G(x_+, t) = c(t) = \frac{\d}{\d t}f^{-1}_{t*} a(t).\]
Интегрируя это равенство и учитывая, что $a(0) = 0$, получаем
\[\int_0^t\sgrad G(x_+, s)\,\d s = f_{t*}^{-1}a(t).\]
Отсюда следует, что $v_+(t) = a(t)$.
То же рассуждение работает для $x_-$.
Поскольку $x_+$ и $x_-$ различны, построенные нами вариации не мешают
друг другу, 
что завершает доказательство.
\qeds

\parit{Доказательство теоремы \ref{12.3.A}.}
Приняв в расчёт \ref{12.4.A} и \ref{12.5.A}, приходим к выводу, что
путь $\{f_t\}$ устойчив тогда и только тогда, когда оба
функционала $Q_+$ и $-Q_-$ неотрицательны. 
Классические методы вариационного исчисления дадут нам точные условия,
когда такое происходит. 
Мы сосредоточимся на $Q_+$, но те же рассуждения применимы и к $Q_-$.

Рассмотрим лагранжиан $\L(\dot v, v) = -\Omega(C_+^{-1}\dot v, \dot v)
- \Omega(\dot v, v)$. 
Требуется исследовать критические точки следующего функционала на~$V_+$:
\[\int_0^1 \L(\dot v(t), v(t))\,\d t.\]

Прежде всего мы утверждаем, что $\L$ удовлетворяет условию Лежандра,
то есть $\tfrac{\partial^2\L}{\partial\dot v^2}$ положительно
определена. 
Действительно, 
\[\L( \dot v, v) = -\langle iC_+^{-1} \dot v,\dot v\rangle - \langle i\dot v, v\rangle\]
и, следовательно, $\tfrac{\partial^2\L}{\partial\dot v^2}= - 2i\cdot C_+^{-1}$.

Так как 
\[C_+=i\frac{\partial^2F}{\partial x^2}(x_+)
\quad\text{получаем, что}\quad
-iC_+^{-1}=i\left(\frac{\partial^2F}{\partial x^2}(x_+)\right)^{-1}i\]
и поэтому 
\[
\left\langle i\left(\frac{\partial^2F}{\partial x^2}(x_+)\right)^{-1}i\xi,\xi\right\rangle
=
-\left\langle \left(\frac{\partial^2F}{\partial x^2}(x_+)\right)^{-1}i\xi,i\xi\right\rangle
\]
Точка $x_+$ является невырожденным максимумом $F$.
Таким образом, матрица $\frac{\partial^2F}{\partial x^2}(x_+)$
отрицательно определена, как и $\left(\frac{\partial^2F}{\partial x^2}(x_+)\right)^{-1}$.
Значит квадратичная форма $\tfrac{\partial^2\L}{\partial\dot v^2}$ положительно определена и утверждение следует.

Так как лагранжиан $\L$ квадратичен по $v$, постоянный путь $v_0(t) \z\equiv 0$ экстремален.
Для квадратичных функционалов, удовлетворяющих условию Лежандра, классическое вариационное исчисление говорит нам следующее.
Если путь $v_0$ — точка минимума, то при всех $T \in [0;1)$ уравнение Эйлера — Лагранжа (совпадающее с уравнением Якоби) не имеет нетривиальных решений с $v(0) = v(T) = 0$.
И наоборот, если таких решений нет при $T\in[0;1]$, то $v_0$ — точка минимума.
Остаётся связать уравнение Эйлера — Лагранжа с периодическими орбитами линеаризованного потока $f_{t*}$ на $\T_{x_+}M$.
Уравнение Эйлера — Лагранжа имеет вид
\begin{align*}
\frac{\d}{\d t}\frac{\partial \L}{\partial \dot v}
&=\frac{\partial \L}{\partial v}
\\
&\Updownarrow
\\
\frac{\d}{\d t}(-2iC_+^{-1}\dot v+iv)
&=-i\dot v
\\
&\Updownarrow
\\
\frac{\d}{\d t}(-2iC_+^{-1}\dot v)
&=
\frac{\d}{\d t}(-2iv)
\end{align*}
Если $v$ является решением, то $C_+^{-1}\dot v = v + \const$.
Пусть $w = v + \const$, тогда 
\[
\begin{cases}
\quad C_+^{-1}\dot w=w,
\\
\quad w(0)=w(T).
\end{cases}
\]
Но это как раз и означает, что $w$ является $T$-периодической орбитой линеаризованного потока, что доказывает теорему \ref{12.3.A}.
\qeds

\section{Кратчайшие}
Пусть $I \subset \RR$ — интервал.
Рассмотрим геодезическую $\{f_t\}$, $t\in I$ с единичной скоростью, то есть
\[\|\tfrac{\d}{\d t}f_t\|=1\]
при всех $t\in I$.
Такая геодезическая называется \rindex{локально кратчайшая}\emph{локально кратчайшей} 
(см. \ref{sec:8.2}), если для каждого $t\in I$ существует окрестность $U$ точки $t$ в $I$ такая, что  
\[\rho(f_a,f_b)=|a-b|\]
при всех $a,b\in U$.

Каждая геодезическая на римановом многообразии локально кратчайшая.
По существу это следует из того, что экспоненциальное отображение по
связности Леви-Чивиты является диффеоморфизмом в окрестности нуля. 
Подобного утверждения в хоферовской геометрии нет.
Тем не менее различные примеры служат подверждением следующей гипотезы.

\begin{thm}{Гипотеза}\label{12.6.A}
Каждая геодезическая хоферовской метрики локально кратчайшая.
\end{thm}

Впервые это было доказано для стандартной симплектической структуры на
$\RR^{2n}$ в \rindex{Бялый}\cite{BP1}. 
Позже эта гипотеза была подтверждена для некоторых других
симплектических многообразий \rindex{Лалонд}\cite{LM2}, включая
кокасательные расслоения, замкнутые ориентированные поверхности и
$\CP^2$. 

Полезно ослабить понятие локально кратчайшей следующим образом.
Будем говорить, что геодезическая $\{f_t\}$, $t\in I$, \rindex{локально
  слабая кратчайшая}\emph{локально слабая кратчайшая}, если для
каждого $t\in I$ существует такая окрестность $U$ точки $t$ в $I$, что
для любых $a, b\in U$ с $a < b$ геодезический отрезок $\{f_t\}$, $t\in
[a;b]$, минимизирует длину в гомотопическом классе путей с
фиксированными концами. 

\begin{thm}{Гипотеза}\label{12.6.B}
Каждая геодезическая хоферовской метрики является локально слабой кратчайшая.
\end{thm}

В конечномерной римановой геометрии легко доказать, что любая локально
слабая кратчайшая является локально кратчайшей. 
В хоферовской геометрии это верно при следующем дополнительном предположении.
Обозначим $S \subset [0; +\infty)$ спектр длин $\Ham(M, \Omega)$ (см. \ref{sec:7.3}). 
Положим $a = \inf S\setminus\{0\}$.
Как обычно, предполагается, что точная нижняя грань пустого множества
равена $+\infty$. 

\begin{thm}[(\cite{LM2})]{Предложение}\label{12.6.C}
Предположим, что $a>0$.
Тогда любая геодезическая является локально кратчайшей.
\end{thm}

\parit{Доказательство.}
Пусть $\{f_t\}$ — локально слабая кратчайшая.
Не умаляя общности, можно предположить, что её гамильтониан имеет единичную норму в каждый момент времени $t$.
Таким образом, существует $\epsilon \in (0; a/2)$ со следующим свойством.
Неравенство $\length\alpha \ge \epsilon$ выполняется для всякого пути $\alpha$ в $\Ham(M, \Omega)$, соединяющего $\1$ с $f_\epsilon$ и гомотопного $\{f_t\}$, $t\in[0; \epsilon]$.
Достаточно доказать, что $\rho(\1, f_\epsilon) = \epsilon$.
Рассуждая от противного, предположим, что $\length \beta < \epsilon$ для некоторого пути $\beta$, соединяющего $\1$ с $f_\epsilon$.
Обозначим через $\gamma$ путь $\{f_t\}$, $t\in[0; \epsilon]$.
Рассмотрим петлю, образованную $\gamma$ и $\beta$.
Её длина строго меньше $a$.
Следовательно, по определению $a$ эту петлю можно стянуть в сколь угодно короткую.
Но это означает, что путь $\gamma$ гомотопен пути $\gamma'$, длина которого сколь угодно близка к длине $\beta$.
Приходим к противоречию с неравенствами $\length \gamma' \ge\epsilon$ и $\epsilon > \length \beta$.
\qeds

В частности, если для многообразия $(M,\Omega)$ величина $a$ строго положительна, то \ref{12.6.B} влечёт \ref{12.6.A}. 
Условие $a > 0$ проверить трудно. 
Оно выполняется для многообразий Лиувилля (см. \ref{7.3.B}), в случае, когда $\pi_1(\Ham(M, \Omega)$) конечно (например, для поверхностей или $\CP^2$), и в некоторых других примерах в размерности 4. 
Пока нет известных примеров многообразий с $a = 0$. 

Гипотеза \ref{12.6.B} подтверждена для следующего интересного класса  симплектических многообразий \rindex{Лалонд}\cite{LM2}.
Опишем этот класс подробно. 
Предположим, что $M$ полумонотонно (см. \ref{11.3.A}), то есть существует константа $k \ge 0$ такая, что $c_1(A) \z= k([\Omega],A)$ для любого сферического класса гомологий $A \in H_2 (M;\ZZ)$. 
Далее, если $M$ открыто, предположим дополнительно, что $M$ «геометрически ограничено» на бесконечности. 
Например, любое кокасательное расслоение $\T^\ast N$ со стандартной симплектической структурой, геометрически ограничено.
Точное определение дано в \cite{AL}.
Тогда выполняется \ref{12.6.B}.

Вернёмся к вопросу, который мы начали обсуждать в \ref{sec:8.2}.
Пусть дана геодезическая.
Что можно сказать о длине временного интервала, на котором она минимальна? 
Существует красивый подход к этому вопросу, основанный на тщательном изучении замкнутых орбит соответствующего гамильтонова потока. 
Он становится особенно прозрачен в случае автономных геодезических, то есть однопараметрических подгрупп группы $\Ham(M, \Omega)$. 
Пусть $\{f_t\}$ — однопараметрическая подгруппа $\Ham(M, \Omega)$.
Стандартная техника обычных дифференциальных уравнений даёт, что существует строго положительное число $\tau > 0$ со следующим свойством \cite[Sec. 5.7]{HZ}). 
\rindex{замкнутая орбита}\emph{Каждая замкнутая орбита потока $\{f_t\}$ на временн\'{о}м интервале $[0; \tau]$ есть неподвижная точка потока.} 

\begin{thm}{Гипотеза}\label{12.6.D}
Сужение любой однопараметрической подгруппы на $[0;\tau]$ является
слабой кратчайшей.%
\footnote{
При некоторых дополнительных предположениях это утверждение доказал Ён Гон O \cite{O01}.
Он ипользовал гомологии Флоера, описанные в следующей главе.
Полное доказательство получено Макдафф \cite[Proposition 1.5(i)]{McD01}, где результат выведен из обобщённой теоремы о несжимаемости и использует методы, разработанные Ф. Лалондом и Д. Макдафф \cite{LM2}.
Доказательство Макдафф использует теорию $J$-голоморфных кривых.
\dpp}

\end{thm}

Это гипотеза доказана \rindex{Хофер}Хофером для случая $\RR^{2n}$ \cite{H2}. 
Её можно рассматривать как глобальную версию критерия сопряжённых точек \ref{12.3.A}.%
\footnote{На самом деле, открытие сопряжённых точек обосновано результатом Хофера.}  
Для обобщений на неавтономные потоки, а также на другие симплектические многообразия мы отсылаем к
\rindex{Бялый}\rindex{Зибург}\rindex{Шварц}\rindex{Слимвиц}\cite{BP1,
  Si1,LM2,Sch3,MSl}. 
Гипотеза \ref{12.6.D} даёт интересный способ доказательства
существования нетривиальных замкнутых орбит автономных гамильтоновых
потоков.  
Действительно, предположим, что мы заранее знаем, что автономный поток
$\{f_t\}$ \textbf{не} является минимальным на некотором временном
интервале. 
Это может следовать, например, из наличия сопряжённых точек или
процедуры укорачивания, описанных в разделах \ref{sec:8.2} и
\ref{sec:8.3}.  
Тогда \ref{12.6.D} гарантирует существование нетривиальных замкнутых
орбит на этом интервале.  
Мы отсылаем читателя к \rindex{Лалонд}\cite{LM2} и \cite{P8} для
дальнейшего ознакомления.  

Изучение кратчайших привело к пониманию следующей удивительной
особенности хоферовской геометрии.

\begin{thm}[\cite{BP1}.]{$\bm{C^1}$-уплощение}\rindex{уплощение}
\label{12.6.E}
Существуют $C^1$-окрестность $\E$ единицы в группе $\Ham(\RR^{2n})$ и $C^2$-окрестность $\C$ нуля в её алгебре Ли $\A$ такие, что пространство $(\E,\rho)$ изометрична $(\C, \|\ \|)$. 
\end{thm}

Некоторые обобщения даются в \cite{LM2}.
Поучительно сравнить явление $C^1$-уплощения с теоремой \rindex{Сикорав}Сикорава \ref{8.2.A}, утверждающей, что каждая однопараметрическая подгруппа $\Ham(\RR^{2n})$ лежит на ограниченным расстоянием от единицы. 
Она похожа на положительность кривизны, которая интуитивно противоречит уплощению! 
Пока неясно, как разрешить этот парадокс.
Мы отсылаем читателя к \cite{BP1,HZ} за доказательством $C^1$-уплощения (см. также \cite{P8} для дальнейшего обсуждения). 

Существует несколько различных подходов к представленным выше результатам о локальных кратчайших. 
Все они достаточно сложные.
В следующей главе мы представляем одну из них, которую, наверное, можно довести до доказательства гипотезы \ref{12.6.B} в полной общности. 
Мы обсудим лишь следующий простейший случай.

\begin{thm}{Теорема}\label{12.6.F}
Пусть $(M, \Omega)$ — замкнутое симплектическое многообразие с $\pi_2(M) = 0$.
Тогда любая однопараметрическая подгруппа группы $\Ham(M, \Omega)$, порождённая независящим от времени гамильтонианом общего положения на $M$ является локально слабой кратчайшей. 
\end{thm}

Набросок доказательства дан в разделе \ref{sec:13.4}.

\chapter[Гомологии Флоера]{В гостях у гомологий Флоера}\label{13}

В этой главе мы дадим набросок доказательства теоремы \ref{12.6.F}, в которой утверждается, что любая однопараметрическая подгруппа группы $\Ham(M, \Omega)$, порождённая гамильтонианом общего положения, является локально минимальной, если $\pi_{2}(M) = 0$.
Наш подход основан на теории \rindex{гомологии Флоера}гомологий Флоера.
Изложение не является
ни полным, ни стопроцентно строгим.
Его цель дать представление об этом сложном и всё ещё развивающемся наборе инструментов, а не
предоставить систематическое введение в теорию.
Мы довольно точно следуем двум статьям \rindex{Шварц}М. Шварца \cite{Sch2,
  Sch3}.

\section{У входа}\label{13.1}

Гомологии Флоера — один из мощнейших инструментов современной симплектической топологии. 
Его создание вызвал следующий вопрос:
\textit{Какие бывают инварианты гамильтоновых диффеоморфизмов?}
План для ответа приблизительно следующий.
Будем работать на связном замкнутом симплектическом многообразии
$(M,\Omega)$, для простоты предполагая асферичность многообразия:
$\pi_{2}(M) = 0$.
Роль этого предположения скоро прояснится.
Введём некоторые обозначения.
Обозначим через $\widetilde\Ham(M,\Omega)$ универсальное накрытие
группы гамильтоновых диффеоморфизмов с отмеченной точкой $\1$.
Будем обозначать через $\L M$ пространство стягиваемых петель $S^{1}\to
M$.
Через $\L\Ham(M,\Omega)$ обозначим группу стягиваемых петель
$\{h_{t}\}$ гамильтоновых 
диффеоморфизмов, 
начинающихся в $\1$
и порождённых гамильтонианами из $\H$.
Чтобы упростить обозначения, мы часто опускаем зависимость от $M$
и $\Omega$ и пишем $\Ham$ вместо $\Ham(M, \Omega)$ и т. д. 

Группа $\L\Ham$ канонически действует на $\L M$
следующим образом
\[
T_{h}: \{z(t)\}\mapsto \{h_{t}z(t)\}
\]

\parbf{Первое наблюдение:} существует естественное
отображение $\widetilde\Ham\z\to C^{\infty}(\L M)/\L\Ham$.
Опишем его.
Выберем элемент $\phi\in\widetilde\Ham$.
Обозначим
через $\F(\phi)$ множество всех гамильтонианов $F\in\F$,
порождающих $\phi$.

Группа $\L\Ham$ действует транзитивно на $\F(\phi)$. Это
действие определяется следующим образом. Рассмотрим петлю $h\in\L\Ham$ и
гамильтониан $F\in\F(\phi)$. Обозначим через $\{f_{t}\}$
гамильтонов поток функции $F$.
Тогда $h(F)$ определяется как нормализованный гамильтониан, порождающий поток
$h^{-1}_{t}\circ f_{t}$.
Из формулы \ref{1.4.D} вытекает, что
\[
h(F)(x,t) = -H(h_{t}x,t) + F(h_{t}x,t)
\]
где $H$ это гамильтониан порождающий $\{h_{t}\}$.

Для $F\in\F(\phi)$ определим функцию $A_{F}:\L M\to\RR$,
называющуюся \textit{функционалом симплектического действия}:
\[
A_{F}(z)=\int_{0}^{1}F(z(t),t)\d t - \int_{D}\Omega
\]
где $D$ это диск, затягивающий петлю $z$.
Так как $M$ асферично, результат не зависит от выбора $D$.

\begin{ex}{Упражнение}\label{13.1.A}
Докажите, что
\[
(T_{h}^{-1})^{*}A_{F}= A _{h(F)}
\]
для всех $h\in\L\Ham$ и $F\in\F(\phi)$.
\textit{Подсказка:}
Воспользуйтесь тем, что для каждой стягиваемой петли $\{h_{t}\}\in\L\Ham$, порождённой некоторым $H\in\H$, действие тождественно равно нулю на её орбитах, то есть $A_{H}(\{h_{t}x\}) \z= 0$ для любого $x\in M$.
На самом деле на асферических многообразиях это верно даже для нестягиваемых петель $\{h_{t}\}$.
Этот трудный результат недавно доказал \rindex{Шварц}Шварц [Sch3].
\end{ex}

Таким образом, мы получили естественное отображение, которое переводит
$\phi\in\widetilde\Ham$ в класс эквивалентности $[A_{F}]\in C^{\infty} (\L M)/\L\Ham$.

Функция с точностью до диффеоморфизма — очень богатый объект.
Например, в конечномерном случае можно узнать многое посмотрев на критические точки и топологию поверхностей уровня.
Мощным инструментом для получения этих данных является теория Морса.
Нам придётся работать с бесконечномерным многообразием — пространством петель $\L M$.
Ключевое наблюдение, сделанное \rindex{Флоер}Флоером, состоит в том, что существует подходящая версия теории Морса, работающая в бесконечномерных пространствах.
Эта версия теории будет описана в следующих разделах.
Теория Морса — Флоера порождает довольно сложную структуру, естественным образом связанную 
с группой $\widetilde\Ham$.
Наша основная задача — выяснить роль
хоферовской нормы гамильтонова сиплектоморфизма $\phi$ в этой структуре.
Как мы увидим, она тесно связана со значениями функционала действия
в так называемых гомологически существенных критических точках. 

Начнём с экскурса в конечномерный случай.%
\footnote{Смотри книжку \rindex{Шварц}\cite{Sch1}.}
Следующее замечание может помочь читателю развить
правильную интуицию.  Многообразие $M$ естественным образом
отождествляется с подмножеством $\L M$, состоящим из постоянных
петель.
Если функция $F\in\F$ не зависит от времени, то сужение $A_{F}$ на $M$ равно $F$.
Таким образом, обычная теория функций на $M$ «сидит» внутри теории функционалов действия на $\L M$.

\section[Конечномерный случай]{Гомологии Морса в конечномерном\\ случае}\label{13.2}
\rindex{гомологии Морса}

Пусть $F$ — функция Морса на замкнутом связном $N$-мерном
многообразии $M$. 
Множество критических точек $F$ будем обозначать через $\Crit F$.
Индекс Морса%
\footnote{То есть число отрицательных квадратов в канонической форме $d^{2}_{x}F$} критической точки $x$ будет обозначаться через $i(x)$, а через $\Crit_{m} F$ — множество критических точек с индекса~$m$.
Обозначим через $C(F)$ векторное пространство над $\ZZ_{2}$,
порождённое $\Crit F$, а через $C_{m}(F)$ его подпространство,
порождённое критическими точками индекса $m$.

Выберем риманову метрику \index{общее положение}\textbf{общего положения}%
\footnote{Здесь и далее понятие «общего положения» д\'{о}лжно разуметь
  так же как и в последнем абзаце раздела \ref{sec:4.2}: метрика
  общего положения это элемент некоторого плотного подмножества
  пространства всех метрик, которое является счётным пересечением
  открытых всюду плотных подможеств.} 
$r$ на $M$ и рассмотрим отрицательный градиентный поток
\[
\frac{\d u}{\d s} (s) = -\nabla_{r} F(u(s)).
\]
Выберем пару точек $x_{-},x_{+} \in \Crit F$.

\begin{thm}{Факт}\label{13.2.A}
Пространство орбит $u(s)$ градиентного потока, удовлетворяющих условиям $u(s)\z\to x_{-}$ при $s\z\to-\infty$ и $u(s)\z\to x_{+}$ при $s\z\to+\infty$ является гладким многообразием размерности $i(x_{-})-i(x_{+})$.
\end{thm}
  
Заметим, что это пространство допускает естественное свободное
$\RR$-действие. 
В самом деле, если $u(t)$ это решение, то и $u(t+\const)$ — тоже решение.
Таким образом, если $i(x_{-})-i(x_{+}) = 1$, то факторпространство
является нуль-мерным многообразием.
На самом деле можно показать, что \textit{оно состоит из конечного числа точек.}
Обозначим через $k_{r}(x_{-},x_{+})\in \ZZ_{2}$ чётность этого числа.
Определим линейный оператор
\[
\partial_{r}: C_{m}(F)\to C_{m-1}(F)
\]
следующим образом. Для каждого $x \in \Crit_{m}(F)$ определим
\[
\partial_{r}x = \sum_{y\in\Crit_{m-1}(F)}k_{r}(x,y)y
\]

\begin{thm}{Факт}\label{13.2.B}
  Оператор $\partial_{r}$ является дифференциалом: $\partial_{r}^{2}=0$.
  Таким образом, $(C(F),\partial_{r})$ это цепной комплекс. 
\end{thm}

\begin{thm}{Факт}\label{13.2.C}
  Группа $H_{m}(C(F),\partial_{r})$ гомологий размерности $m$ этого
  комплекса изоморфна группе гомологий $H_{m}(M; \ZZ_{2})$ многообразия.
\end{thm}

В частности, хотя эти группы и зависят от дополнительных параметров
$F$ и $r$, все они изоморфны друг-другу.
Замечательно то, что эти изоморфизмы организуются в естественное семейство следующим образом. 

Рассмотрим пространство пар $(F, r)$, где $F$ — функция, а $r$ —
риманова метрика. 
Выберем две пары $\alpha = (F_{0}, r_{0})$ и
$\beta = (F_{1},r_{1})$ в общем положении.
Выберем также путь общего положения
$(F_{s},r_{s})$, $s\in\RR$ такой, что
$(F_{s}, r_{}) = (F_{0}, r_{0})$ при $s\le0$ и
$(F_{s}, r_{s}) = (F_{1},r_{1})$ при $s\ge1$.
Рассмотрим уравнение
\begin{equation}\label{13.2.D}
\frac{\d u}{\d s}(s)=-\nabla_{r_{s}}F_{s}(u(s)).
\end{equation}

Пусть $x_{-}\in\Crit F_{0}$ и $x_{+}\in\Crit F_{1}$ — две критические точки.
Как и прежде, в общем положении пространство решений $u(s)$,
удовлетворяющих условиям $u(s)\to x_{-}$ при $s\to-\infty$ и $u(s)\to
x_{+}$ при  $s\to+\infty$, --- это гладкое многообразие размерности
$i(x_{-})-i(x_{+})$.
Далее, если $i(x_{-}) = i(x_{+})$, то существует лишь конечное число
решений.%
\footnote{Существенное различие между уравнением~\ref{13.2.D} и
  градиентным потоком состоит в том, что пространство его решений
  не допускает $\RR$-действия, если семейство $(F_{s},r_{s})$ зависит
  от $s$ нетривиальным образом.}
Обозначим через $b(x_{-},x_{+})\in\ZZ_{2}$ чётность этого числа.
Определим линейный оператор
$I^{\beta,\alpha} : C_{*}(F_{0})\to C_{*}(F_{1})$ формулой
\[
I^{\beta,\alpha}(x) = \sum_{i(y)=i(x)}b(x, y)y.
\]

\begin{thm}{Факт}\label{13.2.E}
  \begin{itemize}
  \item
    Каждый из операторов $I^{\beta,\alpha}$ является цепным
    отображением и индуцирует изоморфизм
    \[
    I^{\beta,\alpha}_{*} :
    H_{*}(C(F_{0}),\partial_{r_{0}}) \to
    H_{*}(C(F_{1}),\partial_{r_{1}})
    \]
  \item
    Oператор $I^{\beta,\alpha}_{*}$ не зависит от выбора пути
    $(F_{s},r_{s})$ общего положения.
  \item
    $I^{\alpha,\alpha}=\1$ и $I^{\gamma,\beta}_{*}\circ
    I^{\beta,\alpha}_{*}=I^{\gamma,\alpha}_{*}$, где $\alpha$, $\beta$ и
    $\gamma$ находятся в общем положении.
  \end{itemize}
\end{thm}

Назовём семейство операторов $I^{\beta,\alpha}_{*}$ удовлетворяющее
последнему свойству, естественным семейством. 

\begin{ex}{Определение}\label{13.2.F}
Пусть $(C, \partial)$ — цепной комплекс над полем $\ZZ_{2}$ с заданным
базисом $B = \{e_{1},\dots,e_{k}\}$. 
Элемент $e\in B$ называется \rindex{гомологически
  существенный}\emph{гомологически существенным}, если для 
любого $\partial$-инвариантного подпространства $K\subset
\Span(B\setminus \{e\})$ индуцированное вложением отображение 
\[
H_{*}(K,\partial)\to H_{*}(C,\partial)
\]
\textbf{не является} сюръективным.
\end{ex}

Следующее утверждение играет решающую роль в наших дальнейших рассуждениях.
Пусть $F$ — функция Морса общего положения.
Предположим, что $x_{+}\in M$ является его единственной точкой
абсолютного максимума. 
Для римановой метрики $r$ общего положения на $M$ рассмотрим комплекс
$(C(F),r)$ с базисом $\Crit F$. 

\begin{thm}{Предложение}\label{13.2.G}
  Точка $x_{+}\in\Crit F$ гомологически существенна.  
\end{thm}
\parbf{Набросок доказательства:} Поскольку функция $F$
убывает вдоль траекторий потока с отрицательным градиентом,
пространство
$Q \z= \Span(\Crit F \setminus \{x_{+}\})$ является $\partial_{r}$-инвариантным.
Более того, $H_{N}(Q,\partial_{r})$ обращается в нуль, где $N = \dim M$.
Это отражает то, что многообразие $\{F<a\}$ открыто для $a\le
\max F$ и, таким образом, у него нет фундаментального цикла. 
Мы заключаем, что для любого $\partial$-инвариантного подпространства
$Q$ образ его $N$-й группы гомологий в $H_{N}(C(F),\partial_{r})$ равен
нулю. Так как
\[
H_{N}(C(F),\partial_{r})=H_{N}(M;\ZZ_{2})=\ZZ_{2}
\]
критическая точка $x_{+}$ является гомологически существенной.
\qeds

Это предложение даёт следующее важное свойство коэффициентов $b(x,
y)$, которые считают чётность числа решений
уравнения~\ref{13.2.D}. 
Пусть $F$ — функция Морса общего положения с единственной точкой абсолютного максимума $x_{+}$.
Рассмотрим семейство $(F_{s},r_{s})$, $s\in\RR$, как в
уравнении~\ref{13.2.D}, для которого $F_{s}$ равна некоторой функции
Морса $F_{0}$ для всех $s\le0$ и $F_{s} = F$ для всех $s\ge1$.

\begin{thm}{Следствие}\label{13.2.H}
Существует такой $x\in\Crit F_{0}$, что $b(x, x_{+})\z\neq0$.  
\end{thm}

\parit{Доказательство.}
Рассмотрим оператор
\[
I:(C(F_{0}),\partial_{r_{0}})\to (C(F),\partial_{r_{1}})
\]
определённый нашими данными.
Если все $b(x,x_{+})$ равны нулю, то образ оператора $I$ содержится в
$\Span(\Crit F\setminus\{x_{+}\})$. 
Это $\partial_{r_{1}}$-ин\-ва\-ри\-ант\-ный подкомплекс.
Более того, его гомологии совпадают с $H(C(F),\partial_{r_{1}})$,
так как $I_{*}$ это изоморфизм.
Получаем противоречие с тем, что точка $x_{+}$ является гомологически существенной.
\qeds

\section{Гомологии Флоера}\label{13.3}
Существенные аспекты теории, описанной в предыдущем разделе,
обобщаются на бесконечномерный случай.
Заменим многообразие $M$ пространством $\L M$ стягиваемых в нём
петель, а функциональное пространство $C^{\infty}(M)$ — пространством
симплектических функционалов действия $A_{F}$, где $F\in\F$. 

Начнём с описания критических точек $A_{F}$.
Обозначим через $\Pcal(F)\subset\L M$ множество стягиваемых
$1$-периодических орбит гамильтонова потока, порождённого $F$.

\begin{ex}{Упражнение}\label{13.3.A}
Докажите, что критические точки $A_{F}$ это в точности элементы $\Pcal(F)$.
Более того, если отображение $f_{1}$ в момент времени $t=1$
невырождено в том смысле, что его график трансверсален диагонали в
$M\times M$, то каждая критическая точка $A_{F}$ невырождена. 
\end{ex}

Обозначим через $\J$ пространство всех почти комплексных структур
на $M$ совместимых с $\Omega$.
Выберем элемент $J\in C^{\infty}(S^{1},\J)$.
Каждый такой $J$ определяет риманову метрику на $\L M$ следующим образом.
Касательное пространство к $\L M$ в петле $z\in\L M$
состоит из
векторных полей на $M$ вдоль $z$. 
Для двух таких векторных полей, скажем $\xi(t)$ и $\eta(t)$ определим
их скалярное произведение как 
\[
\int_{S^{1}}\Omega\big(\xi(t),J(t)\eta(t)\big)\d t
\]

Для выбранного $J\in\J$ обозначим через $\nabla_{J}$ градиент
относительно римановой метрики $\Omega(\xi,J\eta)$ на $M$.
\begin{ex}{Упражнение}\label{13.3.B}
  Докажите, что градиент функционала $A_{F}$ по отношению к римановой метрике на $\mathcal{L}M$,
  описанной выше, определяется выражением 
  \[
  \grad A_{F}(u)=J\frac{\d u}{\d t}(t)+\nabla_{J(t)}F(u(t),t).
  \]
\end{ex}
Таким образом, для определения комплекса Морса функционала $A_{F}$
необходимо исследовать решения следующей задачи.
\begin{quote}
  Найти гладкое отображение $u \: \RR(s)\times S^{1}(t)\to M$ такое,
  что
  \begin{equation}\label{13.3.C}
    \frac{\partial u}{\partial s}(s,t) +
    J(t)\frac{\partial u}{\partial t}(s,t) +
    \nabla_{J(t)}F(u(s,t),t) = 0
  \end{equation}
  и удовлетворяющее условиям $u(s,t)\to z_{\pm}$ при $s\to\pm\infty$,
  где $z_{\pm}\z\in\Pcal(F)$.
\end{quote}
Обозначим через $\M_{F,J}(z_{-},z_{+})$ пространство решений
уравнения~\ref{13.3.C}. 
Флоер установил, что это пространство имеет очень красивую структуру,
которую мы сейчас опишем.
Важно заметить, что уравнение~\ref{13.3.C} --- это задача Фредгольма.
И в самом деле, видно, что с точностью до членов нулевого порядка это
хорошо знакомое нам уравнение Коши — Римана. 
Оказывается, что для $F$ и $J$ общего положения пространство
$\M_{F,J}(z_{-},z_{+})$, снабжённое естественной топологией, является
гладким многообразием.  
Более того, размерность \index[symb]{$d(z_{-}, z_{+})$}$d(z_{-}, z_{+})$ этого многообразия не меняется ни при гомотопиях гамильтонова пути $\{f_{t}\}$, $t\in [0;1]$, с фиксированными конечными точками, ни при возмущении $J$.
Таким образом, это число зависит только от поднятия отображения $f_{1}$ на универсальное накрытие группы $\Ham(M,\Omega)$. 

Мы подошли к самому удивительному во всей теории.
Конечномерная интуиция подсказывает, что число $d(z_{-},z_{+})$ это не что иное, как разность индексов Морса критических точек $z_{-}$ и $z_{+}$ функционала $A_{F}$.
Однако легко увидеть, что эти индексы бесконечны.
Тем не менее их разность имеет смысл!
Общая формула для
$d(z_{-},z_{+})$ довольно сложна — она требует понятия индекса Конли — Цендера, который мы не обсуждаем в этой книге.
Однако всё резко упрощается в случае, когда
$F\in\F(g_{\epsilon})$, где $\{g_{t}\}$ --- это некоторый гамильтонов
поток, порождённый нормализованной, независящей от времени функцией
Морса $G$, а $\epsilon>0$ достаточно мало.
Ниже мы описываем гомологии Флоера для гамильтонианов
$F\in\F(g_{\epsilon})$.

\begin{ex}{Упражнение}\label{13.3.D}
  \begin{itemize}
  \item
    Покажите, что при достаточно малом $\epsilon > 0$ каждая
    $\epsilon$-пе\-ри\-оди\-чес\-кая орбита автономного потока $\{g_{t}\}$
    является неподвижной точкой потока. Таким образом, $\Pcal(\epsilon
    G)$ совпадает с множеством $\Crit G$ критических точек функции $G$.
  \item
    Докажите, что для $F\in\F(g_{\epsilon})$ любая 1-периодическая
    орбита потока $\{f_{t}\}$ стягиваема. 
    Более того, отображение
    \[
    j : \Pcal(F)\to\Crit G,
    \quad
    \{z(t)\}\mapsto z(0)
    \]
    устанавливает взаимно однозначное соответсвие между
    1-пе\-ри\-оди\-чес\-кими орбитами потока $\{f_{t}\}$ и критическими
    точками функции $G$.
  \end{itemize}
\end{ex}


Напомним, что $i(x)$ обозначает индекс Морса точки $x\in\Crit G$.

\begin{thm}[(Формула размерности)]{Факт}\label{13.3.E}
Для  функции $F\in\F(g_{\epsilon})$ общего положения и для всех $z_{\pm}\in\Pcal(F)$ выполняется следующее равенство:
\[d(z_{-}, z_{+}) = i(z_{-}(0)) - i(z_{+}(0)).\]
\end{thm}

Заметим, что пространство решений $\M_{F}(z_{-},z_{+})$ допускает
естественное $\RR$-действие сдвигами $u(s, t)\mapsto u(s+c, t)$, где
$c\in\RR$.
Граничные условия гарантируют, что такое действие свободно.
Предположим дополнительно, что $ i(z_{-}(0)) - i(z_{+}(0)) = 1$.
Тогда из формулы размерности следует, что
факторпространство $\M_{F,J}(z_{-},z_{+})/\RR$ является нульмерным
многообразием.
Вариант теоремы Громова о компактности решений
уравнения~\ref{13.3.C} утверждает, что это многообразие компактно,
следовательно, состоит из конечного числа точек.
Обозначим через $k_{J}(z_{-},z_{+})\in\ZZ_{2}$ чётность этого числа.

Далее поступим точно так же, как и в конечномерном случае.
Напомним, что критические точки функционала $A_{F}$ отождествляются с $\Crit G$
посредством отображения $j$ (упражнение~\ref{13.3.D}).
Положим $C_{m}(A_{F}) \z= C_{m}(G)$ и $C(A_{F}) = C(G)$.
Зафиксируем петлю $\{J(t)\}$ и определим линейный оператор
$\partial_{J}:C(A_{F})\to C(A_{F})$ следующим образом.
Для элемента $x\in\Crit_{m} G$ положим
\[
\partial_{J}(x)=\sum_{y\in\Crit_{m-1}G}k_{J}(j^{-1}x,j^{-1}y)y
\]

\begin{thm}{Факт}\label{13.3.F}
Для функции $F$ и почти комплексной структуры $J$ общего положения оператор $\partial_{J}$ является дифференциалом: $\partial_{J}^{2}=0$.
Таким образом $(C_{A},\partial_{J})$ это цепной комплекс.
\end{thm}

\begin{thm}{Пример}\label{13.3.G}
Рассмотрим уравнение~\ref{13.3.C} для $F = \epsilon G$ и для постоянной петли $J(t)\equiv J$.
Каждое решение $u(s)$ уравнения градиентного потока
\[
\frac{\d u}{\d s} (s) = -\epsilon\nabla_{J} G(u(s))
\]
также является решением уравнения~\ref{13.3.C}.
Довольно сложное рассуждение \rindex{Саламон}\cite[Lemma 7.1]{HS} показывает, что при достаточно малом $\epsilon>0$ это единственные решения при условии, что $G$ и $J$ находятся в общем положении.
Поэтому $(C(A_{F}),\partial_{J})$ совпадает с обычным комплексом Морса функции $F$!
\end{thm}

Следующее рассуждение иллюстрирует важную тонкость в теории Флоера.
Попробуем доказать сформулированное выше утверждение для случая, когда
$d(z_{-}, z_{+})=1$.  
Мы должны показать, что каждое решение уравнения~\ref{13.3.C} не
зависит от времени.
Поскольку $G$ не зависит от времени, пространство решений
$\M_{\epsilon G,J}(z_{-},z_{+})$ допускает естественное $\RR\times
S^{1}$-действие $u(s, t) \mapsto (s\z+c_{1},t+c_{2})$, где
$(c_{1},c_{2})\in\RR\times S^{1}$.
Предположим, что $u(s, t)$ — решение, нетривиально зависящее от $t$.
Тогда действие группы $\RR\times S^{1}$ свободно в окрестности точки $u$,
а значит, размерность $d(z_{-}, z_{+})$ не меньше двух.
Это противоречие доказывает, что $u$ не может зависеть от $t$.
К сожалению, в этом рассуждении есть дыра.
А именно, \emph{независящие от времени} функции образуют «очень тонкое» множество в $\F$, поэтому никакие рассуждения с использованием общего положения не могут гарантировать, что пространство решений $\M_{\epsilon G,J}(z_{-}, z_{+})$ является многообразием!
Полное доказательство приводится в~\cite{HS}.

Пусть $F_{0}\in\F(g_{\delta})$ и $F_{1}\in\F(g_{\epsilon})$ две
функции общего положения. 
Рассмотрим семейство $(F_{s}, J_{s})$, $s\in\RR$, где
$F_{s}\in\F$, $J_{s}\in C^{\infty}(S^{1},\J)$, удовлетворяющее
условиям 
$(F_{s},J_s)=(F_{0},J_0)$ при $s\le0$
и $(F_{s},J_s)\z= (F_{1},J_1)$ при $s\z\ge1$.
Точно так же, как и в конечномерном случае, мы определим
естественный гомоморфизм
\[
I:(C(A_{F_{0}}),\partial_{J_{0}})\to(C(A_{F_{1}}),\partial_{J_{1}})
\]
Следующая задача похожа на уравнение~\ref{13.2.D}.

\begin{quote}
  Найти гладкое отображение $u:\RR(s)\times S^{1}(t)\to M$
  такое что 
  \begin{equation}\label{13.3.H}
    \frac{\partial u}{\partial s}(s,t)+
    J_{s}(t)\frac{\partial u}{\partial t}(s,t) +
    \nabla_{J_{s}(t)}F_{s}(u(s,t),t) = 0
  \end{equation}
  и удовлетворяющее $u(s,t)\to z_{\pm}$ при $s\to\pm\infty$, где
  $z_{-}\in\Pcal(F_{0})$ и $z_{+}\in\Pcal(F_{1})$.
\end{quote}
Анализ этого уравнения подобен тому, что мы обсудили выше.
В частности, если $i(z_{-})=i(z_{+})$, то при некоторых предположениях
общего положения пространство решений нульмерно и компактно. 
Таким образом, оно состоит из конечного числа точек.
Обозначим через $b(z_{-},z_{+})\in\ZZ_{2}$ чётность этого числа.
В отличие от уравнения~\ref{13.3.C} решения
уравнения~\ref{13.3.H} не инвариантны относительно сдвигов
$u(s,t)\mapsto u(s+c,t)$, если $(F_{s}, J_{s})$ нетривиально зависит
от переменной $s$.
Определим теперь линейный оператор $I$ следующим образом.
Для $x\in\Crit_{m}G$ положим
\begin{equation}
I(x) = \sum_{y\in\Crit_{m}G}b(j^{-1}x,j^{-1}y)\,y
\label{13.3.I}
\end{equation}
Те же рассуждения, что в~\ref{13.2.E} показывают, что $I$ индуцирует изоморфизм $I_{*}$ в гомологиях, независящий от выбора пути общего положения $(F_{s}, J_{s})$. 
Причём точно так же, как и в разделе~\ref{13.2}, семейство отображений
$I_{*}$, связанных с различным выбором параметров $(F,J)$, можно
организовать в естественное семейство.
В частности, из утверждений~\ref{13.3.G} и~\ref{13.2.C} следует,
что 
\[
H_{*}(C(A_{F}),\partial_{J}) = H_{*}(M;\ZZ_{2}).
\]
для $F$ и $J$ общего положения.

\section{Приложение к геодезическим}\label{sec:13.4}

Теперь у нас есть всё необходимое, чтобы доказать теорему~\ref{12.6.F}.
Пусть $G$ это нормализованная функция Морса на $M$ с
единственной
точкой абсолютного максимума%
\footnote{Единственность точек абсолютного
максимума и минимума можно предположить, так как гамильтониан из
теоремы~\ref{12.6.F} находится в общем положении. — \textit{Прим. ред.}} $x_{+}$ и единственной точкой абсолютного
минимума $x_{-}$. 
Обозначим через $\{g_{t}\}$ её гамильтонов поток.
Возьмём любой гамильтонов путь $\{f_{t}\}$, $t\in[0;1]$ с $f_{0}=\1$,
$f_{1}=g_{\epsilon}$ такой, что $\{f_{t}\}$ гомотопен $\{g_{\epsilon
t}\}$, $t\in[0;1]$, с неподвижными концами.
Пусть $F\in\F$ — его нормализованный гамильтониан.

\begin{thm}{Предложение}\label{13.4.A}
Имеют место следующие неравенства: для любого достаточно малого $\epsilon>0$
\begin{align*}
\int_{0}^{1}\max_{x}F(x,t) &\ge \epsilon\max G\\
\int_{0}^{1}\min_{x}F(x,t) &\le \epsilon\min G
\end{align*}
\end{thm}

Теорема~\ref{12.6.F} немедленно следует из предложения.

Для доказательства~\ref{13.4.A} начнём со следующего вспомогательного утверждения.
Рассмотрим петлю $h_{t} = f_{t} \circ g_{\epsilon t}^{-1}$.
В силу~\ref{13.1.A} $A_{F} = T_{h}^{*}A_{\epsilon G}$.
Выберем \emph{независящую от времени} почти комплексную структуру
$J_0\in\J$ общего положения на $M$ и рассмотрим петлю $J(t) \z= h_{t*}^{-1}
J_{0}h_{t*}$. 
Обозначим через $r_{0}$ и $r$ римановы метрики на $\L M$, связанные с $J_{0}$ и $J(t)$ соответственно. 
Тогда выполняется $r = T_{h}^{*}r_{0}$.
Таким образом, с помощью $T_{h}$ можно отождествить комплексы
$(C(A_{\epsilon G}), \partial_{J_{0}})$ и $(C(A_{F}), \partial_{J})$.
Кроме того, это отождествление сохраняет базис $\Crit G$ обоих
комплексов поэлементно.

\begin{thm}{Лемма}\label{13.4.B}
  Точка $x_{+}$ гомологически существенна в $(C(A_{F}), \partial_{J})$.  
\end{thm}

\parit{Доказательство}.
Утверждения~\ref{13.3.G} и~\ref{13.2.G} влекут, что
элемент $x_{+}$ гомологически существенен в $(C(A_{\epsilon G}), \partial_{J_{0}})$. 
Результат немедленно следует из приведённого выше отождествления.
\qeds

\parit{Доказательство~\ref{13.4.A}}.
Выберем $\delta\in(0,\epsilon)$.
Пусть $a(s)$, $s\in\RR$, — неубывающая функция, для которой
$a(s)\equiv0$ при $s\le0$ и $a(s)\equiv1$ при $s\ge1$.
Положим
\[
F_{s}(x, t) = (1 - a(s))\cdot\delta\cdot G(x) + a(s)F(x)
\]

Выберем такой путь $J_{s}\:S^{1}\to\J$, что $J_{s}(t)\equiv J_{0}$
при $s\le0$ и $J_{s}(t)\equiv J(t)$ при $s\ge1$
Рассмотрим гомоморфизм $I$, определённый уравнением~\ref{13.3.I}.
Положим $z_{+}(t)=f_{t}x_{+}\in\Pcal(F)$.
Пользуясь тем, что точка $x_{+}$ гомологически существенна
(см.~\ref{13.4.B}), и рассуждая в точности как в доказательстве
следствия~\ref{13.2.H}, получаем, что $b(z_{-}, z+)\z\neq 0$ для
некоторого $z_{-}\in\Pcal(\delta G)$.
Таким образом, существует решение $u(s, t)$ задачи~\ref{13.3.H}.
Рассмотрим интеграл энергии
\[
E=
\int_{-\infty}^{+\infty}\d s
\int_{0}^{1}
\Omega\left(\frac{\partial u}{\partial s}(s,t),
J_{s}(t)\frac{\partial u}{\partial s}(s,t)\right)
\,\d t
\]
Используя уравнение~\ref{13.3.H} и то, что
$\sgrad F=J\nabla_{J}F$ для каждого $J\in\J$, мы можем вычислить
\begin{equation}\label{13.4.C}
  \Omega\left(\frac{\partial u}{\partial s},
  J_{s}(t)\frac{\partial u}{\partial s}\right)
  =
  \Omega\left(\frac{\partial u}{\partial s},
         \frac{\partial u}{\partial t}\right)
  -
  \d F_{s}\left(\frac{\partial u}{\partial s}\right)
\end{equation}
Пусть $D^2_{-}$ и $D^2_{+}$ — ориентированные диски в $M$, затягивающие $z_{-}$ и
$z_{+}$ соответственно.
Поскольку $\pi_{2}(M)=0$, мы получаем, что
\begin{equation}\label{13.4.D}
  \int_{-\infty}^{+\infty}\d s
  \int_{0}^{1}
  \Omega\left(\frac{\partial u}{\partial s},
  \frac{\partial u}{\partial t}\right)
  \,\d t
  =
  \int_{D^2_{+}}\Omega-\int_{D^2_{-}}\Omega
\end{equation}
Кроме того,
\begin{equation}\label{13.4.E}
  \d F_{s}\left(\frac{\partial u}{\partial s}\right)
  =
  \frac{\d}{\d s}\big(F_{s}(u(s,t),t)\big) -
  \frac{\partial F_{s}}{\partial s}(u(s,t),t)  
\end{equation}
Интегрируя равенство~\ref{13.4.C} и подставляя~\ref{13.4.D}, \ref{13.4.E}
получаем
\[
E=A_{\delta G}(z_{-})-A_{F}(z_{+})+Q
\]
где
\[
Q=\int_{-\infty}^{+\infty}\d s
  \int_{0}^{1}\d t\;
  \frac{\d a}{\partial s}(s)\big(F(u(s,t),t)-\delta G(u(s,t))\big),
\]
Ясно, что
\begin{align*}
  Q\le
  \int_{0}^{1} \max_{x}\big(F(x,t)-\delta G(x)\big)\,\d t\\
  A_{F}(z_{+})=A_{\epsilon G}(x_{+})=\epsilon\max G
\end{align*}
и $E\ge0$.
Более того, $A_{\delta G}(z_{-})\le \delta\max G$, так как все
замкнутые орбиты потока $\{g_{\delta t}\}$ являются просто
критическими точками функции~$G$.

Следовательно,
\[
\int_{0}^{1}\max_{x}\big(F(x,t)-\delta G(x)\big)\,\d t
\ge
(\epsilon-\delta)\max G
\]
Поскольку это верно для всех  $\delta > 0$, мы получаем, что
\[
\int_{0}^{1}\max_{x} F(x,t)\,\d t
\ge
\epsilon\max G.
\]
Это завершает доказательство первого неравенства в~\ref{13.4.A}.
Второе доказывается точно так же.
\qeds

\section{К выходу}\label{13.5}

Для начала просуммируем всё сказанное в этой главе.
Пусть $(M,\Omega)$ — асферическое симплектическое многообразие.

\let\subsectionsave=\subsection
\makeatletter
\renewcommand{\subsection}{%
  \@startsection{subsection}%
  {2}%
  {0pt}%
  {1ex}%
  {0pt}%
  {\it}}
\makeatother
\def\thesubsection{\thesection.\Alph{subsection}}

\begin{ex}{}\label{13.5.A}
Каждой функции $F\in\F$ общего положения ставится в соответствие
векторное пространство $C(A_{F})$ над полем $\ZZ_{2}$.
Это пространство снабжено выделенным базисом $\Pcal(F)$, состоящим из
всех стягиваемых периодических орбит соответствующего гамильтонова
потока. 
Кроме того, $C(A_{F})$ имеет $\ZZ$-градуировку в терминах индекса
Конли — Цендера. 
Мы явно описали эту градуировку в терминах индекса Морса в простейшем
случае, когда $F$ принадлежит $\F(g_{\epsilon})$.
\end{ex}

\begin{ex}{}\label{13.5.B}
Выбор петли совместимых почти комплексных структур общего положения
$J\in C^{\infty}(S^{1}, \J)$ определяет дифференциал
$\partial_{J}\:C_{*}(A_{F})\to C_{*-1}(A_{F})$.
Гомологии комплекса $(C(A_{F}),\partial_{J})$ изоморфны
$H_{*}(M;\ZZ_{2})$. 
\end{ex}

\begin{ex}{}\label{13.5.C}
Для заданных $(F_{0},J_{0})$ и $(F_{1},J_{1})$ и соединяющего их пути
$(F_{s},J_{s})$ общего положения существует естественное с точностью
до цепной гомотопии отображение 
комплексов $I\:(C(A_{F_{0}}),\partial_{J_{0}})\z\to(C(A_{F_{1}}),\partial_{J_{1}})$.
Отображение $I$ индуцирует изоморфизм $I_{*}$ на гомологиях, независящий от выбора пути. 
Изоморфизмы $I_{*}$ образуют естественное семейство в соответствии в смысле~\ref{13.2.E}. 
\end{ex}

\begin{ex}{}\label{13.5.D}
Предположим, что $F_{0}$ и $F_{1}$ порождают один и тот же элемент в
$\widetilde\Ham(M,\Omega)$. 
Тогда $F_{0} = h(F_{1})$ для некоторой петли $h\z\in\L\Ham(M,\Omega)$
(см.~\ref{13.1.A}).
\end{ex}

Для заданного $J_{0}\in C^{\infty}(S^{1},\J)$ положим
$J_{1}(t)=h_{t*}^{-1}J_{0}(t)h_{t*}$.
Тогда отображение $T_{h}:\L M\to\L M$ введённoe в разделе~\ref{13.1} отождествляет
комплекс $(C(A_{F_{0}}),\Pcal(F_{0}),\partial_{J_{0}})$ с $(C(A_{F_{1}}),\Pcal(F_{1}),\partial_{J_{1}})$ 
Это отождествление намного сильнее, чем описанное в п.~\ref{13.5.C} — оно работает на уровне цепных комплексов, тогда как $I$ индуцируют изоморфизм только на гомологическом уровне. 

\medskip

Структура, описанная в \ref{13.5.A}—\ref{13.5.D} — это простейшая
часть так называемой теории гомологий 
Флоера, связанной с симплектическим многообразием. 
Следующее утверждение связывает эту теорию с нормой Хофера.

\begin{ex}{}\label{13.5.E}
Предположим, что для некоторого $J\in C^{\infty}(S^{1},\J)$ стягиваемая периодическая орбита $z\in\Pcal(F)$ является гомологически существенной в $(C(A_{F}),\Pcal(F),\partial_{J})$.
Тогда
\[
A_{F}(z)\le\int_{0}^{1}\max_{x}F(x,t)\,\d t
\]
\end{ex}

Оно доказывается так же как~\ref{13.4.A}.
Далее, обозначим через $\phi\z\in\widetilde\Ham(M,\Omega)$ элемент, порождённый $F$.
С учётом п.~\ref{13.5.D} гомологическая
существенность орбиты $z$ для некоторого $J$ является внутренним свойством
неподвижной точки $z(0)$ гамильтонова автоморфизма $\phi$.
Значение $A_{F}(z)$ не зависит от конкретного выбора $F\in\F(\phi)$.
Таким образом, неравенство~\ref{13.5.E} даёт способ оценки хоферовского расстояния от $\1$ до $\phi$ на универсальном накрытии $\widetilde\Ham(M,\Omega)$.
Остаётся разработать механизм, позволяющий решить, какие орбиты
гомологически существенны. 
В~\ref{13.4.B} мы рассмотрели простейший случай, когда $z$ соответствует абсолютному максимуму малого неавтономного гамильтониана. 
В~\cite{Sch3} \rindex{Шварц}Шварц использовал более изощрённое
рассуждение, позволившее доказать локальную минимальность широкого
класса геодезических для асферических симплектических многообразий.  
На самом деле существует важная дополнительная структура,
канонически связанная с $(C(A_{F}),\Pcal(F))$, а именно
\textit{каноническая вещественная фильтрация} комплекса $C(A_{F})$.
Для $a\in\RR$ обозначим через $C^{a}$ подпространство $C(A_{F})$,
порождённое теми $z\in\Pcal(F)$, которые удовлетворяют неравенству
$A_{F}(z)\le a$.
Поскольку функционал действия убывает вдоль траекторий его
отрицательного градиентного потока, это подпространство является
$\partial_{J}$-инвариантным. 
Таким образом, можно определить относительные группы гомологий
$H(C^{a}/C^{b},\partial_{J})$.
Оказывается, что эти гомологии знают о $\phi$ много интересного.
Такая фильтрация впервые систематически рассматривалась \rindex{Витербо}Витербо \cite{V1}.
Дальнейшие результаты были получены Ён Гон О~\cite{O4} и Шварцем~\cite{Sch3}.

Замечу, что существуют два естественных направления обобщения
обрисованной выше теории. 
Первый — распространить её на симплектические многообразия с
нетривиальной $\pi_{2}$.
Для таких многообразий функционалы действия $A_{F}$ являются
\emph{многозначными} функциями $\L M\to\RR$.
На самом деле их дифференциалы $\d A_{F}$ являются хорошо определёнными замкнутыми 1-формами на $\L M$.
Гомологии Флоера в этом случае устроены как обобщение гомологий
Морса — \rindex{Новиков}Новикова замкнутых 1-форм (см.~\rindex{Хофер}\cite{HS}).
Второе направление состоит в том, чтобы распространить
когомологические операции (такие как произведение в когомологиях) на
гомологии Флоера (см.~\rindex{Пиунихин}\cite{PSS}). 

Наконец, замечу, что гомологии Флоера
представляют собой очень сложную структуру, связанную с
$\widetilde\Ham(M,\Omega)$. 
Интересная задача, которая ещё далека от решения, состоит в том, чтобы
разработать алгебраический язык, пригодный для прозрачного описания
этой структуры. 
Мы ссылаемся на~\rindex{Фукая}\cite{Fu} для знакомства с первыми шагами в
этом направлении.

\chapter[Негамильтоновы диффеоморфизмы]{Геометрия негамильтоновых
  диффеоморфизмов}\label{chap:14}

В этой главе мы обсудим связь между группами $\Ham(M,\Omega)$ и
$\Symp_0(M,\Omega)$ и роль негамильтоновых симплектоморфизмов в
хоферовской геометрии. 

\section{Гомоморфизм потока}\label{sec:14.1}

Пусть $(M,\Omega)$ — замкнутое симплектическое многообразие.
Напомним, что $\Symp_0(M,\Omega)$ обознчает компоненту единицы в
группе симплектоморфизмов (см. \ref{1.4.C}).
Для заданного пути $\{f_t\}$ симплектоморфизмов с $f_0 = \1$
рассмотрим векторное поле $\xi_t$ на $M$, порождающее этот путь.
Заметим, что $\L_{\xi_t}\Omega=0$ и значит
$i_{\xi_t}\Omega=\lambda_t$ — семейство замкнутых 1-форм.
Подчеркнём, что эти формы не обязательно точны.
Напомним базовый пример, который мы начали обсуждать в начале книги
(см. \ref{1.4.C}).

\begin{thm}{Пример}\label{14.1.A}
Рассмотрим двумерный тор $(\TT^2=\RR^2/\ZZ^2$, $\d p\z\wedge \d q)$ и
следующую систему уравненинй 
\[
\begin{cases}
\quad\dot p=0,
\\
\quad\dot q=1.
\end{cases}
\]
соответствующую пути симплектоморфизмов $f_t(p, q) = (p, q \z+ t)$.
Этот путь порождён автономным семейством замкнутых 1-форм $\lambda_t = dp$.
\end{thm}

Наш первый вопрос заключается в следующем.  Пусть дан симплектоморфизм
$f\in\Symp_0(M,\Omega)$ (скажем, заданный явной формулой).
\textit{Как понять, является ли $f$ гамильтоновым?}  Мощным
инструментом, позволяющем ответить на этот вопрос в приведённом выше
примере, является понятие гомоморфизма \emph{потока};
его ввёл \rindex{Калаби}Калаби и оно изучалось далее
\rindex{Баньяга}Баньягой \cite{B1}.
Это понятие является предметом настоящей главы.

Обратим внимание, что если гамильтонов поток $\{f_{t}\}$ задаётся неточными формами, то, вообще говоря, $f_{1}$ может быть гамильтоновым.
Так, в приведённом выше примере $f_1 = \1$ является гамильтоновым.

Введём следующее полезное понятие.
Пусть $\{f_t\}$, $t\in[0;1]$ — петля симплектоморфизмов с $f_0 = f_1 = \1$.
Пусть ${\lambda_t}$ — семейство замкнутых 1-форм, порождающих эту петлю.

\begin{ex*}{Определение}
\rindex{поток петли}\emph{Поток петли} $\{f_t\}$ задаётся формулой
\index[symb]{$\flux$} 
\[\flux(\{f_t\}) = \int_0^1 [\lambda_t]\,\d t \in H^1(M;\RR).\]
\end{ex*}

Приведём более геометрическое описание.
Пусть $C$ — 1-цикл на~$M$.
Определим 2-цикл $\partial[C] = \bigcup_t f_t(C)$, являющийся образом $C$ относительно потока $f_t$. 
Заметим, что $\partial$ является линейным отображением $H_1(M;\ZZ)\to H_2(M;\ZZ)$.

\begin{ex*}{Упражнение}
Покажите, что 
\[(\flux(\{f_t\}), [C]) = ([\Omega], \partial[C])\]
для всех $[C]\in H_1(M;\ZZ)$.
\end{ex*}

В частности, поскольку правая часть сохраняется при гомотопиях, это показывает, что $\flux(\{f_t\})$ зависит только от гомотопического класса петли $\{f_t\}$ в $\pi_1 (\Symp_0(M,\Omega))$.
Таким образом, мы получаем гомоморфизм потока
\[\flux\:\pi_1(\Symp_0(M,\Omega))\to H^1(M;\RR).\]

\begin{ex*}{Definition}
Образ $\Gamma \subset H^1(M;\RR)$ гомоморфизма потока называется
\rindex{группа потока}\emph{группой потока}. 
\end{ex*}

Если петля $\gamma$ представима гамильтоновой петлёй, то $\flux(\gamma) = 0$.
Фундаментальный результат \rindex{Баньяга}О. Баньяги \cite{B1} говорит,
что верно и обратное. 
То есть если $\flux(\gamma) = 0$ то петля $\gamma$ гомотопна гамильтоновой петле.

Понятие потока полезно обобщить на произвольные (не обязательно
замкнутые) гладкие пути симплектоморфизмов. 
Для заданного пути $\{f_t\}$, порождённого семейством замкнутых 1-форм
$\lambda_t$, положим  
\[\flux(\{f_t\})=\int_0^1[\lambda_t]\,\d t\in H^1(M;\RR).\]
Возьмём симплектоморфизм $f \in \Symp_0(M,\Omega)$ и выберем любой
путь $\{f_t\}$, для которого $f_0 = \1$ и $f_1 = f$.
Ясно, что $\flux (\{f_t\})$ зависит от выбора пути, соединяющего $\1$
с $f$, но разность между потоками любых двух таких путей принадлежит
$\Gamma$. 
Таким образом, мы получаем отображение
\[\Delta\:\Symp_0(M,\Omega)\to H^1(M;\RR)/\Gamma.\]
Проверку того, что $\Delta$ является гомоморфизмом оставляем читателю.
Баньяга показал, что $\ker(\Delta) = \Ham(M,\Omega)$, и поэтому 
\[\Symp_0(M,\Omega)/\Ham(M,\Omega) = H^1(M;\RR)/\Gamma.\]

\begin{ex}{Упражнение}\label{14.1.B}
  Воспользуйтесь последним равенством вместе с теоремой \ref{1.5.A},
  чтобы доказать, что всякая нормальная подгруппа $\Symp_0(M,\Omega)$
  имеет вид $\Delta^{-1}(K)$ для некоторой подгруппы $K$ в
  $H^1(M;\RR)/\Gamma$.
\end{ex}

Вычислим группу потока для $(\TT^2,\d p\wedge \d q)$ (см. \ref{14.1.A}).
Отождествим следующие группы: 
\begin{align*}
H_1(\TT^2;\ZZ)&=\ZZ^2,
&
H_2(\TT^2;\ZZ)&=\ZZ,
&
H^1(\TT^2;\ZZ)&=\ZZ^2
\subset
\RR^2=
H^1(\TT^2;\RR).
\end{align*}
Мы утверждаем, что на этом языке $\Gamma = \ZZ^2$.
В самом деле, для каждого $\gamma\in\pi_1(\Symp_0(\TT^2))$ имеем 
\[(\flux(\gamma), a) = ([dp \wedge dq], \partial a).\]
Здесь $\partial$ — функционал, переводящий $\ZZ^2$ в $\ZZ$, значение
класса $[\d p \wedge \d q]$ на образующей $H_2 (\TT^2;\ZZ)$ равно $1$
и $a \in H_1 (\TT^2;\ZZ)$.
Таким образом, $\flux(\gamma)\in H^1(\TT^2;\ZZ)$, откуда следует, что
$\Gamma\subset \ZZ^2$.
С другой стороны, потоки полных оборотов тора задаются формулами 
\begin{align*}
\flux( (p, q) \mapsto (p , q+ t))
&=
[dp],
\\
\flux( (p, q) \mapsto (p+ t , q))
&=
-[dq],
\end{align*}
и мы заключаем, что $\ZZ^2 \subset  \Gamma$.
Поэтому $\Gamma = \ZZ^2$.
В частности, мы видим, что $f_T(p, q) = (p, q + T)$ является
гамильтоновым ровно тогда когда $T \in \ZZ$. 
Действительно, $\Delta(f_T) = T[dp] \pmod{\ZZ^2}$.

\section{Гипотеза потока}

В общем случае вычислить группу потока не просто. 
Простой вопрос, является ли подгруппа $\Gamma$ дискретной, оказывается важным по
следующей причине. 
Напомним (см. \ref{1.4.F}), что гипотеза потока%
\footnote{Точнее, гипотеза $C^\infty$-потока.
Обратите внимание, что она эквивалентна гипотезе $C^1$-потока.
Однако гипотеза $C^0$-потока является гораздо более сильным утверждением (см. обсуждение в \cite{LMP1}).}
утверждает, что подгруппа $\Ham(M,\Omega)$ является $C^\infty$-замкнутой в $\Symp_0(M,\Omega)$.

\begin{thm}{Теорема}\label{14.2.A}
Если подгруппа $\Gamma$ дискретна, то гипотеза потока верна. 
\end{thm}

\parit{Доказательство.}
Пусть $\phi_k$ — последовательность гамильтоновых диффеоморфизмов, $C^\infty$-сходящаяся к $\phi\in\Symp_0(M,\Omega)$.
Выберем путь симплектоморфизмов, соединяющий $\phi$ с единицей.
Обозначим через $\lambda \in H^1 (M;\RR)$ его поток.
Мы утверждаем, что существует последовательность путей, соединяющих $\phi_k$ с $\phi$ и таких, что $\epsilon_k=\flux \phi_k\to 0$ при $k\to +\infty$.
Предположим, что это утверждение доказано.
Тогда (см. рис. \ref{pic-13}) для каждого $k$ существует петля с потоком $\lambda + \epsilon_k\in \Gamma$.
\begin{figure}[ht!]
\centering
\includegraphics{mppics/pic-13}
\caption{}\label{pic-13}
\vskip0mm
\end{figure}
Поскольку $\Gamma$ дискретно, мы заключаем, что $\epsilon_k=0$ при всех достаточно больших $k$.

Таким образом, диффеоморфизм $\phi$ гамильтонов, так как $\Delta(\phi) \z= 0 \pmod{\Gamma}$ и теорема доказана.

Остаётся доказать наше утверждение.
Для этого достаточно найти последовательность путей $\gamma_k$, от
$\1$ до $\phi_k^{-1}\phi$, потоки которых сходятся к $0$.

\begin{figure}[ht!]
\centering
\includegraphics{mppics/pic-14}
\caption{}\label{pic-14}
\vskip0mm
\end{figure}

Заметим, что $\phi_k^{-1}\phi \to \1$ при $t\to+\infty$ в $C^{1}$-топологии, поэтому график
$\phi_k^{-1}\phi$ близок к диагонали $L$ в $(M \times M, \Omega \oplus
-\Omega)$. 
Теперь воспользуемся следующим приёмом.
Напомним, что $L$ — лагранжево подмногообразие в $(M \times M,
\Omega \oplus -\Omega)$.

\begin{thm*}{Лемма}
Пусть $(P, \omega)$ — симплектическое многообразие, а $L\subset P$
— замкнутое лагранжево подмногообразие. 
Тогда существуют окрестность $U$ многообразия $L$ в $P$ и вложение
$f\: U \to \T^\ast L$ со следующими свойствами: 
\begin{itemize}
\item $f^\ast\omega_0 = \omega$, где $\omega_0$ — стандартная симплектическая структура на $\T^\ast L$;
\item $f(x) = (x,0)$ для всех $x\in L$.
\end{itemize}
\end{thm*}
Это вариант \rindex{теорема Дарбу}\emph{теоремы Дарбу} для лагранжевых
подмногообразий \cite{MS}. 

Применим эту лемму к диагонали $L \subset M \times M$.
Она позволяет отождествить трубчатую окрестность $L$ с трубчатой
окрестностью нулевого сечения в $\T^\ast L$. 
Поскольку $\T^\ast L \simeq \T^\ast M$, мы видим, что при больших $k$
график $\phi_k^{-1}\phi$ соответствует сечению расслоения $\T^\ast M$, которое
является замкнутой $1$-формой $\alpha_k$ на
$M$.\footnote{См. упражнение на странице
  \pageref{1-form-lagrange}. — \textit{Прим. ред.}} 

Ясно, что $\alpha_k \to 0$ при $k\z\to +\infty$.
Зафиксируем достаточно большое $k$ и рассмотрим деформацию лагранжевых
подмногообразий  
\[\graph (s\alpha_k) \subset \T^\ast M\]
для $s \in [0;1]$.
Для каждого $s$ график формы $s\alpha_k$ отождествляется с графиком симплектоморфизма $M$.
Обозначим этот симплектоморфизм через $\gamma_k(s)$.
Ясно, что поток пути $\gamma_k(s)$ равен $[\alpha_k]$.
Мы получаем требуемый путь и завершаем доказательство.
\qeds

Из приведённого выше доказательства становится ясно, в чём
сложность гипотезы потока. 
Она заключается в том, что для произвольной гамильтоновой изотопии
$\{f_t\}$, график $\{f_t\}$ может выйти из трубчатой окрестности $L$
при некотором $t$. 
Он вернётся в эту окрестность, но неясно, останется ли он графиком
точной 1-формы. 
Отметим также, что верна и теорема, обратная к теореме~\ref{14.2.A}(см. \cite{LMP1}).

\begin{thm*}{Следствие}
  Если $[\Omega] \in H^2(M;\QQ)$, то гипотеза потока верна.
\end{thm*}

\parit{Доказательство.}
Мы можем переписать это предположение как $[\Omega] \z\in
\tfrac1kH^2(M;\ZZ)$ для некоторого $k\in\ZZ$.
Таким образом,
\[(\flux(\gamma), [C]) \z= ([\Omega], \partial[C]) \z\in \frac1k\ZZ\]
для каждого $\gamma \in \pi_1(\Symp_0(M,\Omega))$. 
Поэтому $\Gamma\z\subset \tfrac1k H^1(M;\ZZ)\z\subset H^1(M;\RR)$, и мы
заключаем, что $\Gamma$ дискретно. 
Остаётся применить приведённую выше теорему.
\qeds

\section[Жёсткая симплектическая топология]{Связь с жёсткой симплектической топологией}

Приведём более концептуальное доказательство гипотезы потока для $\TT^2$.
Ясно, что $C^k$-замыкание группы $\Ham(M,\Omega)$ в $\Symp_0(M,\Omega)$
является нормальной подгруппой в $\Symp_0(M,\Omega)$. 
Следовательно, учитывая \ref{14.1.B}, достаточно показать, что для любого 
$\alpha\z\in H^1(\TT^2;\RR)/\Gamma \setminus\{0\}$
\emph{существует}
$f \in \Symp_0(M,\Omega)$
с
$\Delta(f)=\alpha$
который непредставим как предел гамильтоновых диффеоморфизмов.
Не умаляя общности, можно считать, что $f$ является сдвигом $(p, q) \mapsto (p, q+T)$, где $T \notin\ZZ$.
Знаменитая гипотеза \rindex{Арнольд}Арнольда, доказанная в \rindex{Зендер}\rindex{Конли}\cite{CZ}, утверждает, что каждый $\phi \in \Ham(\TT^2)$ имеет неподвижную точку.
Таким образом, если $\phi_k \to f$, это означало бы, что $f$ также имеет неподвижную точку — противоречие.
Это рассуждение принадлежит \rindex{Эрман}М. Эрману (1983), и оно работает и для $\TT^{2n}$.
Отметим также, что оно доказывает гипотезу $C^0$-потока на торе.

Можно попытаться обобщить это рассуждение на другие симплектические многообразия.
Идея в том, чтобы вместо предельного поведения неподвижных точек рассматривать предельное поведение гомологий Флоера.
На этом пути получается следующий результат.

\begin{thm}[(\cite{LMP1})]{Теорема}\label{14.3.A}
Предположим, что первый класс Черна $c_1(\T M)$ обращается в нуль на $\pi_2(M)$.
Тогда гипотеза потока верна.
\end{thm}

И вот ещё одно приложение теории псевдоголоморфных кривых к гипотезе потока.
Оно было найдено в \rindex{Лалонд}\cite{LMP2} для 4-мерных
симплектических многообразий и позже доказано в \cite{McD00} в полной
общности. 

\begin{thm}{Теорема}\label{14.3.B}
Ранг группы потока $\Gamma$ конечен и удовлетворяет неравенству
\[\rank_\ZZ \Gamma\le b_1(M) = \dim H^1(M;\RR).\]
\end{thm}

Как следствие мы получаем, что для симплектических многообразий, у
которых первое число Бетти равно $1$, группа $\Gamma$ дискретна, и,
значит, для них верна гипотеза потока.

\section{Изометрии в хоферовской геометрии }

На группе $\Symp_0(M,\Omega)$ нет известной естественной метрики .
Однако негамильтоновы симплектоморфизмы можно включить в рамки
хоферовской геометрии следующим образом (см. \cite{LP}). 
Заметим, что группа $\Symp(M,\Omega)$ всех симплектоморфизмов
$(M,\Omega)$ действует на $\Ham(M,\Omega)$ изометриями.
Для $\phi \in \Symp(M,\Omega)$ определим 
\[T_\phi\: \Ham(M,\Omega) \to \Ham(M,\Omega),
\quad
f \mapsto \phi f\phi^{-1}.\]
Легко проверить, что $T_\phi$ корректно определено и является
изометрией для хоферовской метрики $\rho$. 

\begin{ex*}{Определение}
Изометрия $T$ называется \rindex{ограниченная
  изометрия}\emph{ограниченной}, если $\sup \rho(f, Tf)\z<\infty$, где
точная верхняя грань берётся по всем $f\z\in\Ham(M,\Omega)$. 
\end{ex*}

Если $\phi$ гамильтонов, то изометрия $T_\phi$ ограничена.
Действительно, следующее неравенство даёт оценку, независящую от $f$:
\[\rho(f,T_\phi f) = \rho(f,\phi f\phi^{-1}) = \rho(\1,\phi
f\phi^{-1}f^{-1})\le 2\rho(\1, \phi).\] 
Определим множество $BI_0 \subset \Symp_0(M,\Omega)$ как множество
всех таких $\phi \z\in \Symp_0(M,\Omega)$, что $T_\phi$ ограничено. 
Как мы видели выше, $\Ham(M,\Omega)\subset BI_0$.

\begin{ex}{Упражнение}\label{14.4.A}
Докажите, что $BI_0$ — нормальная подгруппа в $\Symp_0(M,\Omega)$.
\end{ex}

Следующая гипотеза даёт возможную характеристику гамильтоновых
диффеоморфизмов в метрических терминах. 

\begin{ex*}{Гипотеза}
$\Ham(M,\Omega) = BI_0$.
\end{ex*}

\begin{thm}[(\cite{LP})]{Теорема}\label{14.4.B}\rindex{Лалонд}
Гипотеза верна для поверхностей рода $\ge1$ и их произведений.
\end{thm}

Приведём идею доказательства для случая, когда $M=\TT^2$.
Как следует из \ref{14.4.A}, для доказательства теоремы достаточно
показать следующее.
Пусть
\[(a,b)\in H^1(\TT^2;\RR)\setminus\Gamma=\RR^2\setminus\ZZ^2,\]
тогда существует такой симплектоморфизм $\phi\in\Symp_0(\TT^2)$ с $\Delta
(\phi) \z= (a, b) \pmod \ZZ^2$, что соответствующая изометрия
$T_\phi$ неограничена. 
Не умаляя общности можно предположить, что $a= 0$, $b \in (0; 1)$ и $\phi(p, q) = (p, q + b)$.
Заметим, что для кривой $C = \{q = 0\}$ выполняется $C \cap \phi(C) = \emptyset$.
Пусть $F = F(q)$ — нормализованный гамильтониан на $\TT^2$, носитель
которого лежит в малой окрестности $C$ и такой, что $F|_C \equiv 1$. 
Обозначим через $f_t$ соответствующий гамильтонов поток и рассмотрим
поток, образованный коммутаторами $g_t = \phi^{-1}f_t^{-1}\phi f_t$. 
Этот поток порождается гамильтонианом $G(q) = F(q) - F(q + b)$.
Поскольку $G|_C \equiv 1$ {и $\pi_{1}\Ham(\TT^{2})=0$}, из \ref{7.4.A} следует, что $\rho(\1, g_t)$
стремится к бесконечности при $t \to \infty$. 
Таким образом, изометрия $T_\phi$ неограничена.

Задача описания всех изометрий пространства $(\Ham(M,\Omega), \rho)$ открыта,
она выглядит сложной даже для поверхностей. 
В связи с этим напомню следующий классический \rindex{теорема Мазура
  — Улама}результат Мазура — Улама \cite{MU} (1932): всякая
изометрия линейного нормированного пространства с фиксированной точкой
в нуле является линейным отображением. 
Было бы интересно доказать или опровергнуть нелинейный вариант этого
утверждения: всякая изометрия группы $\Ham(M,\Omega)$ с фиксированной
точкой $\1$ является изоморфизмом групп (с точностью до композиции с
инволюцией $f\mapsto f^{-1}$). 
Предположим на мгновение, что это действительно так.
Тогда теорема \rindex{Баньяга}Баньяги \ref{1.5.D} влекла бы, что
каждая такая изометрия (с точностью до упомянутой выше инволюции)
совпадает с $T_\phi$, где $\phi\: M \to M$ — либо симплектоморфизм,
либо антисимплектоморфизм (то есть $\phi^\ast\Omega = -\Omega$).
Это означало бы, что хоферовская геометрия определяет симплектическую топологию.

{

\small

\printindex
\printindex[symb]

}

{

\sloppy

\printbibliography[heading=bibintoc]

\fussy

}

\end{document}